\documentclass[11pt]{article}

\usepackage{amsmath, amssymb, algorithm, algorithmicx, mathtools, algpseudocode, graphicx, float, authblk, listings, xcolor, bm, subfigure, amsthm, bbm, enumerate}

\usepackage{natbib}

\usepackage{caption}

\usepackage[title]{appendix}

\newtheorem{theorem}{Theorem}[section]
\newtheorem{corollary}{Corollary}[section]
\newtheorem{lemma}{Lemma}[section]
\newtheorem{proposition}[theorem]{Proposition}
\newtheorem{definition}{Definition}

\theoremstyle{remark}

\newcommand{\E}{\operatorname{\mathbb{E}}}
\newcommand{\Cov}{\operatorname{\textrm{Cov}}}

\newcommand{\Var}{\textrm{Var}}

\renewcommand{\P}{\operatorname{\mathbb{P}}}

\newcommand{\vct}[1]{\bm{#1}}
\newcommand{\mtx}[1]{\bm{#1}}

\newcommand{\trace}{\operatorname{trace}}

\newcommand{\diag}{\operatorname{diag}}
\newcommand{\Diag}{\operatorname{Diag}}

\usepackage[top=1.25in,bottom=1.25in,left=1.25in,right=1.25in]{geometry}
\floatname{algorithm}{Algorithm}

\usepackage[colorlinks, linkcolor=blue, anchorcolor=blue, citecolor=blue]{hyperref}   
\graphicspath{{./}}
\usepackage{verbatim}

\begin{document}
\title{Misspecification Analysis of High-Dimensional Random Effects Models for Estimation of Signal-to-Noise Ratios}
\author[1]{\rm Xiaohan Hu}
\author[2]{\rm Xiaodong Li}
\affil[1]{Department of Mathematics, University of California, Davis}
\affil[2]{Department of Statistics, University of California, Davis}
	
\date{}
	
\maketitle
	
\begin{abstract}
Estimation of signal-to-noise ratios and residual variances in high-dimensional linear models has various important applications, including heritability estimation in bioinformatics. One widely used estimator is the Gaussian random-effects maximum likelihood estimator (MLE), based on the likelihood of the homogeneous Gaussian random-effects model in which both the regression coefficients and the noise variables are assumed to be i.i.d. Gaussian. This paper studies the behavior of this likelihood estimator under model misspecification. For isotropic random designs with independent, symmetric, sub-Gaussian entries, we establish consistency and asymptotic normality of the SNR MLE for fixed dense coefficient vectors and independent, centered, heteroscedastic finite-moment noise, allowing moderately heavy-tailed errors. We also give parallel consistency and central limit results for correlated Gaussian noise as a benchmark. The asymptotic variance depends on the limiting aspect ratio, the true SNR, and a scalar noise-square fluctuation parameter. This explicit form yields feasible plug-in confidence intervals under independent noise in two cases where the fluctuation parameter can be estimated from response fourth moments: heterogeneous Gaussian noise and homogeneous non-Gaussian noise. Numerical simulations compare likelihood-based and method-of-moments confidence intervals under heterogeneous and non-Gaussian noise, and a real-data illustration demonstrates the resulting calibrations on high-dimensional text features.
\end{abstract}

\section{Introduction}
Estimation and inference for signal-to-noise ratios (SNR) as well as residual variances in high-dimensional linear models are fundamental statistical problems with various important applications. A notable application of SNR estimation is heritability estimation \citep{falconer1961} in genome-wide association studies (GWAS), which aims to study how much phenotypic variance can be explained by genetic variation. Another important application concerns tuning-parameter selection in regularized regression such as Lasso and Ridge regression \citep{sun2012scaled,dicker2014variance,janson2017eigenprism,dicker2016maximum,dobriban2018}. A common method for estimating the SNR in modern high-dimensional applications is the Gaussian random-effects maximum likelihood estimator (MLE) \citep{yang2011gcta, Gusev2014, DeLos2015, yang2017concepts, REMforhertability, dicker2019}. Since the formulation studied below has no fixed-effect projection and uses the marginal likelihood of $\vct{y}$ under the postulated random-effects model, we refer to this likelihood-based variance-component estimator as MLE throughout. Asymptotic analysis for likelihood estimators under linear mixed-effects models is a well-studied topic in the statistical literature; see e.g. \cite{10.2307/2333854}, \cite{jiang1996reml}, \cite{rao1997variance}, and \cite{jiang}.

An interesting line of work in the literature investigates the asymptotic behavior of random-effects likelihood estimators under misspecified models, i.e., when the true model for the coefficient vector does not follow the postulated i.i.d. Gaussian model. Going back to \cite{jiang1996reml}, consistency and asymptotic normality have been established for Gaussian random-effects likelihood estimators even if the coefficient vector consists of i.i.d. but non-Gaussian components. A recent notable paper \cite{jiang2016high} shows that such estimators can be consistent and asymptotically normal even when the true model follows a sparse random-effects model. Model misspecification analysis for random-effects MLE has also been extended to the case where the coefficient vector can be a general fixed one, at the cost of assuming that the design matrix consists of i.i.d. Gaussian entries \citep{dicker2016maximum}. That analysis relies crucially on the rotational invariance of the Gaussian design matrix, and also employs some general normal approximation tools developed in \cite{dicker2017flexible}.

Beyond random-effects likelihood estimators of the SNR or residual variances, other methods have also been proposed in the literature. Examples include the method of moments \citep{haseman1972investigation, dicker2014variance}, EigenPrism \citep{janson2017eigenprism}, and Lasso and sparsity-based methods \citep{sun2012scaled, fan2012variance, bayati2013estimating}. In high-dimensional settings, unless the coefficient vector is very sparse, the empirical performance of the random-effects MLE for estimating SNR, heritability, or noise variance is often comparable to, and sometimes substantially better than, that of the above alternative methods; see the extensive simulation studies conducted in \cite{dicker2016maximum}.

The main contribution of this paper is to analyze the homogeneous Gaussian random-effects MLE when both the coefficient model and the noise model are misspecified. We treat the true coefficient vector as fixed and dense, allow the independent noise variables to be centered, heteroscedastic, and finite-moment, and relax the Gaussian-design assumption to independent symmetric sub-Gaussian design entries. Under this setting, we prove consistency and asymptotic normality for the SNR MLE. The limiting variance depends on a scalar parameter $\kappa_\varepsilon$ that captures the fluctuation of the squared noise variables; for Gaussian heteroscedastic noise, this parameter reduces to the usual variance-profile index. We also record a parallel benchmark result for correlated Gaussian noise, where the analogous variance inflation is governed by the Frobenius norm of the noise covariance. For feasible inference, we focus on independent noise and give plug-in Wald confidence intervals in two settings where $\kappa_\varepsilon$ can be estimated from fourth moments of the response: heterogeneous Gaussian noise and homogeneous non-Gaussian noise.

This paper is organized as follows. In Section \ref{sec:methods}, we introduce the high-dimensional linear-model formulation, clarify the SNR target under isotropic and general feature covariances, and define the likelihood-based MLE studied in the paper. In Section \ref{sec:theory}, we present the main consistency and asymptotic normality results for the MLE and describe plug-in inference for independent noise under heterogeneous Gaussian and homogeneous non-Gaussian calibrations. In Section \ref{sec:experiments}, we conduct systematic simulation comparisons of confidence intervals for $\gamma_0$ under heterogeneous Gaussian and homogeneous non-Gaussian noise. The simulation section also compares likelihood-based intervals with method-of-moments baselines and includes a real-data illustration using high-dimensional text features. Supplementary design-distribution checks and additional diagnostic figures for consistency and sampling normality are deferred to the appendix. The proofs of our main results are given in Section \ref{sec:proofs}, while some preliminary tools as well as proofs of important supporting lemmas are deferred to the appendix. In Section \ref{sec:discussion}, we summarize the contributions and discuss several remaining questions for future work.

\section{Methods}
\label{sec:methods}

\subsection{High-dimensional Linear Models with Heteroscedastic and Correlated Noise}
We begin by clarifying the signal-to-noise ratio that is considered in this
paper. Suppose first that the observed covariate vector before standardization
is $\vct{x}\in\mathbb{R}^p$ with mean zero and covariance
$\mtx{\Sigma}_x=\E(\vct{x}\vct{x}^{\top})$. For the linear model
$y=\vct{x}^{\top}\vct{\beta}+\varepsilon$, the natural population signal
variance is
\[
\Var(\vct{x}^{\top}\vct{\beta})=\vct{\beta}^{\top}\mtx{\Sigma}_x\vct{\beta}.
\]
Thus, if the average noise level is $\sigma_0^2$, the SNR under a general
feature covariance is
\[
\gamma_{\Sigma}
\coloneqq
\frac{\vct{\beta}^{\top}\mtx{\Sigma}_x\vct{\beta}}{\sigma_0^2}.
\]
When the design is isotropic, $\mtx{\Sigma}_x=\mtx{I}_p$, this reduces to
\[
\gamma_0
\coloneqq
\frac{\|\vct{\beta}\|^2}{\sigma_0^2}.
\]

The theoretical analysis in this paper focuses on this isotropic-design target.
Equivalently, if $\mtx{\Sigma}_x$ is known and positive definite in a
non-isotropic problem, one may whiten the covariates by setting
$\vct{z}=\mtx{\Sigma}_x^{-1/2}\vct{x}$ and
$\vct{\theta}=\mtx{\Sigma}_x^{1/2}\vct{\beta}$. Then
\[
\vct{x}^{\top}\vct{\beta}
=
\vct{z}^{\top}\vct{\theta},
\qquad
\|\vct{\theta}\|^2
=
\vct{\beta}^{\top}\mtx{\Sigma}_x\vct{\beta},
\]
so inference for the isotropic SNR
$\|\vct{\theta}\|^2/\sigma_0^2$ is inference for
$\gamma_{\Sigma}$. Throughout the main theory, we therefore work with the
whitened or already-isotropic representation
\begin{equation}
\label{eq:Fixed_effect}
\vct{y} = \mtx{Z}\vct{\beta} + \vct{\varepsilon},
\end{equation}
where $\mtx{Z}$ is an $n \times p$ isotropic design matrix with $p$ allowed to
exceed $n$, $\vct{\beta}$ is the corresponding coefficient vector, and
$\vct{y}$ is the response vector. In the non-isotropic notation above, this
$\vct{\beta}$ should be read as the whitened coefficient
$\vct{\theta}=\mtx{\Sigma}_x^{1/2}\vct{\beta}_{\mathrm{orig}}$.

For the noise vector $\vct{\varepsilon}$, the main theory assumes independent,
centered, heteroscedastic finite-moment coordinates, independent of
$\mtx{Z}$. This allows moderately heavy-tailed noise and does not require
sub-Gaussian tails. We write
$\Var(\varepsilon_i)=\sigma_i^2$ and define the average noise level
$\sigma_0^2=n^{-1}\sum_{i=1}^n\sigma_i^2$. This formulation does not require
the true noise variables to be Gaussian. Under the isotropic representation,
our target SNR is
\[
\gamma_0 \coloneqq \|\vct{\beta}\|^2/\sigma_0^2.
\]
When the original noise has a known row covariance, one may first whiten the
response and design as a preprocessing step. The independent-noise theory below
then applies to the whitened model provided the transformed design satisfies
the same isotropic random-design assumptions. Separately, we also record an
unwhitened correlated Gaussian benchmark in Theorems
\ref{thm:consistency_correlated} and \ref{thm:asymp_normality_correlated}; a
feasible confidence interval for that setting is not developed because
estimating the corresponding covariance-inflation parameter requires additional
structure. The case of unknown feature or noise covariance, where covariance or
precision matrices must be estimated before whitening, is not covered by the
theory in this paper.

\subsection{MLE Based on a Homogeneous Gaussian Random-Effects Model}
As noted above, one common SNR estimator in practice is based on the likelihood of the Gaussian random-effects model, in which the coefficient vector is modeled as $p^{-1/2}\vct{\alpha}$, where $\vct{\alpha}$ is assumed to consist of i.i.d. $\mathcal{N}(0, \sigma_\alpha^2)$ variables. In addition, the noise terms are assumed to be independent and follow the same distribution $\mathcal{N}(0, \sigma_\varepsilon^2)$. Comparing the true model and the postulated model, it is clear that $\sigma_0^2$ corresponds to $\sigma_\varepsilon^2$, $\|\vct{\beta}\|^2$ corresponds to $\sigma_\alpha^2$, and $\gamma_0 = \|\vct{\beta}\|^2/\sigma_0^2$ corresponds to $\gamma \coloneqq \sigma_\alpha^2/\sigma_\varepsilon^2$. Based on this postulated homogeneous Gaussian random-effects model, maximum likelihood estimation can be derived for the variance components $\sigma_\alpha^2$ and $\sigma_\varepsilon^2$ \citep{jiang,jiang2016high,dicker2016maximum}. Under the above Gaussian random-effects model, $\vct{y} \sim \mathcal{N}_n(\vct{0}, \mtx{\Omega})$, where
\[
\mtx{\Omega} =\mtx{\Omega}(\sigma_{\varepsilon}^2, \sigma_\alpha^2) \coloneqq
\sigma_\varepsilon^2 \mtx{I}_n + \frac{\sigma_\alpha^2}{p}\mtx{Z}\mtx{Z}^\top \coloneqq \sigma_{\varepsilon}^2 \mtx{V}_\gamma,
\]
and
\begin{equation}
\label{eq:V_gamma}
\mtx{V}_\gamma = \mtx{I}_n + \frac{\gamma}{p}\mtx{Z}\mtx{Z}^\top. 
\end{equation}
The log-likelihood function for $(\sigma_\varepsilon^2, \sigma_\alpha^2)$ is
\begin{equation}
\label{eq:loglikelihood_single}
l(\sigma_\varepsilon^2, \sigma_\alpha^2) = c -\frac{1}{2}\log \det \left( \mtx{\Omega} \right) - \frac{1}{2} \vct{y}^\top \mtx{\Omega}^{-1} \vct{y},
\end{equation}
where $c$ is a constant. By taking the partial derivatives of the log-likelihood with respect to $\sigma_\varepsilon^2$ and $\sigma_\alpha^2$ to obtain the score functions, we obtain the following likelihood equations:
\begin{equation*}
\begin{cases}
S_{\sigma_\varepsilon^2}(\sigma_\varepsilon^2, \sigma_\alpha^2)
\coloneqq \frac{1}{2}\vct{y}^\top \mtx{\Omega}^{-2} \vct{y} - \frac{1}{2}\trace\left( \mtx{\Omega}^{-1} \right)=0   
\\
S_{\sigma_\alpha^2}(\sigma_\varepsilon^2, \sigma_\alpha^2)
\coloneqq
\frac{1}{2}\vct{y}^\top \mtx{\Omega}^{-1} \frac{1}{p}\mtx{Z} \mtx{Z}^\top  \mtx{\Omega}^{-1}  \vct{y} - \frac{1}{2}\trace\left( \mtx{\Omega}^{-1}\frac{1}{p} \mtx{Z} \mtx{Z}^\top \right)=0. 
\end{cases}
\end{equation*}
Using $\frac{1}{p}\mtx{Z}\mtx{Z}^\top = \frac{1}{\gamma}(\mtx{V}_\gamma - \mtx{I}_n)$, the two score equations reduce to a single estimating equation for the SNR $\gamma = {\sigma_{\alpha}^2}/{\sigma_{\varepsilon}^2}$:

\begin{equation}
\label{eq:mle_gamma}
\Delta(\gamma) \coloneqq \vct{y}^\top \mtx{B}_\gamma \vct{y} = 0,
\end{equation}
where
\begin{equation}
\label{eq:B_gamma}
\mtx{B}_\gamma = \frac{\mtx{V}_\gamma^{-1}}{n} - \frac{\mtx{V}_\gamma^{-2} }{\trace(\mtx{V}_\gamma^{-1})}.
\end{equation}
Let $\hat{\gamma}$ be a solution to \eqref{eq:mle_gamma}, which is referred to as the (misspecified) MLE of the true SNR $\gamma_0 = \|\vct{\beta}\|^2/\sigma_0^2$.

\subsection{Misspecification Analysis of the MLE}
We aim to study the consistency and asymptotic distribution of $\hat{\gamma}$ when the Gaussian random-effects model is significantly misspecified, i.e., the true coefficient vector $\vct{\beta}$ is a general fixed one, and the true noise $\vct{\varepsilon}$ is centered, heteroscedastic, and not necessarily Gaussian. There is a trade-off between the allowed misspecification in $\vct{\beta}$ and $\vct{\varepsilon}$ and the assumption imposed on the design matrix $\mtx{Z}$. Our main results, which will be presented in the next section, assert that the consistency and asymptotic distribution of $\hat{\gamma}$ can be rigorously established when the entries in $\mtx{Z}$ are independent, symmetric, standardized sub-Gaussian random variables and the noise coordinates are independent, centered, and have uniformly bounded $4+\delta$ moments. We also record a parallel consistency and asymptotic normality result for correlated Gaussian noise as a benchmark. The symmetry assumption on the design entries is imposed for technical reasons; supplementary numerical checks in Appendix \ref{app:additional_simulation_figures} suggest some empirical robustness beyond the exact symmetric-design setting. Our misspecification analysis is conducted under the asymptotically proportional setting $n,p \rightarrow \infty$ such that $n/p \rightarrow \tau>0$, where $1/\tau$ is usually referred to as the limiting aspect ratio in the literature.

In our main independent-noise result, we will show that the asymptotic variance of $\sqrt{n}(\hat{\gamma}-\gamma_0)$ depends only on the aspect ratio $1/\tau$, the true SNR $\gamma_0$, and a noise-square fluctuation parameter $\kappa_\varepsilon$. In the correlated Gaussian benchmark, the analogous parameter is $\kappa_\Sigma$. In order to estimate the variance and thereby make inference on the true SNR $\gamma_0$, we also need to estimate the average noise level $\sigma_0^2$ and, in the independent-noise setting, $\kappa_\varepsilon$. Given the SNR estimate $\hat{\gamma}$, the postulated Gaussian random-effects model yields the estimator
\begin{equation}
\label{eq:variance_mle_modified}
\hat{\sigma}^2 = \frac{1}{n} \vct{y}^\top \mtx{V}_{\hat{\gamma}}^{-1} \vct{y}.
\end{equation}
One intuition for this estimator comes from the postulated homogeneous Gaussian random-effects model: conditional on $\mtx{Z}$,
\[
\E[\vct{y}^\top \mtx{V}_\gamma^{-1} \vct{y}\mid \mtx{Z}]
=
\trace(\mtx{V}_\gamma^{-1}\mtx{\Omega})
=
n\sigma_\varepsilon^2,
\]
when $\mtx{\Omega}=\sigma_\varepsilon^2\mtx{V}_\gamma$.

The estimation of $\kappa_\Sigma$ is in general difficult under correlated noise. Therefore, our Wald-type confidence interval is stated for independent noise. We focus on two settings in which the noise-square fluctuation parameter $\kappa_\varepsilon$ can be estimated from the response fourth moment: heterogeneous Gaussian noise and homogeneous non-Gaussian noise. We will state these two plug-in calibrations formally in the next section.

\section{Main Results}
\label{sec:theory}
Our main goal in this section is to study the consistency and asymptotic distribution of the SNR MLE $\hat{\gamma}$, which is the solution to the estimating equation \eqref{eq:mle_gamma} derived from the homogeneous Gaussian random-effects model. The results are organized as follows. We first establish consistency under independent, centered, heteroscedastic finite-moment noise, and then record the corresponding consistency result for correlated Gaussian noise. We then introduce the variance parameters and give two central limit theorems: one for the independent finite-moment setting and one for correlated Gaussian noise. Finally, we state a feasible confidence interval result only for the independent-noise setting.

\subsection{Consistency}
We begin with the baseline consistency result under independent noise. The true coefficient vector is fixed and the true noise need not be Gaussian.

\begin{theorem}
\label{thm:consistency}
Consider the linear model \eqref{eq:Fixed_effect} with the asymptotic setting $n,p \rightarrow \infty$ such that $\sqrt{n}\left| \frac{n}{p} - \tau\right| \rightarrow 0$, where $\tau>0$ is a fixed constant. Assume that the entries of the design matrix $\mtx{Z}$ are independent, symmetric, sub-Gaussian, and unit-variance random variables, and their maximum sub-Gaussian norm is uniformly upper bounded by some numerical constant $C_0$. Let $\vct{\varepsilon}$ be independent of $\mtx{Z}$ and consist of independent centered variables with $\Var(\varepsilon_i)=\sigma_i^2$. Assume that for some constants $\delta>0$ and $C_\varepsilon<\infty$,
\begin{enumerate}[(i)]
\item $\max_{i\in[n]}\E|\varepsilon_i|^{4+\delta}\leq C_\varepsilon$;
\item $\frac{1}{n}\sum_{i=1}^n \sigma_i^2 = \sigma_0^2$, where $\sigma_0^2$ is set to be fixed for all $n$;
\end{enumerate}
Let $\vct{\beta}$ be the coefficient vector with fixed two-norm $\|\vct{\beta}\|^2 >0$ for all $n$, which implies the SNR $\gamma_0\coloneqq \|\vct{\beta}\|^2/\sigma_0^2$ is fixed for all $n$. 

Under the above conditions, there is a sequence of estimates $\hat{\gamma}_n$ as solutions to \eqref{eq:mle_gamma} satisfying $\hat{\gamma}_n \stackrel{P}{\longrightarrow} \gamma_0$ as $n \rightarrow \infty$. Moreover, the corresponding sequence of noise variance estimates in \eqref{eq:variance_mle_modified} satisfies $\hat{\sigma}^2 \stackrel{P}{\longrightarrow} \sigma_0^2$.
\end{theorem}

\begin{theorem}[Consistency under correlated Gaussian noise]
\label{thm:consistency_correlated}
Consider the same asymptotic, design, and coefficient assumptions as in Theorem \ref{thm:consistency}, except that the noise vector is allowed to be correlated Gaussian:
$\vct{\varepsilon}\sim \mathcal{N}_n(\vct{0},\mtx{\Sigma}_\varepsilon)$, independent of $\mtx{Z}$. Let the diagonal entries of $\mtx{\Sigma}_\varepsilon$ be $\sigma_1^2,\ldots,\sigma_n^2$, and assume
\[
\max_{i\in[n]}\sigma_i^2=O(1),\qquad
\frac{1}{n}\sum_{i=1}^n\sigma_i^2=\sigma_0^2,\qquad
\|\mtx{\Sigma}_\varepsilon\|_F=o(n).
\]
Then there is a sequence of estimates $\hat{\gamma}_n$ as solutions to \eqref{eq:mle_gamma} satisfying $\hat{\gamma}_n\stackrel{P}{\longrightarrow}\gamma_0$. Moreover, the corresponding sequence of noise variance estimates in \eqref{eq:variance_mle_modified} satisfies $\hat{\sigma}^2\stackrel{P}{\longrightarrow}\sigma_0^2$.
\end{theorem}

\subsection{Asymptotic Normality}
We next introduce the deterministic quantities that appear in the limiting variance. Let $f_\tau$ denote the density of the Mar$\check{c}$enko-Pastur law with parameter $\tau>0$:
\begin{align*}
f_{\tau}(x) &= \frac{1}{2 \pi \tau x} \sqrt{\left( b_{+}(\tau) - x\right) \left(x - b_{-}(\tau) \right) } 1_{\{b_{-}(\tau) \leq x \leq b_{+}(\tau)\}},
\end{align*}
where $b_{\pm}(\tau) = (1 \pm \sqrt{\tau})^2$. Note that the Mar$\check{c}$enko-Pastur law also has a point mass $1 - \tau^{-1}$ at the origin when $\tau > 1$.  For any $\tau,\gamma>0$ and positive integer $k$, define
\begin{align}
\label{eq:h_MP}
h_{k}(\gamma,\tau) = &\int_{b_{-}(\tau)}^{b_{+}(\tau)} \frac{1}{(1+\gamma x)^k} f_{\tau}(x)\,\mathrm{d}x + \left( 1 - \frac{1}{\tau} \right)1_{\{\tau > 1\}}. 
\end{align}
For independent noise, the relevant scalar parameter is the fluctuation of the squared noise variables:
\begin{equation}
\label{eq:kappa_epsilon}
\kappa_{\varepsilon}
\coloneqq
\frac{1}{2n\sigma_0^4}\sum_{i=1}^n \Var(\varepsilon_i^2).
\end{equation}
In the Gaussian special case, this reduces to the usual variance-profile heterogeneity index,
\begin{equation}
\label{eq:kappa}
\kappa_\varepsilon
=\frac{1}{n\sigma_0^4}\sum_{i=1}^n\sigma_i^4.
\end{equation}
For non-Gaussian noise, however, $\kappa_\varepsilon$ also reflects the fourth-moment behavior of the noise and may be smaller than one. With these quantities defined, we obtain the following asymptotic distribution under independent finite-moment noise:

\begin{theorem}
\label{thm:asymp_normality}
In addition to the assumptions in Theorem \ref{thm:consistency}, we further assume $\| \vct{\beta} \|_\infty = o(p^{-1/4})$ and that $\kappa_\varepsilon$ in \eqref{eq:kappa_epsilon} is fixed for all $n$.
Then, with $h_k(\gamma_0, \tau)$ as in \eqref{eq:h_MP}, as $n \rightarrow \infty$,
\begin{align}
\label{eq:asym_dist}
    \sqrt{n}\left( \hat{\gamma} - \gamma_0 \right) \Longrightarrow \mathcal{N}\left(0, 2\gamma_0^2\left(\frac{1}{ h_2(\gamma_0,\tau)-h_1^2(\gamma_0,\tau) } + \kappa_{\varepsilon}  - \tau - 1 \right)\right).
\end{align}
\end{theorem}

\begin{theorem}[Correlated Gaussian noise]
\label{thm:asymp_normality_correlated}
Under the assumptions of Theorem \ref{thm:consistency_correlated}, if, in addition, $\| \vct{\beta} \|_\infty = o(p^{-1/4})$,
$\|\mtx{\Sigma}_\varepsilon\|=O(1)$, and
\begin{equation}
\label{eq:kappa_gaussian_correlated}
\kappa_{\Sigma}
\coloneqq
\frac{1}{n\sigma_0^4}\|\mtx{\Sigma}_\varepsilon\|_F^2
\end{equation}
is fixed for all $n$, then any consistent sequence of solutions $\hat\gamma$ to \eqref{eq:mle_gamma} satisfies
\begin{align}
\label{eq:asym_dist_correlated}
    \sqrt{n}\left( \hat{\gamma} - \gamma_0 \right) \Longrightarrow \mathcal{N}\left(0, 2\gamma_0^2\left(\frac{1}{ h_2(\gamma_0,\tau)-h_1^2(\gamma_0,\tau) } + \kappa_{\Sigma}  - \tau - 1 \right)\right).
\end{align}
\end{theorem}

\subsection{Feasible Inference Under Independent Noise}
The correlated Gaussian theorem above is useful as a theoretical benchmark, but estimating $\kappa_\Sigma$ is generally difficult without additional structure. We therefore state feasible confidence intervals only for independent noise. The fourth-moment plug-in estimators below are used to estimate the scalar noise-square fluctuation parameter $\kappa_\varepsilon$ in the two settings considered here.
For independent noise and fixed $i$, the fourth moment identity
\begin{align*}
\E[y_i^4] &= \sum_{j=1}^p \left(\E[z_{ij}^4] - 3\right)\beta_j^4 + 3\|\vct{\beta}\|_2^4 + 6\|\vct{\beta}\|_2^2 \sigma_i^2 + \E[\varepsilon_i^4]
\end{align*}
holds. Under the additional assumption $\|\vct{\beta}\|_\infty = o(p^{-1/4})$, we have $\sum_{j=1}^p \beta_j^4 = o(1)$, so the first term above is asymptotically negligible. This leads to feasible estimation of $\kappa_\varepsilon$ in two common cases. First, under heterogeneous Gaussian noise, $\E[\varepsilon_i^4]=3\sigma_i^4$, so we use
\begin{equation}
\label{eq:kappa_hat_gaussian}
\hat{\kappa}_{\varepsilon,\mathrm{G}}
\coloneqq
\frac{1}{3}\left\{
\frac{1}{n\hat{\sigma}^4}\sum_{i=1}^n y_i^4
-3\hat{\gamma}^2-6\hat{\gamma}
\right\}.
\end{equation}
Second, under homogeneous non-Gaussian noise, $\sigma_i^2=\sigma_0^2$ for all $i$, and we use
\begin{equation}
\label{eq:kappa_hat_homogeneous}
\hat{\kappa}_{\varepsilon,\mathrm{H}}
\coloneqq
\frac{1}{2}\left\{
\frac{1}{n\hat{\sigma}^4}\sum_{i=1}^n y_i^4
-3\hat{\gamma}^2-6\hat{\gamma}-1
\right\}
\end{equation}
to allow non-Gaussian marginal noise distributions.

The following result records these consistency properties.

\begin{proposition}
\label{thm:consistency_kappa}
Under the assumptions in Theorem \ref{thm:asymp_normality}, suppose in addition that
\[
\max_{i\in[n]}\E|\varepsilon_i|^8=O(1).
\]
Then $\hat{\kappa}_{\varepsilon,\mathrm{G}}\stackrel{P}{\longrightarrow}\kappa_{\varepsilon}$ under independent heterogeneous Gaussian noise, and $\hat{\kappa}_{\varepsilon,\mathrm{H}}\stackrel{P}{\longrightarrow}\kappa_{\varepsilon}$ under independent homogeneous non-Gaussian noise.
\end{proposition}

The next corollary records the confidence intervals implied by Theorem \ref{thm:asymp_normality} and Proposition \ref{thm:consistency_kappa}. Let
\begin{equation}
\label{eq:plugin_var}
\mathcal{V}(\gamma,\kappa,\tau)
\coloneqq
2\gamma^2\left(\frac{1}{h_2(\gamma,\tau)-h_1^2(\gamma,\tau)}+\kappa-\tau-1\right),
\end{equation}
and denote $\tau_n=n/p$. Let $\bar\kappa_n$ be $\max\{\hat{\kappa}_{\varepsilon,\mathrm{G}},1\}$ under heterogeneous Gaussian noise or $\hat{\kappa}_{\varepsilon,\mathrm{H}}$ under homogeneous non-Gaussian noise. The truncation in the heterogeneous Gaussian case reflects the population constraint $\kappa_\varepsilon\geq 1$ and does not affect consistency. For a nominal level $1-\alpha\in(0,1)$, define
\begin{equation}
\label{eq:plugin_se}
\hat{s}_n^2
\coloneqq
\frac{1}{n}\mathcal{V}(\hat{\gamma},\bar\kappa_n,\tau_n).
\end{equation}

\begin{corollary}[Plug-in confidence interval for $\gamma_0$]
\label{cor:ci_gamma}
Under the assumptions of Proposition \ref{thm:consistency_kappa}, in either the heterogeneous Gaussian case or the homogeneous non-Gaussian case,
\[
n\hat{s}_n^2 \stackrel{P}{\longrightarrow} \mathcal{V}(\gamma_0,\kappa_{\varepsilon},\tau).
\]
Consequently, if $z_{1-\alpha/2}$ denotes the $(1-\alpha/2)$ quantile of $\mathcal{N}(0,1)$, then the Wald interval
\begin{equation}
\label{eq:ci_gamma}
\mathrm{CI}_{1-\alpha}(\gamma_0)
\coloneqq
\left[\hat{\gamma}-z_{1-\alpha/2}\hat{s}_n,\ \hat{\gamma}+z_{1-\alpha/2}\hat{s}_n\right]
\end{equation}
satisfies
\[
\P\bigl(\gamma_0\in \mathrm{CI}_{1-\alpha}(\gamma_0)\bigr) \longrightarrow 1-\alpha.
\]
\end{corollary}

\subsection{Remarks}
We discuss several implications of Theorems \ref{thm:consistency}--\ref{thm:asymp_normality_correlated}, Proposition \ref{thm:consistency_kappa}, and Corollary \ref{cor:ci_gamma}.

A prominent distinction between our results and previous work \cite{jiang2016high,dicker2016maximum} concerns noise-model misspecification. In \cite{jiang2016high}, the true coefficient vector is assumed to follow a sparse random-effects model, whereas in \cite{dicker2016maximum} it is treated as fixed. However, both works assume i.i.d. Gaussian noise. In contrast, Theorems \ref{thm:consistency} and \ref{thm:asymp_normality} allow independent, centered, heteroscedastic finite-moment noise, including moderately heavy-tailed noise. The resulting variance parameter $\kappa_\varepsilon$ summarizes the fluctuation of the squared noise variables rather than only the variance profile. Theorems \ref{thm:consistency_correlated} and \ref{thm:asymp_normality_correlated} separately record the correlated Gaussian benchmark, where the asymptotic variance depends on $\kappa_\Sigma$.

Both this paper and \cite{dicker2016maximum} consider fixed coefficient vectors, but their analysis relies crucially on the assumption that the design matrix consists of i.i.d. Gaussian entries. We relax this condition to independent symmetric, non-Gaussian entries.

\begin{figure}[hbtp] 
\centering 
\includegraphics[width=1\columnwidth]{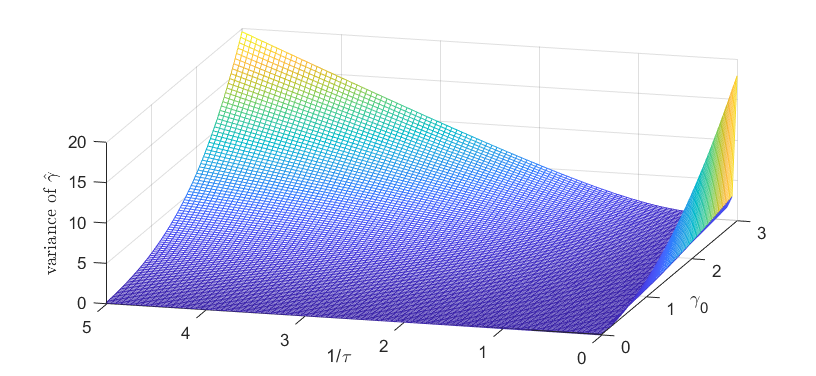} 
\caption{Asymptotic variance of $\hat{\gamma}$ with $\kappa_\varepsilon =1$ and $n=100$.}
\label{fig:mesh} 
\end{figure}

It is worth emphasizing that when the noise variables are independent, homogeneous, and Gaussian, we have $\kappa_\varepsilon =1$, and the asymptotic distribution given in \eqref{eq:asym_dist} is consistent with the result derived from i.i.d. Gaussian design in \cite{dicker2016maximum}. Setting $\kappa_\varepsilon=1$ also recovers the homogeneous Gaussian calibration underlying the confidence interval in \cite{dicker2016maximum}; our plug-in intervals can therefore be viewed as recalibrations of that likelihood-based interval for heterogeneous or non-Gaussian noise. For homogeneous but non-Gaussian noise, $\kappa_\varepsilon=\Var(\varepsilon_1^2)/(2\sigma_0^4)$ can be different from one. An explicit formula can be derived for the asymptotic variance based essentially on the Stieltjes transform of the Mar$\check{c}$enko-Pastur distribution, see e.g. Lemma 3.11 in \cite{bai2010spectral}.
Define
\begin{align*}
    m_\tau(z) & = \int_{b_{-}(\tau)}^{b_{+}(\tau)} \frac{1}{x+z} f_\tau (x) \mathrm{d}x  + \frac{1}{z}\left(1-\frac{1}{\tau}\right)1_{\{\tau>1\}}= \frac{(\tau - z - 1) + \sqrt{(\tau - z- 1)^2 + 4z\tau}}{2z \tau}.
\end{align*}
Then we can obtain
\begin{align*}
    h_1(\gamma,\tau ) = \frac{1}{\gamma} m_\tau \left(\frac{1}{\gamma}\right) = \frac{(\tau \gamma - 1 -\gamma)  + \sqrt{(\tau \gamma - 1 -\gamma)^2 + 4\tau \gamma }}{2\tau \gamma},
\end{align*}
and
\begin{align*}
    h_2(\gamma,\tau ) = -\frac{1}{\gamma^2} m'_\tau \left(\frac{1}{\gamma}\right) = -  \frac{(\tau\gamma-\tau + \gamma+1)\left(-\gamma - 1 + \sqrt{(\tau \gamma - 1 -\gamma)^2 + 4\tau \gamma }\right)}{2\gamma^2 \tau^2 \sqrt{(\tau \gamma - 1 -\gamma)^2 + 4\tau \gamma } }.
\end{align*}
We illustrate the asymptotic variance in Figure \ref{fig:mesh} with $\kappa_\varepsilon =1$ and $n = 100$. From this figure, fixing the aspect ratio $1/\tau$, the variance of $\hat{\gamma}$ increases in the true SNR $\gamma_0$; while fixing $\gamma_0$, the variance of $\hat{\gamma}$ first decreases and then increases in the aspect ratio $1/\tau$.

\section{Experiments}
\label{sec:experiments}
The purpose of this section is to compare inference procedures for the signal-to-noise ratio $\gamma_0$ under the noise misspecifications emphasized by the theory. The diagnostic simulations that visually illustrate consistency and asymptotic normality are still useful, but they are secondary to the comparison of confidence intervals. We therefore report the main simulation evidence through tables that compare Monte Carlo mean, bias, empirical coverage, and average interval length under strong independent heterogeneity and homogeneous non-Gaussian noise. Supplementary design-distribution checks and supporting figures are deferred to Appendix \ref{app:additional_simulation_figures}.

Throughout the numerical experiments, we use the Minorization-Maximization (MM) algorithm given in \cite{doi:10.1080/10618600.2018.1529601} to maximize \eqref{eq:loglikelihood_single} and obtain the random-effects MLEs $\hat{\gamma}$ and $\hat{\sigma}^2$. Unless otherwise stated, the baseline setting is $n=1200$, $p=2000$, $\gamma_0=2$, $\sigma_0^2=0.5$, the strong independent geometric heterogeneity setting described below, and $200$ independent Monte Carlo replications for coverage comparisons.

\noindent\textbf{Code availability.}
The Python code used for the simulation tables and the E2006-tfidf real-data illustration is publicly available at \url{https://github.com/xdgli/SNR-Inference}. 
The repository contains the simulation and real-data analysis scripts, together with instructions for running them. 

\subsection{Simulation Design}

We generate the coefficient vector $\vct{\beta}$ in the form
\begin{align}
\label{eq:alpha}
\vct{\beta} \propto \left( 1, 2^{-g},3^{-g},\cdots, p^{-g} \right)^\top,
\end{align}
where $g\geq 0$ controls the decay of the coefficients. Since the focus of this paper is the dense signal regime, we use a very small decay parameter, $g=0.1$, in the main comparison tables. The vector is rescaled so that $\|\vct{\beta}\|_2^2=\gamma_0\sigma_0^2$. Sensitivity checks over the SNR level $\gamma_0$ are summarized in Table \ref{tab:ci_gamma_sensitivity}, while diagnostic plots are deferred to Appendix \ref{app:additional_simulation_figures}.

For the heterogeneity-focused confidence-interval simulations, we use the Gaussian special case of the independent finite-moment noise model. Conditional on the variances, $\varepsilon_i\sim\mathcal{N}(0,\sigma_i^2)$ independently, so the variance-profile index below coincides with $\kappa_\varepsilon$. We generate the noise variances through a geometric sequence:
\begin{align}
\label{eq:sigma_i}
(\tilde\sigma_1^2,\tilde\sigma_2^2,\ldots,\tilde\sigma_n^2) \propto \left( 1, q , q^2,\cdots, q^{n-1} \right),
\end{align}
where $q>0$ and $\sum_{i=1}^n\tilde\sigma_i^2=n\sigma_0^2$. We then randomly shuffle $(\tilde\sigma_1^2,\ldots,\tilde\sigma_n^2)$ to obtain $(\sigma_1^2,\ldots,\sigma_n^2)$. The reported heterogeneous-noise comparison tables use a strong independent geometric heterogeneity setting: at the baseline dimension $n=1200$, this sets $q=0.95$, which gives the heterogeneity index $\kappa_0=30.77$. Homogeneous Gaussian noise corresponds to $\sigma_i^2=\sigma_0^2$ and $\kappa_0=1$, while the milder geometric setting $q=0.995$ gives $\kappa_0=2.54$ at the baseline dimension. In the dimension-scaling comparison, $q$ is chosen separately for each $n$ so that the same heterogeneity index $\kappa_0=30.77$ is maintained. Table \ref{tab:ci_gamma_hom_nongaussian} instead uses homogeneous rare-shock Gaussian scale-mixture noise to examine the non-Gaussian plug-in calibration under high-kurtosis errors. Correlated-noise simulations are excluded from the confidence-interval comparison because the current theory does not provide a feasible estimator of the corresponding variance-inflation parameter.

For the design matrix, the main comparisons use Gaussian or Rademacher designs as specified in each table. Additional checks for standardized $t_7$ and standardized genotype designs are reported in Appendix \ref{app:additional_simulation_figures}.

\subsection{Methods and Evaluation Metrics}

The comparison focuses on whether the confidence intervals correctly account for the additional uncertainty caused by heterogeneous noise or non-Gaussian noise. Table \ref{tab:method_grid} lists the interval constructions used across the simulation tables, with the relevant rows selected according to the noise setting. The likelihood-based rows use the same MLE point estimate $\hat{\gamma}$ and differ only in the calibration of the standard error. The row labeled MLE (homo-Gauss) fixes $\kappa_\varepsilon=1$, corresponding to the homogeneous Gaussian likelihood calibration used for the confidence interval in \cite{dicker2016maximum}; it therefore serves as the classical likelihood-based benchmark. The rows labeled MLE (hetero-Gauss) and MLE (homo-non-Gauss) use the corresponding plug-in estimates under heterogeneous Gaussian noise and homogeneous non-Gaussian noise, respectively. We also include method-of-moments calibrations as external baselines.

For the method-of-moments baseline, we use the unit-response version of the fixed-effects moment estimator in \cite{dicker2014variance,hu2025multivariate}, adapted to the present scaling. Let
\[
K=\frac{1}{p}\mtx{Z}\mtx{Z}^\top,\qquad
m_1=\frac{1}{n}\vct{y}^\top\vct{y},\qquad
m_2=\frac{1}{n}\vct{y}^\top K\vct{y},
\]
and define the empirical trace quantities
\[
c_1=\frac{1}{n}\trace(K),\qquad c_2=\frac{1}{n}\trace(K^2).
\]
The empirical-trace MoM point estimate used in the simulations is
\[
\hat{\eta}_{\mathrm{MoM}}=\frac{m_2-c_1m_1}{c_2-c_1},
\qquad
\hat{\sigma}^2_{\mathrm{MoM}}=m_1-\hat{\eta}_{\mathrm{MoM}},
\qquad
\hat{\gamma}_{\mathrm{MoM}}=\frac{\hat{\eta}_{\mathrm{MoM}}}{\hat{\sigma}^2_{\mathrm{MoM}}}.
\]
All MoM rows use this same point estimate and differ only in the covariance calibration used to form a Wald interval: MoM (homo-Gauss) uses the homogeneous Gaussian moment covariance, MoM (hetero-Gauss) uses the Gaussian moment covariance with the variance-profile shape used in the simulation, and MoM (homo-non-Gauss) adds a homogeneous fourth-cumulant correction estimated from the response fourth moment. The explicit formulas for these three intervals are collected in Appendix \ref{app:mom_derivation}. We use these moment methods only as comparators for the single-response SNR target; under the nonstandard noise settings below, their coverage is evaluated empirically rather than supported by the main theory of this paper.

\begin{table}[t!]
\centering
\small
\begin{tabular}{p{0.32\textwidth}p{0.62\textwidth}}
\hline
Method & Interval construction \\
\hline
MLE (homo-Gauss) & Plug-in interval from \eqref{eq:ci_gamma} with $\kappa_\varepsilon$ fixed at $1$, matching the homogeneous Gaussian likelihood calibration in \cite{dicker2016maximum} \\
MLE (hetero-Gauss) & Plug-in interval from \eqref{eq:ci_gamma} with $\kappa_\varepsilon$ estimated by $\max\{\hat{\kappa}_{\varepsilon,\mathrm{G}},1\}$ \\
MLE (homo-non-Gauss) & Plug-in interval from \eqref{eq:ci_gamma} with $\kappa_\varepsilon$ estimated by $\hat{\kappa}_{\varepsilon,\mathrm{H}}$ under homogeneous non-Gaussian noise \\
MoM (homo-Gauss) & MoM interval with homogeneous Gaussian moment covariance \\
MoM (hetero-Gauss) & MoM interval with heterogeneous Gaussian moment covariance \\
MoM (homo-non-Gauss) & MoM interval with homogeneous fourth-cumulant correction \\
\hline
\end{tabular}
\caption{Interval constructions compared in the simulation tables for $\gamma_0$.}
\label{tab:method_grid}
\end{table}

For an estimator $\tilde{\gamma}$ and a nominal $95\%$ confidence interval $[\mathrm{L},\mathrm{U}]$, we report
\[
\text{Mean}=\frac{1}{R}\sum_{r=1}^R\tilde{\gamma}^{(r)},\quad
\text{Bias}=\frac{1}{R}\sum_{r=1}^R(\tilde{\gamma}^{(r)}-\gamma_0),
\]
together with empirical coverage
$R^{-1}\sum_{r=1}^R\mathbbm{1}\{\mathrm{L}^{(r)}\leq \gamma_0\leq \mathrm{U}^{(r)}\}$
and average interval length
$R^{-1}\sum_{r=1}^R(\mathrm{U}^{(r)}-\mathrm{L}^{(r)})$.

\subsection{Confidence Intervals of \texorpdfstring{$\gamma_0$}{gamma0}}

The confidence-interval comparison below uses the strong independent geometric heterogeneity setting as the default. In this Gaussian-noise special case, MLE intervals calibrated with $\kappa_\varepsilon=1$, equivalently the homogeneous Gaussian calibration of \cite{dicker2016maximum}, are expected to be too narrow, while the intervals calibrated with $\hat{\kappa}_{\varepsilon,\mathrm{G}}$ account for heterogeneity by inflating the standard error according to \eqref{eq:asym_dist}.

Table \ref{tab:ci_gamma_scaling} checks the finite-sample behavior of the four displayed interval constructions across several high-dimensional regimes with $n<p$ and feature dimensions at the baseline scale or above. We use the strong independent geometric heterogeneity setting and the Gaussian design. In all rows of this dimension-scaling comparison, the geometric parameter is chosen separately for each $n$ so that the heterogeneity index remains fixed at $\kappa_0=30.77$. To keep the scaling comparison compact, the table reports only the Monte Carlo mean, bias, empirical coverage, and average interval length.

\begin{table}[t!]
\centering
\small
\resizebox{\textwidth}{!}{%
\begin{tabular}{llcccc}
\hline
$(n,p)$ & Method & Mean & Bias & Coverage & Length \\
\hline
\multicolumn{6}{l}{Settings with $n<p$} \\
$(1200,2000)$ & MLE (homo-Gauss) & 2.164 & 0.164 & 0.820 & 1.648 \\
$(1200,2000)$ & MLE (hetero-Gauss) & 2.164 & 0.164 & 0.915 & 2.453 \\
$(1200,2000)$ & MoM (homo-Gauss) & 2.209 & 0.209 & 0.860 & 2.923 \\
$(1200,2000)$ & MoM (hetero-Gauss) & 2.209 & 0.209 & 0.925 & 3.520 \\
\hline
$(1600,2000)$ & MLE (homo-Gauss) & 2.107 & 0.107 & 0.760 & 1.126 \\
$(1600,2000)$ & MLE (hetero-Gauss) & 2.107 & 0.107 & 0.940 & 1.880 \\
$(1600,2000)$ & MoM (homo-Gauss) & 2.170 & 0.170 & 0.915 & 2.170 \\
$(1600,2000)$ & MoM (hetero-Gauss) & 2.170 & 0.170 & 0.945 & 2.728 \\
\hline
$(2400,4000)$ & MLE (homo-Gauss) & 2.098 & 0.098 & 0.815 & 1.117 \\
$(2400,4000)$ & MLE (hetero-Gauss) & 2.098 & 0.098 & 0.950 & 1.696 \\
$(2400,4000)$ & MoM (homo-Gauss) & 2.146 & 0.146 & 0.920 & 1.878 \\
$(2400,4000)$ & MoM (hetero-Gauss) & 2.146 & 0.146 & 0.950 & 2.302 \\
\hline
\end{tabular}
}
\caption{Dimension-scaling comparison of $95\%$ confidence intervals for $\gamma_0$ under strong independent geometric heterogeneity, retaining only settings with $n<p$. In all rows, the geometric noise parameter is chosen separately for each $n$ so that $\kappa_0=30.77$, with $\gamma_0=2$, $\sigma_0^2=0.5$, $g=0.1$, and a Gaussian design. Each row uses $200$ Monte Carlo replications.}
\label{tab:ci_gamma_scaling}
\end{table}

Supplementary checks under additional non-Gaussian design distributions are reported in Appendix \ref{app:additional_simulation_figures}.

Table \ref{tab:ci_gamma_sensitivity} summarizes robustness to the SNR level under a Rademacher design. These rows keep $n=2400$, $p=4000$, the dense coefficient setting $g=0.1$, and the strong independent geometric heterogeneity setting fixed, and compare the low-SNR case $\gamma_0=0.5$ with the high-SNR case $\gamma_0=5$.

\begin{table}[t!]
\centering
\small
\resizebox{\textwidth}{!}{%
\begin{tabular}{llcccc}
\hline
SNR & Method & Mean & Bias & Coverage & Length \\
\hline
$\gamma_0=0.5$ & MLE (homo-Gauss) & 0.524 & 0.024 & 0.915 & 0.352 \\
$\gamma_0=0.5$ & MLE (hetero-Gauss) & 0.524 & 0.024 & 0.945 & 0.473 \\
$\gamma_0=0.5$ & MoM (homo-Gauss) & 0.516 & 0.016 & 0.900 & 0.394 \\
$\gamma_0=0.5$ & MoM (hetero-Gauss) & 0.516 & 0.016 & 0.955 & 0.507 \\
\hline
$\gamma_0=5$ & MLE (homo-Gauss) & 5.261 & 0.261 & 0.880 & 3.320 \\
$\gamma_0=5$ & MLE (hetero-Gauss) & 5.261 & 0.261 & 0.965 & 4.612 \\
$\gamma_0=5$ & MoM (homo-Gauss) & 5.670 & 0.670 & 0.910 & 9.317 \\
$\gamma_0=5$ & MoM (hetero-Gauss) & 5.670 & 0.670 & 0.940 & 10.005 \\
\hline
\end{tabular}
}
\caption{SNR-robustness comparison of $95\%$ confidence intervals for $\gamma_0$ under a Rademacher design. Except for the SNR level shown in the first column, all rows use $n=2400$, $p=4000$, $\sigma_0^2=0.5$, $g=0.1$, and strong independent geometric heterogeneity with $\kappa_0=30.77$. Each row uses $200$ Monte Carlo replications.}
\label{tab:ci_gamma_sensitivity}
\end{table}

Table \ref{tab:ci_gamma_hom_nongaussian} complements Table \ref{tab:ci_gamma_sensitivity} by replacing heterogeneous Gaussian noise with homogeneous high-kurtosis non-Gaussian noise. Specifically, we generate
\[
\varepsilon_i=\sigma_0 u_i,\qquad
u_i=\sqrt{\frac{V_i}{1.38}}\,\xi_i,\qquad
\xi_i\sim\mathcal{N}(0,1),
\]
independently across $i$, where $V_i=1$ with probability $0.98$ and $V_i=20$ with probability $0.02$. The normalizing constant $1.38$ is $\E V_i$, so $\E u_i^2=1$. This rare-shock mixture has $\E u_i^4\approx 14.15$ and $\kappa_\varepsilon=(\E u_i^4-1)/2\approx 6.57$, giving a high-kurtosis non-Gaussian noise example. The likelihood-based plug-in row uses $\hat{\kappa}_{\varepsilon,\mathrm{H}}$ from \eqref{eq:kappa_hat_homogeneous}; the MoM rows compare the homogeneous Gaussian covariance calibration with the homogeneous non-Gaussian fourth-moment plug-in calibration.

\begin{table}[t!]
\centering
\small
\resizebox{\textwidth}{!}{%
\begin{tabular}{llcccc}
\hline
SNR & Method & Mean & Bias & Coverage & Length \\
\hline
$\gamma_0=0.5$ & MLE (homo-Gauss) & 0.498 & -0.002 & 0.935 & 0.341 \\
$\gamma_0=0.5$ & MLE (homo-non-Gauss) & 0.498 & -0.002 & 0.955 & 0.364 \\
$\gamma_0=0.5$ & MoM (homo-Gauss) & 0.493 & -0.007 & 0.950 & 0.380 \\
$\gamma_0=0.5$ & MoM (homo-non-Gauss) & 0.493 & -0.007 & 0.960 & 0.401 \\
\hline
$\gamma_0=5$ & MLE (homo-Gauss) & 5.115 & 0.115 & 0.930 & 3.178 \\
$\gamma_0=5$ & MLE (homo-non-Gauss) & 5.115 & 0.115 & 0.950 & 3.427 \\
$\gamma_0=5$ & MoM (homo-Gauss) & 5.479 & 0.479 & 0.910 & 9.019 \\
$\gamma_0=5$ & MoM (homo-non-Gauss) & 5.479 & 0.479 & 0.910 & 9.185 \\
\hline
\end{tabular}
}
\caption{SNR-robustness comparison of $95\%$ confidence intervals for $\gamma_0$ under homogeneous high-kurtosis non-Gaussian noise. All rows use $n=2400$, $p=4000$, $\sigma_0^2=0.5$, $g=0.1$, a Rademacher design, and homogeneous rare-shock Gaussian scale-mixture noise with raw variance components $1$ and $20$. Each row uses $200$ Monte Carlo replications.}
\label{tab:ci_gamma_hom_nongaussian}
\end{table}

For applications where the target parameter is heritability, the SNR estimates can be transformed through
\[
h=\frac{\gamma}{\gamma+1},\qquad
\hat{h}=\frac{\hat{\gamma}}{\hat{\gamma}+1}.
\]
We keep the main tables on the $\gamma_0$ scale, since this is the primary theoretical target in the paper. The main conclusions are therefore based on Tables \ref{tab:ci_gamma_scaling}--\ref{tab:ci_gamma_sensitivity} and Table \ref{tab:ci_gamma_hom_nongaussian}.
Across these comparisons, both adjusted MLE and adjusted MoM calibrations improve coverage relative to the homogeneous Gaussian calibration when the noise is strongly heterogeneous or non-Gaussian. The likelihood-based intervals are generally shorter and, especially at high SNR, the MLE point estimates are less biased than the MoM estimates, so the MoM rows are best viewed as external benchmarks for the likelihood-based inference developed here.

\subsection{Real-data Illustration}
\label{sec:real_data_illustration}

As a real-data illustration, we analyze the public E2006-tfidf regression data \citep{Kogan2009}, downloaded from the LIBSVM data repository at \url{https://www.csie.ntu.edu.tw/~cjlin/libsvmtools/datasets/regression.html}. The full training set contains $16087$ observations and $150360$ sparse text features extracted from 10-K filings. The response is a volatility measure. For a computationally manageable illustration, we take a fixed-seed subsample of $n=3000$ observations, keep the $p=4000$ features with the largest sample variances, standardize the response and selected columns, and apply right-side shrinkage whitening to the selected design with shrinkage $0.20$. The standardized response has empirical kurtosis $4.88$, which makes this example informative for comparing Gaussian and non-Gaussian calibrations. This whitening step is a preprocessing device motivated by the general-covariance discussion in Section \ref{sec:methods}; the intervals below should therefore be read as descriptive real-data uncertainty summaries, since they do not account for covariance-estimation error in the whitening step.

Because the sparse random-effects calibration of \citet{jiang2016high} depends on the unknown sparsity proportion $m/p$, we do not apply that adjustment in this real-data example. The homogeneous Gaussian row serves as the classical likelihood benchmark related to \citet{dicker2016maximum}, while the heterogeneous and non-Gaussian rows illustrate the plug-in calibrations developed here.

\begin{table}[t!]
\centering
\normalsize
\renewcommand{\arraystretch}{1.18}
\setlength{\tabcolsep}{8pt}
\begin{tabular*}{0.98\textwidth}{@{\extracolsep{\fill}}lcccc@{}}
\hline
Method & Estimate & SE & 95\% CI & Length \\
\hline
MLE (homo-Gauss) & 3.489 & 0.368 & $[2.767,4.211]$ & 1.444 \\
MLE (hetero-Gauss) & 3.489 & 0.514 & $[2.482,4.496]$ & 2.014 \\
MLE (homo-non-Gauss) & 3.489 & 0.573 & $[2.366,4.612]$ & 2.246 \\
MoM (homo-Gauss) & 3.145 & 1.071 & $[1.047,5.243]$ & 4.197 \\
MoM (hetero-Gauss) & 3.145 & 1.112 & $[0.965,5.324]$ & 4.359 \\
MoM (homo-non-Gauss) & 3.145 & 1.208 & $[0.776,5.513]$ & 4.737 \\
\hline
\end{tabular*}
\caption{E2006-tfidf real-data SNR analysis for standardized volatility. All rows use $n=3000$, $p=4000$, $\tau=0.75$, and the shrinkage-whitened text-feature design. The homo-Gauss rows fix $\kappa_\varepsilon=1$. The hetero-Gauss MLE uses the response fourth-moment plug-in estimate $\hat{\kappa}_{\varepsilon,\mathrm{G}}=16.823$, while the hetero-Gauss MoM row uses a response-second-moment variance-profile plug-in with profile index $10.232$. The homo-non-Gauss MLE uses $\hat{\kappa}_{\varepsilon,\mathrm{H}}=24.734$, and the homo-non-Gauss MoM row uses the corresponding moment plug-in estimate $17.182$.}
\label{tab:e2006_tfidf_real_data}
\end{table}

\section{Proof of the Main Results}
\label{sec:proofs}
 \subsection{Supporting Lemmas}
\begin{lemma}
\label{lem:noise_var_bound}
Under the assumptions of Theorem \ref{thm:consistency},
we have
\begin{align}
\label{eq:max_esp_i}
\max_{i \in [n]} \varepsilon_i^2 = O_P\left(n^{2/(4+\delta)}\right)=o_P\left(\sqrt n\right)
\end{align}
and
\begin{align}
\label{eq:noise_var_concent}
\left|\frac{1}{n}\sum_{i=1}^n   \varepsilon_i^2  - \sigma_0^2 \right| = o_P(1).
\end{align}
Moreover, under the assumptions of Theorem \ref{thm:asymp_normality}, there holds
\begin{equation}
\label{eq:noise_CLT}
\frac{1}{\sqrt{n}} \left(\left(\sum_{i=1}^n \varepsilon_i^2\right) - n\sigma_0^2\right) \Longrightarrow \mathcal{N}(0, 2\kappa_\varepsilon \sigma_0^4).
\end{equation}
Under the assumptions of Theorem \ref{thm:consistency_correlated}, $\max_{i\in[n]}\varepsilon_i^2=O_P(\log n)$ and \eqref{eq:noise_var_concent} also hold. If the additional assumptions of Theorem \ref{thm:asymp_normality_correlated} hold, then \eqref{eq:noise_CLT} holds with $\kappa_\varepsilon$ replaced by $\kappa_\Sigma$.
\end{lemma}

\begin{lemma}[Trace approximation for resolvents]
\label{thm:diff_tr_h}
Under the assumptions of Theorem \ref{thm:consistency}, for $\mtx{V}_{\gamma}$ defined in \eqref{eq:V_gamma} and any fixed integer $k>0$, we have $\|\mtx{V}_{\gamma}^{-k}\|\leq 1$. Moreover, with $\tau_n=n/p$,
\[
\frac{1}{n}\trace(\mtx{V}_{\gamma}^{-k})
-
h_k(\gamma,\tau_n)
=
O_P(n^{-1}),
\]
provided the standard linear spectral statistic rate for analytic test functions is invoked; see, for example, Theorem 9.10 of \cite{bai2010spectral}. Consequently, under $\sqrt n|\tau_n-\tau|\to0$,
\[
\frac{1}{n}\trace(\mtx{V}_{\gamma}^{-k})
-
h_k(\gamma,\tau)
=
o_P(n^{-1/2}).
\]
In particular,
\[
\frac{1}{n}\trace(\mtx{V}_{\gamma}^{-k})
\stackrel{P}{\longrightarrow}
h_k(\gamma,\tau).
\]
\end{lemma}

A key technique in proving Theorem \ref{thm:consistency} is the ``leave-$k$-out'' argument developed in \cite{jiang2016high}. Here we list some useful notations.
\begin{definition}
Denote $\mtx{Z} = [\vct{z}_1, \ldots, \vct{z}_p]$ as a concatenation of column vectors. For any subset $ C \subset \{1, \ldots, p\}$, denote $\mtx{V}_{\gamma,-C} \coloneqq \mtx{V}_{\gamma} - \frac{\gamma}{p} \sum_{k \in C} \vct{z}_k\vct{z}_k^\top$. For example, for any $i \neq j$,
\[
\mtx{V}_{\gamma,-ij} \coloneqq \mtx{V}_{\gamma,-\{ij\}} = \mtx{V}_{\gamma} - \frac{\gamma}{p} \left(\vct{z}_i\vct{z}_i^\top + \vct{z}_j\vct{z}_j^\top\right).
\]
Furthermore, for $1 \leq i,j \leq p$, define
\begin{align}
\label{def:eta}
\eta_{ij,C}^{(l)} \coloneqq \vct{z}_i^\top \mtx{V}_{\gamma,-C}^{-l} \vct{z}_j.
\end{align}
Finally, in the case $C = \emptyset$, simply denote
\begin{align}
\label{def:eta_noC}
\eta_{ij}^{(l)} \coloneqq \vct{z}_i^\top \mtx{V}_{\gamma}^{-l} \vct{z}_j.
\end{align}
\end{definition}

The proofs of Lemmas \ref{lem:leave_one_out}--\ref{lem:leave_two_out} closely follow the leave-out arguments in \cite{jiang2016high}. For completeness, we provide self-contained proofs in the appendix and include the modifications needed for the present heteroscedastic and correlated-noise setting.

\begin{lemma}
\label{lem:leave_one_out}
Under the conditions of Theorem \ref{thm:consistency}, we have 
\begin{align}\label{eq:diff_Vinv_Vinvk}
\max_{k \in [p]} \left| \trace\left( \mtx{V}_{\gamma}^{-l}\right) - \trace\left( \mtx{V}_{\gamma,-k}^{-l}\right)  \right|\leq 2^l-1, \quad l=1,2,3,4,
\end{align}
and for $\eta_{kk,k}^{(l)}$ defined in \eqref{def:eta},
\[
\max_{k \in [p]} \left|\frac{1}{n}  \eta_{kk,k}^{(l)} - \frac{1}{n}\trace\left(  \mtx{V}_{\gamma}^{-l} \right) \right| = O_P\left( \sqrt{\frac{\log n}{n}}\right), \quad l = 1,2.
\]
\end{lemma}

\begin{lemma}
\label{lem:diag_homo}
Under the conditions of Theorem \ref{thm:consistency}, for fixed $\gamma >0$, we have 

\begin{align}
\label{eq:diag_homo_1}
\max_{1\leq k \leq p} \left|\frac{1}{n}\vct{z}_k^\top \mtx{V}_\gamma^{-1} \vct{z}_k - \frac{1}{n}\frac{\trace(\mtx{V}_{\gamma}^{-1})}{1 +  \frac{\gamma}{p} \trace(V_{\gamma}^{-1})} \right|
= O_P\left(\sqrt{\frac{\log n}{n}}  \right), 
\end{align}

\begin{align}
\label{eq:diag_homo_2}
\max_{1\leq k \leq p} \left|\frac{1}{n}\vct{z}_k^\top \mtx{V}_\gamma^{-2}\vct{z}_k - \frac{1}{n}\frac{\trace(\mtx{V}_{\gamma}^{-2})}{\left(1 +  \frac{\gamma}{p} \trace(V_{\gamma}^{-1}) \right)^2}\right|
= O_P\left(\sqrt{\frac{\log n}{n}}\right), 
\end{align}

\begin{align}
\label{eq:diag_homo_3}
\max_{1\leq k \leq p} \left|\frac{1}{n}\vct{z}_k^\top \mtx{V}_\gamma^{-l} \vct{z}_k - \frac{1}{np} \trace\left(\mtx{V}_{\gamma}^{-l} \mtx{Z} \mtx{Z}^\top\right) \right|
= O_P\left(\sqrt{\frac{\log n}{n}}  \right), \quad l = 1,2,
\end{align}

\begin{align}
\label{eq:diag_homo_4}
\max_{1 \leq k \leq p}\left| \left(\vct{z}_k^\top \mtx{B}_{\gamma} \vct{z}_k \right)^l -  \left( \frac{1}{p} \trace\left( \mtx{B}_{\gamma} \mtx{Z} \mtx{Z}^\top \right) \right)^l \right| 
= O_P\left(\sqrt{\frac{\log n}{n}}\right), \quad l = 1,2,
\end{align}

\begin{align}
\label{eq:diag_homo_5}
\left|\frac{1}{np} \trace\left(\mtx{V}_{\gamma}^{-1} \mtx{Z} \mtx{Z}^\top\right) - \frac{1}{n}\frac{\trace(\mtx{V}_{\gamma}^{-1})}{1 +  \frac{\gamma}{p} \trace(V_{\gamma}^{-1})} \right|
= O_P\left( \frac{1}{n}  \right), 
\end{align}

and
\begin{align}
\label{eq:diag_homo_6}
\left|\frac{1}{np} \trace\left(\mtx{V}_{\gamma}^{-2} \mtx{Z} \mtx{Z}^\top\right) - \frac{1}{n}\frac{\trace(\mtx{V}_{\gamma}^{-2})}{\left(1 +  \frac{\gamma}{p} \trace(V_{\gamma}^{-1}) \right)^2}\right|
= O_P\left(\frac{1}{n}\right). 
\end{align}

\end{lemma}

\begin{lemma}
\label{lem:diag_homo_taylor}
Under the conditions of Theorem \ref{thm:consistency}, for fixed $\gamma >0$, we have 

\begin{align} \label{eq:diag_homo_taylor_1}
&\max_{1\leq k \leq p} \Bigg\vert\frac{1}{n}\vct{z}_k^\top \mtx{V}_\gamma^{-1} \vct{z}_k - \frac{1}{np} \trace\left(\mtx{V}_{\gamma}^{-1} \mtx{Z} \mtx{Z}^\top\right) \notag
\\
&~~~~~ - \frac{1}{\left(1 +  \frac{\gamma}{p} \trace(V_{\gamma}^{-1}) \right)^2}\left(\frac{1}{n}\eta_{kk,k}^{(1)} - \frac{1}{n}\trace(V_{\gamma}^{-1})  \right)\Bigg\vert
= O_P\left(\frac{\log n}{n} \right), 
\end{align}

and
\begin{align*}
&\max_{1\leq k \leq p} \Bigg\vert\frac{1}{n}\vct{z}_k^\top \mtx{V}_\gamma^{-2}\vct{z}_k - \frac{1}{np} \trace\left(\mtx{V}_{\gamma}^{-2} \mtx{Z} \mtx{Z}^\top\right)
\\
&~~~~~~~~+\frac{\trace(\mtx{V}_{\gamma}^{-2})}{\left(1 +  \frac{\gamma}{p} \trace(V_{\gamma}^{-1}) \right)^3} \frac{2\gamma}{p}\left(\frac{1}{n}\eta_{kk,k}^{(1)} - \frac{1}{n}\trace(V_{\gamma}^{-1})  \right)
\\
&~~~~~~~~-\frac{1}{\left(1 +  \frac{\gamma}{p} \trace(V_{\gamma}^{-1}) \right)^2} \left(\frac{1}{n}\eta_{kk,k}^{(2)} - \frac{1}{n}\trace(V_{\gamma}^{-2})  \right)\Bigg\vert
= O_P\left(\frac{\log n}{n}\right). 
\end{align*}

\end{lemma}

\begin{lemma}
\label{lem:diag_homo_moment}
Under the conditions of Theorem \ref{thm:consistency}, for fixed $\gamma >0$ and $l = 1,2$, we have 
\begin{align*}
\max_{1\leq k \leq p} \E \left[\left(\frac{1}{n}\eta_{kk,k}^{(l)} - \frac{1}{n}\trace(V_{\gamma,-k}^{-l})  \right)^2\right] 
\leq \frac{C}{n},
\end{align*}
and
\begin{align*}
&\max_{1 \leq i < j \leq p}  \Bigg\vert \E \Bigg[\left(\frac{1}{n}\eta_{ii,i}^{(l)} - \frac{1}{n}\trace(V_{\gamma,-i}^{-l})  \right)\left(\frac{1}{n}\eta_{jj,j}^{(l)} - \frac{1}{n}\trace(V_{\gamma,-j}^{-l})  \right) \Bigg]\Bigg\vert \leq \frac{C}{n^{3/2}},
\end{align*}
where $C$ is a constant independent of $n$.
\end{lemma}

\begin{lemma}
\label{lem:leave_two_out}
Under the conditions of Theorem \ref{thm:consistency}, for fixed $\gamma >0$ and $l = 1,2$, we have 

\begin{align}
\label{eq:off_diag_rough}
\max_{k \neq j} \vert \vct{z}_k^\top \mtx{V}_{\gamma}^{-l} \vct{z}_j|^2 = O_P (n \log n)
\quad \text{and} \quad
\max_{k \neq j} \vert \vct{z}_k^\top \mtx{B}_\gamma \vct{z}_j|^2 = O_P \left(\frac{\log n}{n} \right).
\end{align}
Further, under the assumptions of Theorem \ref{thm:asymp_normality} or Theorem \ref{thm:asymp_normality_correlated}, we have
\begin{align}
\label{eq:off_diag_fine}
\begin{dcases}
\frac{1}{p(p-1)}\sum_{i\neq j} n\left(\vct{z}_i^\top \mtx{B}_{\gamma} \vct{z}_j \right)^2  &= \bar{\theta}_1(\gamma, \tau) + o_P(1)
\\
\sum_{i\neq j} \beta_i^2 \beta_j^2 n\left(\vct{z}_i^\top \mtx{B}_{\gamma} \vct{z}_j \right)^2  &= \|\vct{\beta}\|^4 \bar{\theta}_1(\gamma, \tau) + o_P(1),
\end{dcases}
\end{align}
where $\bar{\theta}_1(\gamma, \tau)>0$ is a constant only depending on $\gamma$ and $\tau$.
\end{lemma}

The above lemmas rely crucially on the ``leave-$k$-column-out'' argument in \cite{jiang2016high}. Because the present analysis allows heteroscedastic and correlated noise, we also need the following results, which rely on a similar ``leave-$k$-row-out'' argument. 

\begin{lemma}
\label{lem:V_inv_diag}
For any fixed $\gamma > 0$, under the assumptions of Theorem \ref{thm:consistency}, we have
\begin{equation}
\label{eq:diff_Vii_tr}
\max_{i \in [n]} \left|   \left(\mtx{V}_{\gamma}^{-l} \right)_{ii}  - \frac{1}{n} \trace\left(  \mtx{V}_{\gamma}^{-l} \right)  \right| = O_P\left(\sqrt{\frac{\log n}{n}} \right), \quad l = 1,2,3,4, 
\end{equation}
which implies
\begin{equation}
\label{eq:diff_Bii_tr}
\max_{i \in [n]} \left| (\mtx{B}_{\gamma})_{ii} - \frac{1}{n}\trace\left(\mtx{B}_{\gamma} \right) \right| =  O_P\left( \sqrt{\frac{\log n}{n^3}}\right)
\end{equation} 
and
\begin{equation}
\label{eq:diff_Bii_tr_sq}
\max_{i \in [n]} \left| (\mtx{B}_{\gamma})_{ii}^2 - \left(\frac{1}{n}\trace\left(\mtx{B}_{\gamma} \right)\right)^2 \right| = O_P\left( \sqrt{\frac{\log n}{n^5}}\right).
\end{equation} 
\end{lemma}

\begin{lemma} 
\label{lem:B_off_diag}
For any fixed $\gamma > 0$, under the assumptions of Theorem \ref{thm:consistency}, we have
\begin{equation}
\label{eq:Bij_max}
\max_{1\leq i < j \leq n}  \left\vert (\mtx{B}_{\gamma})_{ij} \right\vert = O_P\left(\sqrt{\frac{\log n}{n^3}} \right)
\end{equation} 
and under the assumptions of Theorem \ref{thm:asymp_normality} or Theorem \ref{thm:asymp_normality_correlated}

\begin{equation}
\label{eq:Bij_weight_sum}
\begin{dcases}
n\sum_{i\neq j} \left(\mtx{B}_{\gamma}\right)_{ij}^2  =  \bar{\theta}_2(\gamma, \tau) + o_P(1)
\\
n\sum_{i\neq j} \varepsilon_i^2 \varepsilon_j^2 \left(\mtx{B}_{\gamma}\right)_{ij}^2  =  \bar{\theta}_2(\gamma, \tau) \sigma_0^4 + o_P(1),
\end{dcases}
\end{equation} 
where $\bar{\theta}_2(\gamma, \tau)>0$ is a constant only depending on $\gamma$ and $\tau$.
\end{lemma}

\begin{lemma}
\label{lem:V_inv_diag_loo}
For any fixed $\gamma > 0$, under the assumptions of Theorem \ref{thm:consistency}, we have
\begin{equation}
\label{eq:diff_Vii_tr_loo}
\max_{k \in [p]}\max_{i \in [n]} \left|   \left(\mtx{V}_{\gamma, -k}^{-l} \right)_{ii}  - \frac{1}{n} \trace\left(  \mtx{V}_{\gamma}^{-l} \right)  \right| = O_P\left(\sqrt{\frac{\log n}{n}} \right), \quad l = 1,2,3,4.
\end{equation}
\end{lemma}

\subsection{New Representation based on Rademacher Sequences}
Since the entries of $\mtx{Z}$ are independent and symmetric, we can replace the original design matrix $\mtx{Z}$ with $\tilde{\tilde{\mtx{Z}}} = \mtx{\Lambda}_\zeta \mtx{Z} \mtx{\Lambda}_\xi$ with the diagonal matrices
\[
\mtx{\Lambda}_\zeta = \diag(\zeta_1,\ldots,\zeta_n),  \quad \mtx{\Lambda}_\xi =  \diag(\xi_1,\ldots,\xi_p),
\]
where $\zeta_i$'s and $\xi_j$'s are i.i.d. Rademacher random variables that are also independent of $\mtx{Z}$, since $\mtx{Z}$ and $\tilde{\tilde{\mtx{Z}}}$ have the same distribution. We also denote
\[
\vct{\xi} = (\xi_1,\ldots,\xi_p)^\top\quad \text{and} \quad \vct{\zeta} = (\zeta_1,\ldots,\zeta_n)^\top.
\]

Under this new representation of the design matrix, the linear model \eqref{eq:Fixed_effect} becomes
\begin{equation}\label{eq:fixed_effects_rad_Zeps}
    \vct{y} = \mtx{\Lambda}_\zeta \mtx{Z} \mtx{\Lambda}_\xi \vct{\beta} + \vct{\varepsilon}.
\end{equation}
We want to emphasize that under this new representation, we still define $\mtx{V}_{\gamma}$ and $\mtx{B}_{\gamma}$ as before:
\[
\mtx{V}_\gamma = \mtx{I}_n + \frac{\gamma}{p}\mtx{Z}\mtx{Z}^\top, \quad \text{and} \quad \mtx{B}_\gamma = \frac{\mtx{V}_\gamma^{-1}}{n} - \frac{   \mtx{V}_\gamma^{-2} }{\trace(\mtx{V}_\gamma^{-1})}.
\]
However, the representation of the estimating equation \eqref{eq:mle_gamma} changes: the original $\mtx{Z}\mtx{Z}^\top$ is replaced with $\mtx{\Lambda}_\zeta \mtx{Z}\mtx{Z}^\top \mtx{\Lambda}_\zeta$. Therefore, the original $\mtx{V}_\gamma$ defined in \eqref{eq:V_gamma} should be replaced with 
\[
\tilde{\tilde{\mtx{V}}}_\gamma 
=  \mtx{I}_n  + \frac{\gamma}{p} \mtx{\Lambda}_\zeta \mtx{Z}\mtx{Z}^\top \mtx{\Lambda}_\zeta 
= \mtx{\Lambda}_\zeta \left( \mtx{I}_n + \frac{\gamma}{p} \mtx{Z}\mtx{Z}^\top \right) \mtx{\Lambda}_\zeta 
=\mtx{\Lambda}_\zeta \mtx{V}_\gamma \mtx{\Lambda}_\zeta.
\]
Also, it is easy to see that the original $\mtx{B}_\gamma$ should be replaced with
$\tilde{\tilde{\mtx{B}}}_\gamma =  \mtx{\Lambda}_\zeta \mtx{B}_\gamma \mtx{\Lambda}_\zeta$. Therefore, the estimating equation \eqref{eq:mle_gamma} should be rewritten as
\begin{align}\label{eq:fixed_effects_rad_double_delta}
    \Delta(\gamma) & \coloneqq \vct{y}^\top \tilde{\tilde{\mtx{B}}}_\gamma \vct{y} \notag \\
                & =\left( \mtx{\Lambda}_\zeta \mtx{Z} \mtx{\Lambda}_\xi \vct{\beta} + \vct{\varepsilon}  \right)^\top    \mtx{\Lambda}_\zeta \mtx{B}_\gamma \mtx{\Lambda}_\zeta  \left( \mtx{\Lambda}_\zeta \mtx{Z} \mtx{\Lambda}_\xi \vct{\beta} + \vct{\varepsilon}  \right) \notag  \\
                & = \vct{\xi}^\top \mtx{\Lambda}_\beta \mtx{Z}^\top \mtx{B}_\gamma  \mtx{Z} \mtx{\Lambda}_\beta \vct{\xi} + 2 \vct{\xi}^\top \mtx{\Lambda}_\beta \mtx{Z}^\top \mtx{B}_\gamma \mtx{\Lambda}_\varepsilon \vct{\zeta} + \vct{\zeta}^\top \mtx{\Lambda}_\varepsilon \mtx{B}_\gamma \mtx{\Lambda}_\varepsilon \vct{\zeta} \notag
                \\
                & = [\vct{\xi}^\top, \vct{\zeta}^\top]
                \begin{bmatrix}
                \mtx{\Lambda}_\beta \mtx{Z}^\top \mtx{B}_\gamma \mtx{Z} \mtx{\Lambda}_\beta 
                & \mtx{\Lambda}_\beta \mtx{Z}^\top \mtx{B}_\gamma \mtx{\Lambda}_\varepsilon
                \\
                \mtx{\Lambda}_\varepsilon \mtx{B}_\gamma \mtx{Z} \mtx{\Lambda}_\beta
                & \mtx{\Lambda}_\varepsilon \mtx{B}_\gamma \mtx{\Lambda}_\varepsilon
                \end{bmatrix}
                \begin{bmatrix} \vct{\xi} \\ \vct{\zeta} \end{bmatrix},
\end{align}
where
\[
\mtx{\Lambda}_\beta= \diag(\beta_1, \ldots,\beta_p) \quad \text{and} \quad \mtx{\Lambda}_\varepsilon = \diag(\varepsilon_1, \ldots,\varepsilon_n). 
\]
Note that now $\Delta(\gamma)$ is a random variable about $\mtx{Z}$, $\vct{\varepsilon}$, $\vct{\xi}$ and $\vct{\zeta}$. Straightforward calculation gives the conditional mean of $\Delta(\gamma)$ on $\mtx{Z}$ and $\vct{\varepsilon}$:
\begin{align}
\label{eq:tilde_delta_star}
    \widetilde{\Delta}_*(\gamma) 
    &\coloneqq \E \left[\Delta(\gamma)\big\vert \mtx{Z},\vct{\varepsilon}\right]
    = \sum_{k=1}^p \beta_k^2 \vct{z}_k^\top \mtx{B}_{\gamma}\vct{z}_k + \trace\left( \mtx{\Lambda}_\varepsilon^2 \mtx{B}_{\gamma}\right).
\end{align}
Furthermore, the conditional variance of $\sqrt{n}(\Delta(\gamma))$ on $\mtx{Z}$ and $\vct{\varepsilon}$ can also be derived as in the following lemma, the proof of which is deferred to the appendix.
\begin{lemma}
\label{lem:conditional_variance_formula}
The conditional variance of $\Delta(\gamma)$ given $\mtx{Z}$ and $\vct{\varepsilon}$ has the formula
\begin{align}
\label{eq:var_gamma0}
    &~~~\Var\left[\sqrt{n}(\Delta(\gamma) ) \vert \mtx{Z},\vct{\varepsilon} \right] \notag \\
    &=\underbrace{ 2n \sum_{1 \leq k\neq j \leq p} \beta_k^2 \beta_j^2 \left(\vct{z}_k^\top \mtx{B}_{\gamma} \vct{z}_j \right)^2}_{V_1} 
    + \underbrace{ 4 n \sum_{k=1}^p \beta_k^2 \vct{z}_k^\top \mtx{B}_{\gamma}\mtx{\Lambda}_{\varepsilon}^2\mtx{B}_{\gamma} \vct{z}_k}_{V_2} 
    + \underbrace{ 2n \sum_{1 \leq k \neq j \leq n} \varepsilon_k^2 \varepsilon_j^2 (\mtx{B}_{\gamma})_{kj}^2 }_{V_3}.
\end{align}
\end{lemma}

\subsection{Proof of Theorem \ref{thm:consistency}} 
\label{sec:proof}

With $\widetilde{\Delta}_*(\gamma)$ defined in \eqref{eq:tilde_delta_star} and for any fixed $\gamma>0$, we first aim at showing 
\begin{equation}
\label{eq:cond_mean_concentration}
\Delta(\gamma) - \widetilde{\Delta}_*(\gamma) \xrightarrow[n \rightarrow \infty]{P} 0.
\end{equation}
First, by \eqref{eq:off_diag_rough} in Lemma \ref{lem:leave_two_out}, we have
\begin{align*}
\sum_{k \neq j} \beta_k^2 \beta_j^2 \left(\vct{z}_k^\top \mtx{B}_{\gamma} \vct{z}_j\right)^2 
\leq \left(\max_{k \neq j} \vert \vct{z}_k^\top \mtx{B}_{\gamma} \vct{z}_j|^2\right) \|\vct{\beta}\|_2^4 = O_P\left(\frac{\log n}{n}\right).
\end{align*}
Second, since $\vct{z}_k$'s are sub-Gaussian vectors, it is obvious that
\[
\max_{1 \leq k \leq p}  \|\vct{z}_k\|^2 = O_P(n). 
\]
Also, a simple consequence of Lemma \ref{thm:diff_tr_h} gives $\|\mtx{B}_\gamma\| = O_P(1/n)$, and Lemma \ref{lem:noise_var_bound} implies $\|\mtx{\Lambda}_{\varepsilon}^2\|=O_P(n^{2/(4+\delta)})$. Therefore,
\[
\sum_{k=1}^p\beta_k^2\vct{z}_k^\top\mtx{B}_{\gamma}\mtx{\Lambda}_\varepsilon^2 \mtx{B}_{\gamma} \vct{z}_k \leq \|\vct{\beta}\|_2^2 \|\mtx{B}_\gamma\|^2 \|\mtx{\Lambda}_{\varepsilon}^2\| \left(\max_{1 \leq k \leq p}  \|\vct{z}_k\|^2\right) = o_P(1).
\]
Third, by \eqref{eq:max_esp_i} in Lemma \ref{lem:noise_var_bound} and \eqref{eq:Bij_max} in Lemma \ref{lem:B_off_diag}, 
\[
\sum_{1\leq k \neq j \leq n} \varepsilon_k^2 \varepsilon_j^2 (\mtx{B}_{\gamma})_{kj}^2 = o_P(1).
\]
Plug the above bounds to \eqref{eq:var_gamma0}, for any $\delta > 0$, by the conditional Chebyshev's inequality, we have
\[
\P \left\{ \left| \Delta(\gamma) - \widetilde{\Delta}_*(\gamma)  \right| > \delta \Big\vert Z, \vct{\varepsilon} \right\} \leq \frac{\Var\left[\Delta(\gamma)  \vert \mtx{Z}, \vct{\varepsilon} \right] }{\delta^2} \xrightarrow[n \rightarrow \infty]{P} 0.
\]
Then, by the dominated convergence theorem, we have proved \eqref{eq:cond_mean_concentration}.

Now, define
\begin{align}
\label{eq:Delta_starstar}
    \Delta_{**}(\gamma)=&\sigma_0^2 \trace\left(\mtx{B}_\gamma \mtx{V}_{\gamma_0}\right)= \sigma_0^2 \trace\left(\mtx{B}_\gamma \left(\mtx{I}_n + \frac{\gamma_0}{p}\mtx{Z}\mtx{Z}^\top\right)\right).
\end{align}
By \eqref{eq:diag_homo_4} in Lemma \ref{lem:diag_homo}, we can easily obtain
\begin{equation}
\label{eq:Delta_starstar_part1}
\left|\sum_{k=1}^p \beta_k^2 \left(\vct{z}_k^\top \mtx{B}_\gamma \vct{z}_k\right) - \frac{\|\vct{\beta}\|^2}{p} \trace(\mtx{B}_{\gamma} \mtx{Z} \mtx{Z}^\top)  \right| = O_P\left(\sqrt{\frac{\log n}{n}}  \right).
\end{equation}
On the other hand, by Lemma \ref{lem:noise_var_bound} and \eqref{eq:diff_Bii_tr} in Lemma \ref{lem:V_inv_diag}, we have
\begin{equation*}
\left| \trace\left( \mtx{\Lambda}_\varepsilon^2 \mtx{B}_{\gamma}\right) - \frac{1}{n} \trace\left( \mtx{\Lambda}_\varepsilon^2 \right) \trace\left( \mtx{B}_{\gamma}\right) \right| = o_P(1).
\end{equation*}
Furthermore, by Lemmas \ref{lem:noise_var_bound} and \ref{thm:diff_tr_h}, we have
\begin{equation*}
\left| \frac{1}{n} \trace\left( \mtx{\Lambda}_\varepsilon^2 \right) \trace\left( \mtx{B}_{\gamma}\right) - \sigma_0^2 \trace\left( \mtx{B}_{\gamma}\right) \right| = o_P(1).
\end{equation*}
Combine the above two inequalities,
\begin{equation}
\label{eq:Delta_starstar_part2}
\left| \trace\left( \mtx{\Lambda}_\varepsilon^2 \mtx{B}_{\gamma}\right) - \sigma_0^2 \trace\left( \mtx{B}_{\gamma}\right) \right| = o_P(1).
\end{equation}
Then, by \eqref{eq:tilde_delta_star}, \eqref{eq:Delta_starstar}, \eqref{eq:Delta_starstar_part1}, and \eqref{eq:Delta_starstar_part2}, we have
\begin{align*}
\left\vert \widetilde{\Delta}_*(\gamma) - \Delta_{**}(\gamma) \right\vert = o_P(1).
\end{align*}
Combined with \eqref{eq:cond_mean_concentration}, we have
\[
\Delta(\gamma) - \Delta_{**}(\gamma) \xrightarrow[n \rightarrow \infty]{P} 0.
\]
Finally, we have the following result that characterizes the limit of $\Delta_{**}(\gamma)$ for any $\gamma>0$. 
\begin{lemma}[\cite{jiang2016high}]
\label{lem:Delta_starstar}
Under the assumption of Theorem \ref{thm:consistency}, we have 
\[
\Delta_{**}(\gamma)\stackrel{a.s.}{\longrightarrow}c_\gamma,
\]
where $c_{\gamma}>0$ for $\gamma < \gamma_0$, $c_{\gamma_0}=0$, and $c_{\gamma}<0$ for $\gamma>\gamma_0$.
\end{lemma}

This result follows from \cite{jiang2016high}, and we give a detailed proof in the appendix for completeness.

Then, for any $\gamma > 0$, there holds $\Delta(\gamma) \stackrel{P}{\longrightarrow} c_{\gamma}$, which is positive, zero, or negative, depending on whether $\gamma$ is smaller than, equal to, or greater than $\gamma_0$. To make the root-selection argument explicit, fix any $\delta>0$. Since $c_{\gamma}$ is continuous and changes sign only at $\gamma_0$, there exists $\eta_\delta>0$ such that
\[
c_{\gamma_0-\delta}\ge 2\eta_\delta
\qquad\text{and}\qquad
c_{\gamma_0+\delta}\le -2\eta_\delta.
\]
By the pointwise convergence already proved,
\[
\Delta(\gamma_0-\delta)\stackrel{P}{\longrightarrow} c_{\gamma_0-\delta},
\qquad
\Delta(\gamma_0+\delta)\stackrel{P}{\longrightarrow} c_{\gamma_0+\delta},
\]
hence
\[
\P\!\left(\Delta(\gamma_0-\delta)>\eta_\delta,\ \Delta(\gamma_0+\delta)<-\eta_\delta\right)\longrightarrow 1.
\]
Since $\Delta(\gamma)$ is continuous in $\gamma$ for each realization of $(\vct{y},\mtx{Z})$, on the above event the intermediate value theorem yields at least one root $\hat{\gamma}_n\in(\gamma_0-\delta,\gamma_0+\delta)$ of the equation $\Delta(\gamma)=0$. Therefore,
\[
\P\!\left(|\hat{\gamma}_n-\gamma_0|<\delta\right)\longrightarrow 1.
\]
Because $\delta>0$ is arbitrary, we conclude that there exists a sequence of roots $\hat{\gamma}_n$ satisfying $\hat{\gamma}_n\stackrel{P}{\longrightarrow}\gamma_0$.

\subsubsection*{Consistency of $\hat{\sigma}^2$}
Let's turn to show $\hat{\sigma}_{\varepsilon}^2  \stackrel{P}{\longrightarrow} \sigma_0^2$, where the noise variance estimate is defined in \eqref{eq:variance_mle_modified}. Let
\[
s_n(\gamma) = \frac{1}{n}\vct{y}^\top \mtx{V}_\gamma^{-1} \vct{y},
\]
so that $\hat{\sigma}^2 = s_n(\hat{\gamma})$.

We first evaluate $s_n(\gamma_0)$. Using the same decomposition as in the proof of consistency for $\Delta(\gamma)$, together with \eqref{eq:cond_mean_concentration}, \eqref{eq:Delta_starstar_part1}, \eqref{eq:Delta_starstar_part2}, and Lemma \ref{thm:diff_tr_h}, we obtain
\[
s_n(\gamma_0)-\sigma_0^2
=
\frac{1}{n}\vct{y}^\top \mtx{V}_{\gamma_0}^{-1}\vct{y}
-\sigma_0^2
\stackrel{P}{\longrightarrow}0.
\]
Indeed, conditional on $(\mtx{Z},\vct{\varepsilon})$, the fluctuation part is $o_P(1)$, while the conditional mean converges to
\[
\frac{\|\vct{\beta}\|^2}{np}\trace\!\left(\mtx{V}_{\gamma_0}^{-1}\mtx{Z}\mtx{Z}^\top\right)
+\frac{\sigma_0^2}{n}\trace(\mtx{V}_{\gamma_0}^{-1})
=
\frac{\sigma_0^2\gamma_0}{np}\trace\!\left(\mtx{V}_{\gamma_0}^{-1}\mtx{Z}\mtx{Z}^\top\right)
+\frac{\sigma_0^2}{n}\trace(\mtx{V}_{\gamma_0}^{-1}),
\]
which equals $\sigma_0^2+o_P(1)$ because $\mtx{V}_{\gamma_0}= \mtx{I}_n + \frac{\gamma_0}{p}\mtx{Z}\mtx{Z}^\top$.

It remains to pass from $\gamma_0$ to $\hat{\gamma}$. By the resolvent identity,
\[
\mtx{V}_{\hat{\gamma}}^{-1}-\mtx{V}_{\gamma_0}^{-1}
=
-(\hat{\gamma}-\gamma_0)\,
\mtx{V}_{\hat{\gamma}}^{-1}\Big(\frac{1}{p}\mtx{Z}\mtx{Z}^\top\Big)\mtx{V}_{\gamma_0}^{-1}.
\]
Hence
\begin{align*}
|s_n(\hat{\gamma})-s_n(\gamma_0)|
&\le
|\hat{\gamma}-\gamma_0|
\cdot
\frac{1}{n}\|\vct{y}\|^2
\cdot
\Big\|\frac{1}{p}\mtx{Z}\mtx{Z}^\top\Big\|
\cdot
\|\mtx{V}_{\hat{\gamma}}^{-1}\|
\cdot
\|\mtx{V}_{\gamma_0}^{-1}\|.
\end{align*}
Now $\|\mtx{V}_{\hat{\gamma}}^{-1}\|\le 1$ and $\|\mtx{V}_{\gamma_0}^{-1}\|\le 1$ by positive definiteness. Also, under our assumptions,
\[
\Big\|\frac{1}{p}\mtx{Z}\mtx{Z}^\top\Big\| = O_P(1)
\qquad\text{and}\qquad
\frac{1}{n}\|\vct{y}\|^2 = O_P(1).
\]
Since $\hat{\gamma}\stackrel{P}{\longrightarrow}\gamma_0$, it follows that
\[
|s_n(\hat{\gamma})-s_n(\gamma_0)|=o_P(1).
\]
Combining this with $s_n(\gamma_0)\stackrel{P}{\longrightarrow}\sigma_0^2$, we conclude
\[
\hat{\sigma}^2=s_n(\hat{\gamma})\stackrel{P}{\longrightarrow}\sigma_0^2.
\]

\subsection{Proof of Theorem \ref{thm:consistency_correlated}}

The proof follows the same Rademacher symmetrization argument as in
Theorem \ref{thm:consistency}. We give the details to show that no
independence among the noise coordinates is needed for consistency.

Under the correlated Gaussian assumptions, Lemma \ref{lem:noise_var_bound}
gives
\[
\max_{i\in[n]}\varepsilon_i^2=O_P(\log n)=o_P(\sqrt n),
\qquad
\frac{1}{n}\vct\varepsilon^\top\vct\varepsilon
\stackrel{P}{\longrightarrow}\sigma_0^2.
\]
Conditioning on $(\mtx Z,\vct\varepsilon)$, the Rademacher representation
\eqref{eq:fixed_effects_rad_double_delta} remains valid because it uses only
the row and column sign-invariance of the design matrix. Thus the conditional
variance formula in Lemma \ref{lem:conditional_variance_formula} still holds.

We first show
\[
\Delta(\gamma)-\widetilde\Delta_*(\gamma)=o_P(1)
\]
for each fixed $\gamma>0$. The $V_1$ term in \eqref{eq:var_gamma0} is unchanged. For the $V_2$ term, using
$\|\mtx B_\gamma\|=O_P(n^{-1})$,
$\max_k\|\vct z_k\|^2=O_P(n)$, and
$\|\mtx\Lambda_\varepsilon^2\|=\max_i\varepsilon_i^2=O_P(\log n)$, we get
\[
\sum_{k=1}^p
\beta_k^2\vct z_k^\top
\mtx B_\gamma\mtx\Lambda_\varepsilon^2\mtx B_\gamma
\vct z_k
\leq
\|\vct\beta\|^2
\|\mtx B_\gamma\|^2
\|\mtx\Lambda_\varepsilon^2\|
\max_k\|\vct z_k\|^2
=o_P(1).
\]
For the $V_3$ term, Lemma \ref{lem:B_off_diag} gives
\[
\max_{i\neq j}|(\mtx B_\gamma)_{ij}|^2
=
O_P\!\left(\frac{\log n}{n^3}\right),
\]
and hence
\[
\sum_{i\neq j}
\varepsilon_i^2\varepsilon_j^2(\mtx B_\gamma)_{ij}^2
\leq
n^2\max_i\varepsilon_i^4
\max_{i\neq j}|(\mtx B_\gamma)_{ij}|^2
=
O_P\!\left(\frac{(\log n)^3}{n}\right)
=o_P(1).
\]
Therefore the conditional variance of $\Delta(\gamma)$ given
$(\mtx Z,\vct\varepsilon)$ is $o_P(1)$, and the conditional Chebyshev
argument used in Theorem \ref{thm:consistency} gives
\[
\Delta(\gamma)-\widetilde\Delta_*(\gamma)=o_P(1).
\]

It remains to compare $\widetilde\Delta_*(\gamma)$ with
$\Delta_{**}(\gamma)$. The signal part is exactly the same as in
\eqref{eq:Delta_starstar_part1}. For the noise part,
\[
\trace(\mtx\Lambda_\varepsilon^2\mtx B_\gamma)
-
\frac{1}{n}\trace(\mtx\Lambda_\varepsilon^2)\trace(\mtx B_\gamma)
=
\sum_{i=1}^n
\varepsilon_i^2
\left\{
(\mtx B_\gamma)_{ii}
-
\frac{1}{n}\trace(\mtx B_\gamma)
\right\}.
\]
By Lemma \ref{lem:V_inv_diag},
\[
\max_i
\left|
(\mtx B_\gamma)_{ii}
-
\frac{1}{n}\trace(\mtx B_\gamma)
\right|
=
O_P\!\left(\sqrt{\frac{\log n}{n^3}}\right).
\]
Since $\sum_i\varepsilon_i^2=O_P(n)$, the last display is
$o_P(1)$. Moreover,
\[
\left|
\frac{1}{n}\trace(\mtx\Lambda_\varepsilon^2)\trace(\mtx B_\gamma)
-
\sigma_0^2\trace(\mtx B_\gamma)
\right|
\leq
\left|
\frac{1}{n}\sum_{i=1}^n\varepsilon_i^2-\sigma_0^2
\right|
\,|\trace(\mtx B_\gamma)|
=o_P(1).
\]
Therefore
\[
\widetilde\Delta_*(\gamma)-\Delta_{**}(\gamma)=o_P(1).
\]
The deterministic limit of $\Delta_{**}(\gamma)$ is the same
$c_\gamma$ as in Lemma \ref{lem:Delta_starstar}, which changes sign only at
$\gamma_0$. The same intermediate-value argument used in the proof of
Theorem \ref{thm:consistency} therefore yields a sequence of roots
$\hat\gamma_n$ satisfying
\[
\hat\gamma_n\stackrel{P}{\longrightarrow}\gamma_0.
\]

Finally, the consistency of $\hat\sigma^2$ follows from the same resolvent
identity argument as in Theorem \ref{thm:consistency}. The only stochastic
inputs needed there are
\[
\frac{1}{n}\vct y^\top\mtx V_{\gamma_0}^{-1}\vct y
\stackrel{P}{\longrightarrow}\sigma_0^2,
\qquad
\frac{1}{n}\|\vct y\|^2=O_P(1),
\qquad
\left\|\frac{1}{p}\mtx Z\mtx Z^\top\right\|=O_P(1),
\]
all of which remain valid under the correlated Gaussian assumptions above.
Thus
\[
\hat\sigma^2\stackrel{P}{\longrightarrow}\sigma_0^2.
\]

\subsection{Proof of Theorem \ref{thm:asymp_normality}} 
\label{sec:conv_rate}
Through the analysis of asymptotic distribution, we use the shorthand $h_k = h_k(\gamma_0, \tau)$ for $k=1,2,3,4$, where $h_k(\gamma_0, \tau)$ is defined in
\eqref{eq:h_MP}.

\subsubsection{Decomposition of \texorpdfstring{$\Delta(\gamma_0)$}{Delta(gamma0)}}
The following lemma, adapted from \cite{jiang2016high}, reduces the asymptotic distribution of $\hat{\gamma}$ to that of $\Delta(\gamma_0)$. We give a detailed proof in the appendix for completeness.
\begin{lemma}
\label{lem:Taylor_approx} 
Under the conditions of Theorem \ref{thm:consistency}, assume $\hat{\gamma}_n$ is a sequence of roots of $\Delta(\gamma) = 0$, which converges to $\gamma_0$ in probability. Then
\begin{align}
\label{diff:gamma_hat_gamma0}
\sqrt{n}\left( \hat{\gamma} - \gamma_0 \right)  = - \frac{\sqrt{n} \Delta(\gamma_0)  }{\Delta_{\infty}'(\gamma_0) } + o_P(1),
\end{align}
where $\Delta_{\infty}'(\gamma_0)$ is the limit of $\Delta'(\gamma)$ as $\gamma = \gamma_0$ and has the formula
\begin{align*}
    \Delta'_\infty(\gamma_0) =  \frac{\sigma_0^2}{\gamma_0} \frac{h_1^2 - h_2}{h_1}.
\end{align*}
\end{lemma}

To investigate the asymptotic distribution of $\sqrt{n} \Delta(\gamma_0)$, consider the following orthogonal decomposition:
\[
\Delta(\gamma_0) = (\Delta(\gamma_0) - \widetilde{\Delta}_*(\gamma_0)) + \widetilde{\Delta}_*(\gamma_0).
\]
In other words, the expectation is taken with respect to Rademacher random variables $\xi_i$'s and $\zeta_i$'s. We aim to derive the asymptotic joint distribution of 
\[
\left(\sqrt{n}(\Delta(\gamma_0) - \widetilde{\Delta}_*(\gamma_0)), \sqrt{n}\widetilde{\Delta}_*(\gamma_0)\right).
\]

\subsubsection{Conditional Variance of \texorpdfstring{$\Delta(\gamma_0)-\widetilde{\Delta}_*(\gamma_0)$}{Delta(gamma0)-Delta star(gamma0)}}
\label{sec:nor_part2}
In order to derive the asymptotic joint distribution of $\left(\sqrt{n}(\Delta(\gamma_0) - \widetilde{\Delta}_*(\gamma_0)), \sqrt{n}\widetilde{\Delta}_*(\gamma)\right)$, we first need to study the conditional distribution of $\Delta(\gamma_0) - \widetilde{\Delta}_*(\gamma_0)$ given $\mtx{Z}$ and $\vct{\varepsilon}$. This consists of two steps: conditional variance and conditional normality. Let's first study its conditional variance. Note we have
\[
\Var\left[\sqrt{n}(\Delta(\gamma_0) - \widetilde{\Delta}_*(\gamma_0)) \Big\vert \mtx{Z},\vct{\varepsilon} \right] = \Var\left[\sqrt{n}(\Delta(\gamma_0) ) \Big\vert \mtx{Z},\vct{\varepsilon} \right].
\]

\begin{lemma}
\label{lem:conditional_variance_limit}
Under the condition of Theorem \ref{thm:asymp_normality} or Theorem \ref{thm:asymp_normality_correlated}, with $V_1$, $V_2$ and $V_3$ defined in \eqref{eq:var_gamma0}, we have
\begin{equation}
\label{eq:V1}
V_1 \stackrel{P}{\longrightarrow} 2 \sigma_0^4 \left[ (1 - 2h_1 + h_2) -2\frac{h_1 - 2h_2 + h_3}{h_1} + \frac{h_2 - 2h_3 +h_4}{h_1^2} - \tau\left(h_{1} - \frac{h_{2}}{h_{1}}\right)^2 \right],
\end{equation}
\begin{equation}
\label{eq:V2}
V_2 \stackrel{P}{\longrightarrow} 4 \sigma_0^4 \left( (h_{1} - h_{2}) - 2\frac{ h_{2} - h_{3} }{ h_{1}} + \frac{h_{3} - h_{4}}{ h_{1}^2 } \right),
\end{equation}
\begin{equation}
\label{eq:V3}
V_3 \stackrel{P}{\longrightarrow}2\sigma_0^4\left( h_{2}- \frac{2 h_{3}}{h_{1}} + \frac{h_{4}}{h_{1}^2 } - \left(h_{1} - \frac{h_{2}}{h_{1}}\right)^2\right).
\end{equation}
Consequently,
\[
\Var\left[\sqrt{n}(\Delta(\gamma_0) ) \vert \mtx{Z},\vct{\varepsilon} \right] \stackrel{P}{\longrightarrow} 
2\sigma_0^4 \left(\frac{h_2 -  h_1^2}{ h_1^2} - (\tau + 1)\left(h_{1} - \frac{h_{2}}{h_{1}}\right)^2 \right) .
\]
\end{lemma}

\begin{proof}

\subsubsection*{Limit of $V_1$}
The proof is based on reducing the off-diagonal column quadratic forms to a trace expression. Since
\[
\sum_{i\neq j}\beta_i^2\beta_j^2
=\|\vct\beta\|^4-\sum_{i=1}^p\beta_i^4
=\|\vct\beta\|^4+o(1),
\]
by \eqref{eq:off_diag_fine} in Lemma \ref{lem:leave_two_out},
\begin{align}
\label{eq:diff_mean_ziBzj}
&\left|
2n\sum_{i\neq j}\beta_i^2\beta_j^2
\left(\vct z_i^\top\mtx B_{\gamma_0}\vct z_j\right)^2
-2n\left(\|\vct\beta\|^4-\sum_{i=1}^p\beta_i^4\right)
\frac{1}{p(p-1)}\sum_{i\neq j}
\left(\vct z_i^\top\mtx B_{\gamma_0}\vct z_j\right)^2
\right|=o_P(1).
\end{align}
Next,
\begin{align*}
\trace\left[\left(\mtx B_{\gamma}\frac{1}{p}\mtx Z\mtx Z^\top\right)^2\right]
=\frac{1}{p^2}\left\{\sum_{i\neq j}\left(\vct z_i^\top\mtx B_{\gamma}\vct z_j\right)^2
+\sum_{k=1}^p\left(\vct z_k^\top\mtx B_{\gamma}\vct z_k\right)^2\right\}.
\end{align*}
By \eqref{eq:diag_homo_4} in Lemma \ref{lem:diag_homo},
\begin{align}
\label{eq:offdiag_trace_reduction}
&\left|
\frac{1}{p(p-1)}\sum_{i\neq j}
\left(\vct z_i^\top\mtx B_{\gamma_0}\vct z_j\right)^2
-\frac{p}{p-1}\trace\left[\left(\mtx B_{\gamma_0}\frac{1}{p}\mtx Z\mtx Z^\top\right)^2\right]
+\frac{1}{p-1}\left\{\trace\left(\mtx B_{\gamma_0}\frac{1}{p}\mtx Z\mtx Z^\top\right)\right\}^2
\right|=o_P(n^{-1}).
\end{align}
Lemma \ref{thm:diff_tr_h} implies
\begin{align}
\label{eq:trace_BZZ_sq_limit}
n\trace\left[\left(\mtx B_{\gamma_0}\frac{1}{p}\mtx Z\mtx Z^\top\right)^2\right]
&\stackrel{P}{\longrightarrow}
\frac{1}{\gamma_0^2}\left[
(1-2h_1+h_2)-2\frac{h_1-2h_2+h_3}{h_1}
+\frac{h_2-2h_3+h_4}{h_1^2}
\right],\\
\label{eq:trace_BZZ_limit}
\trace\left(\mtx B_{\gamma_0}\frac{1}{p}\mtx Z\mtx Z^\top\right)
&\stackrel{P}{\longrightarrow}
\frac{1}{\gamma_0}\left(h_1-\frac{h_2}{h_1}\right).
\end{align}
Combining \eqref{eq:diff_mean_ziBzj}, \eqref{eq:offdiag_trace_reduction}, \eqref{eq:trace_BZZ_sq_limit}, and \eqref{eq:trace_BZZ_limit}, and using $n/(p-1)\to\tau$ and $\|\vct\beta\|^4/\gamma_0^2=\sigma_0^4$, gives \eqref{eq:V1}.

\subsubsection*{Limit of $V_2$}
For any fixed $\gamma >0$, straightforward calculation gives
\begin{align}
\label{eq:B_eps_B}
   n \vct{z}_k^\top \mtx{B}_{\gamma}\mtx{\Lambda}_{\varepsilon}^2\mtx{B}_{\gamma} \vct{z}_k 
    & =  \frac{1}{n}\vct{z}_k^\top \mtx{V}_{\gamma}^{-1}\mtx{\Lambda}_{\varepsilon}^2\mtx{V}_{\gamma}^{-1}\vct{z}_k  - 2 \frac{\frac{1}{n} \vct{z}_k^\top \mtx{V}_{\gamma}^{-1}\mtx{\Lambda}_{\varepsilon}^2\mtx{V}_{\gamma}^{-2}\vct{z}_k}{ \frac{1}{n}\trace(\mtx{V}_{\gamma}^{-1}) }+ \frac{\frac{1}{n} \vct{z}_k^\top \mtx{V}_{\gamma}^{-2}\mtx{\Lambda}_{\varepsilon}^2\mtx{V}_{\gamma}^{-2}\vct{z}_k}{\left( \frac{1}{n}\trace(\mtx{V}_{\gamma}^{-1})\right)^2 }.
\end{align}
Then using Sherman-Morrison-Woodbury formula (Theorem \ref{thm:matrix_inv}), we can have
\begin{align}
\label{eq:leave_one_out_with_eps}
    \vct{z}_k^\top \mtx{V}_{\gamma}^{-1}\mtx{\Lambda}_{\varepsilon}^2\mtx{V}_{\gamma}^{-1}\vct{z}_k &= \frac{ \vct{z}_k^\top \mtx{V}_{\gamma,-k}^{-1}\mtx{\Lambda}_{\varepsilon}^2\mtx{V}_{\gamma,-k}^{-1}\vct{z}_k}{\left(1 + \frac{\gamma}{p} \eta_{kk,k}^{(1)} \right)^2 },\notag\\
    \vct{z}_k^\top \mtx{V}_{\gamma}^{-1}\mtx{\Lambda}_{\varepsilon}^2\mtx{V}_{\gamma}^{-2}\vct{z}_k &=  \frac{  \vct{z}_k^\top \mtx{V}_{\gamma,-k}^{-1}\mtx{\Lambda}_{\varepsilon}^2\mtx{V}_{\gamma,-k}^{-2}\vct{z}_k  }{\left(1 + \frac{\gamma}{p} \eta_{kk,k}^{(1)} \right)^2} - \frac{\frac{\gamma}{p} \eta_{kk,k}^{(2)}  \vct{z}_k^\top \mtx{V}_{\gamma,-k}^{-1}\mtx{\Lambda}_{\varepsilon}^2\mtx{V}_{\gamma,-k}^{-1}\vct{z}_k}{\left(1 + \frac{\gamma}{p} \eta_{kk,k}^{(1)} \right)^3},\notag\\
     \vct{z}_k^\top \mtx{V}_{\gamma}^{-2}\mtx{\Lambda}_{\varepsilon}^2\mtx{V}_{\gamma}^{-2}\vct{z}_k &= \frac{ \vct{z}_k^\top \mtx{V}_{\gamma,-k}^{-2}\mtx{\Lambda}_{\varepsilon}^2\mtx{V}_{\gamma,-k}^{-2}\vct{z}_k}{\left(1 + \frac{\gamma}{p} \eta_{kk,k}^{(1)} \right)^2} - \frac{2\frac{\gamma}{p} \eta_{kk,k}^{(2)} \vct{z}_k^\top \mtx{V}_{\gamma,-k}^{-1}\mtx{\Lambda}_{\varepsilon}^2\mtx{V}_{\gamma,-k}^{-2}\vct{z}_k }{\left(1 + \frac{\gamma}{p} \eta_{kk,k}^{(1)} \right)^3}\notag\\
     &~~~ + \frac{ \left(  \frac{\gamma}{p} \eta_{kk,k}^{(2)}\right)^2 \vct{z}_k^\top \mtx{V}_{\gamma,-k}^{-1}\mtx{\Lambda}_{\varepsilon}^2\mtx{V}_{\gamma,-k}^{-1}\vct{z}_k   }{\left(1 + \frac{\gamma}{p} \eta_{kk,k}^{(1)} \right)^4}. 
\end{align}
By Lemma \ref{lem:noise_var_bound}, we have
\[
\|\mtx{V}_{\gamma,-k}^{-l}\mtx{\Lambda}_{\varepsilon}^2\mtx{V}_{\gamma,-k}^{-m}\| \leq \| \mtx{\Lambda}_{\varepsilon}^2\|  = O_P\left( n^{2/(4+\delta)} \right)
\]
and hence
\[
\|\mtx{V}_{\gamma,-k}^{-l}\mtx{\Lambda}_{\varepsilon}^2\mtx{V}_{\gamma,-k}^{-m}\|_F = O_P\left( n^{1/2+2/(4+\delta)} \right).
\]
Then, by Hanson-Wright inequality and taking the uniform bound, we can easily get
\begin{align}
\label{diff:zk_eps_V_zk_tr}
\max_{k \in [p]}\left|  \frac{1}{n} \vct{z}_k^\top \mtx{V}_{\gamma,-k}^{-l}\mtx{\Lambda}_{\varepsilon}^2\mtx{V}_{\gamma,-k}^{-m}\vct{z}_k -\frac{1}{n} \trace\left(  \mtx{\Lambda}_{\varepsilon}^2 \mtx{V}_{\gamma,-k}^{-(l+m)} \right) \right| = o_P(1).
\end{align}
By Lemma \ref{lem:V_inv_diag_loo} and Lemma \ref{lem:noise_var_bound}, we have
\begin{align}
\label{diff:zk_eps_V_zk_tr_2}
\max_{k \in [p]}
\left|\frac{1}{n} \trace\left(  \mtx{\Lambda}_{\varepsilon}^2 \mtx{V}_{\gamma,-k}^{-(l+m)} \right) - \frac{1}{n^2}\trace\left(\mtx{\Lambda}_\varepsilon^2\right)\trace\left(\mtx{V}_\gamma^{-(l+m)}\right) \right| = o_P(1).
\end{align}
By Lemma \ref{lem:noise_var_bound} again, we have
\begin{align}
\label{diff:zk_eps_V_zk_tr_3}
\max_{k \in [p]}
\left|\frac{1}{n^2}
\trace\left(\mtx{\Lambda}_\varepsilon^2\right)\trace\left(\mtx{V}_\gamma^{-(l+m)}\right) -\sigma_0^2\frac{1}{n} \trace\left( \mtx{V}_{\gamma}^{-(l+m)}\right)\right| = o_P(1).
\end{align}
Combining \eqref{diff:zk_eps_V_zk_tr}, \eqref{diff:zk_eps_V_zk_tr_2} and \eqref{diff:zk_eps_V_zk_tr_3} gives
\begin{align}
\label{eq:diff_to_tr_with_no_eps}
    \max_{k \in [p]}\left|  \frac{1}{n} \vct{z}_k^\top \mtx{V}_{\gamma, -k}^{-l}\mtx{\Lambda}_{\varepsilon}^2\mtx{V}_{\gamma, -k}^{-m}\vct{z}_k -\sigma_0^2\frac{1}{n} \trace\left( \mtx{V}_{\gamma}^{-(l+m)}  \right)  \right| = o_P\left(1\right).
\end{align}
Combining \eqref{eq:diff_to_tr_with_no_eps}, \eqref{eq:leave_one_out_with_eps}, \eqref{eq:B_eps_B}, Lemma \ref{lem:leave_one_out} and Lemma \ref{thm:diff_tr_h}, there exists some constant $C(\gamma, \tau)$, such that
\begin{align*}
\max_{k \in [p]}\left| n\vct{z}_k^\top \mtx{B}_{\gamma}
\mtx{\Lambda}_{\varepsilon}^2\mtx{B}_{\gamma}\vct{z}_k - C(\gamma, \tau) \right\vert = o_P(1).
\end{align*}
This further implies both
\begin{align*}
&\left| C(\gamma, \tau) - n\trace\left( \mtx{B}_{\gamma}\mtx{\Lambda}_{\varepsilon}^2\mtx{B}_{\gamma} \frac{1}{p} \mtx{Z} \mtx{Z}^\top  \right)  \right\vert 
\\
&\leq \left|C(\gamma, \tau) -  \frac{n}{p} \sum_{k=1}^p  \vct{z}_k^\top \mtx{B}_{\gamma}\mtx{\Lambda}_{\varepsilon}^2\mtx{B}_{\gamma}\vct{z}_k  \right\vert
\\
& \leq \max_{k \in [p]}  \left|C(\gamma, \tau) - n\vct{z}_k^\top \mtx{B}_{\gamma} \mtx{\Lambda}_{\varepsilon}^2\mtx{B}_{\gamma}\vct{z}_k \right| 
=o_P(1),
\end{align*}
and
\begin{align*}
    &\left|n\sum_{k=1}^p \beta_k^2\vct{z}_k^\top \mtx{B}_{\gamma}\mtx{\Lambda}_{\varepsilon}^2\mtx{B}_{\gamma}\vct{z}_k - \| \beta \|^2 C(\gamma, \tau) \right|
    \\
    &\leq \| \beta \|^2 \max_{k \in [p]}  \left| C(\gamma, \tau)- n\vct{z}_k^\top \mtx{B}_{\gamma}\mtx{\Lambda}_{\varepsilon}^2\mtx{B}_{\gamma}\vct{z}_k \right| =o_P(1).
\end{align*}
Combining the above inequalities, there holds
\begin{equation}
\label{eq:trace_surrogate}
 \left|n\sum_{k=1}^p \beta_k^2\vct{z}_k^\top \mtx{B}_{\gamma}\mtx{\Lambda}_{\varepsilon}^2\mtx{B}_{\gamma}\vct{z}_k 
 - n\|\vct{\beta}\|^2
 \trace\left( \mtx{\Lambda}_{\varepsilon}^2\mtx{B}_{\gamma}\frac{1}{p} \mtx{Z}\mtx{Z}^\top\mtx{B}_{\gamma}\right)\right| = o_P(1).
\end{equation}
~\\

Finally, note that 
\begin{align*}
\mtx{B}_{\gamma}\frac{1}{p} \mtx{Z}\mtx{Z}^\top\mtx{B}_{\gamma} 
&=\frac{1}{\gamma}\Bigg(
\frac{\mtx{V}_\gamma^{-1}}{n^2} 
- \left(\frac{2}{n\trace(\mtx{V}_\gamma^{-1})} + \frac{1}{n^2}\right)\mtx{V}_\gamma^{-2}
\\
&~~~~~~~+ \left(\frac{2}{n\trace(\mtx{V}_\gamma^{-1})} + \frac{1}{\left(\trace(\mtx{V}_\gamma^{-1})\right)^2}\right)\mtx{V}_\gamma^{-3}
- \frac{1}{\left(\trace(\mtx{V}_\gamma^{-1})\right)^2} \mtx{V}_\gamma^{-4}
\Bigg).
\end{align*}
By \eqref{eq:diff_Vii_tr} in Lemma \ref{lem:V_inv_diag}, there holds
\[
\max_{i \in [n]}\left\vert
n\left(\mtx{B}_\gamma \frac{1}{p}\mtx{Z}\mtx{Z}^\top\mtx{B}_\gamma\right)_{ii} - \trace\left(\mtx{B}_\gamma^2 \frac{1}{p}\mtx{Z}\mtx{Z}^\top\right)\right\vert = O_P\left(\sqrt{\frac{\log n}{n^3}}\right).
\]
Similar to \eqref{eq:Delta_starstar_part2}, we have
\begin{equation}
\label{eq:trace_surrogate_2}
\left\vert
n\trace\left( \mtx{\Lambda}_{\varepsilon}^2\mtx{B}_{\gamma}\frac{1}{p} \mtx{Z}\mtx{Z}^\top\mtx{B}_{\gamma}\right) - n\sigma_0^2\trace\left(\mtx{B}_\gamma^2 \frac{1}{p}\mtx{Z}\mtx{Z}^\top\right)
\right\vert = o_P(1).
\end{equation}
Combining \eqref{eq:trace_surrogate} and \eqref{eq:trace_surrogate_2} and letting $\gamma = \gamma_0$, we have
\begin{align}
\label{eq:B_eps_B_tr}
\left| 4n \sum_{k=1}^p \beta_k^2 \vct{z}_k^\top \mtx{B}_{\gamma_0}\mtx{\Lambda}_{\varepsilon}^2\mtx{B}_{\gamma_0} \vct{z}_k  -  4 \sigma_0^2 \|\beta\|^2 n\trace\left( \mtx{B}_{\gamma_0}^2 \frac{1}{p} \mtx{Z} \mtx{Z}^\top \right) \right| =o_P(1).
\end{align}
Notice that
\begin{align*}
 &4\sigma_0^2 \|\beta\|^2 n\trace\left( \mtx{B}_{\gamma_0}^2 \frac{1}{p} \mtx{Z} \mtx{Z}^\top \right) 
 \\
 &= \frac{4 \sigma_0^2 \|\beta\|^2 }{\gamma_0} n\trace\left( \mtx{B}_{\gamma_0}^2 \left(\mtx{V}_{\gamma_0} - \mtx{I}_n \right) \right)
 \\
 &=4\sigma_0^4 n \trace\left(
 \frac{\mtx{V}_{\gamma_0}^{-1} - \mtx{V}_{\gamma_0}^{-2}}{n^2}
 -2\frac{\mtx{V}_{\gamma_0}^{-2} - \mtx{V}_{\gamma_0}^{-3}}{n\trace(\mtx{V}_{\gamma_0}^{-1})}
 + \frac{\mtx{V}_{\gamma_0}^{-3} - \mtx{V}_{\gamma_0}^{-4}}{(\trace(\mtx{V}_{\gamma_0}^{-1}))^2}
 \right).
\end{align*}
Therefore, we got \eqref{eq:V2} by Lemma \ref{thm:diff_tr_h}.

\subsubsection*{Limit of $V_3$}
By \eqref{eq:Bij_weight_sum} in Lemma \ref{lem:B_off_diag}, \eqref{eq:diff_Bii_tr_sq} in Lemma \ref{lem:V_inv_diag}, and Lemma \ref{thm:diff_tr_h}, we have
\begin{align*}
\label{eq:cond_mean_sigmaij}
V_3 &= (2n)\sum_{i \neq j} \varepsilon_i^2 \varepsilon_j^2 (\mtx{B}_{\gamma_0})_{ij}^2 
\\
&= 2\bar{\theta}_2(\gamma, \tau) \sigma_0^4 + o_P(1) \notag
\\
&= \sigma_0^4 \left(2n\sum_{i \neq j} (\mtx{B}_{\gamma_0})_{ij}^2\right) + o_P(1) \notag
\\
&= \sigma_0^4  \left(2n \trace\left(\mtx{B}_{\gamma_0}^2 \right) -  2n \sum_{i=1}^n (\mtx{B}_{\gamma_0})_{ii}^2\right) + o_P(1) \notag
\\
&= \sigma_0^4  \left(2n \trace\left(\mtx{B}_{\gamma_0}^2 \right) 
-  2\left(\trace\left(\mtx{B}_{\gamma_0} \right)\right)^2 \right) + o_P(1) \notag
\\
& \stackrel{P}{\longrightarrow}2\sigma_0^4
\left(\left( h_{2}- \frac{2 h_{3}}{h_{1}} + \frac{h_{4}}{h_{1}^2 } \right) - \left(h_{1} - \frac{h_{2}}{h_{1}}\right)^2\right).
\end{align*}

\end{proof}

\subsubsection{Conditional Distribution of \texorpdfstring{$\Delta(\gamma_0)-\widetilde{\Delta}_*(\gamma_0)$}{Delta(gamma0)-Delta star(gamma0)}}
\label{sec:cond_dist}
Recall from \eqref{eq:fixed_effects_rad_double_delta} that 
\begin{align*}
\sqrt{n}\Delta(\gamma_0) = \left(\vct{\xi}^\top \quad \vct{\zeta}^\top \right)
\underbrace{\sqrt{n}
\begin{bmatrix}
\mtx{\Lambda}_\beta \mtx{Z}^\top \mtx{B}_{\gamma_0} \mtx{Z} \mtx{\Lambda}_\beta 
& 
\mtx{\Lambda}_\beta \mtx{Z}^\top \mtx{B}_{\gamma_0} \mtx{\Lambda}_\varepsilon
\\
\mtx{\Lambda}_\varepsilon \mtx{B}_{\gamma_0} \mtx{Z} \mtx{\Lambda}_\beta
& 
\mtx{\Lambda}_\varepsilon \mtx{B}_{\gamma_0} \mtx{\Lambda}_\varepsilon
\end{bmatrix}}_{\mtx{Q}}
\begin{pmatrix}
\vct{\xi}
\\
\vct{\zeta}.
\end{pmatrix}
\end{align*}
Note that we can represent $\mtx{Q}$ as
\begin{align*}
\mtx{Q}
=\begin{pmatrix} 
\mtx{\Lambda}_\beta&0
\\
0&\frac{1}{\sqrt{n}}\mtx{\Lambda}_\varepsilon  
\end{pmatrix} 
\begin{pmatrix}
\mtx{Z}^\top 
\\
\sqrt{n}\mtx{I}_n
\end{pmatrix} 
\sqrt{n}\mtx{B}_{\gamma_0} \left( \mtx{Z} \quad \sqrt{n}\mtx{I}_n  \right) 
\begin{pmatrix} 
\mtx{\Lambda}_\beta&0
\\
0&\frac{1}{\sqrt{n}}\mtx{\Lambda}_\varepsilon  
\end{pmatrix}.
\end{align*}
Since $\widetilde{\Delta}_*(\gamma_0) = \E[\Delta(\gamma_0)\vert \mtx{Z},\vct{\varepsilon}]$, we have
\[
\sqrt{n}(\Delta(\gamma_0) - \widetilde{\Delta}_*(\gamma_0)) = \left(\vct{\xi}^\top \quad \vct{\zeta}^\top \right) (\mtx{Q} - \breve{\mtx{Q}}) \begin{pmatrix}
\vct{\xi}
\\
\vct{\zeta} 
\end{pmatrix}.
\]
where $\breve{\mtx{Q}}$ is a diagonal matrix that maintains the diagonal part of $\mtx{Q}$. In other words, conditional on $\mtx{Z}$ and $\vct{\varepsilon}$, $\sqrt{n}(\Delta(\gamma_0) - \widetilde{\Delta}_*(\gamma_0))$ is a quadratic form about $\vct{\xi}$ and $\vct{\zeta}$.

Here we use Theorem \ref{thm:normality} in the appendix to establish a normal approximation of the conditional distribution of $\sqrt{n}(\Delta(\gamma_0) - \widetilde{\Delta}_*(\gamma_0))$ given $\mtx{Z}$ and $\vct{\varepsilon}$. In other words, we show
\begin{equation}
\label{eq:oper_frobe_ratio}
\frac{\| \mtx{Q} - \breve{\mtx{Q}}\|}{\|\mtx{Q} - \breve{\mtx{Q}}\|_F}   = o_P(1).
\end{equation}

To establish a lower bound on $\|\mtx{Q}-\breve{\mtx{Q}}\|_F$, we use the conditional variance identity for Rademacher quadratic forms. Since
\[
\sqrt n\{\Delta(\gamma_0)-\widetilde\Delta_*(\gamma_0)\}
=\left(\vct\xi^\top\quad \vct\zeta^\top\right)(\mtx Q-\breve{\mtx Q})
\begin{pmatrix}\vct\xi\\ \vct\zeta\end{pmatrix},
\]
we have
\[
\Var\left[\sqrt n\{\Delta(\gamma_0)-\widetilde\Delta_*(\gamma_0)\}\mid \mtx Z,\vct\varepsilon\right]
=2\|\mtx Q-\breve{\mtx Q}\|_F^2.
\]
By Lemma \ref{lem:conditional_variance_limit}, the left-hand side converges in probability to $\nu_2$. Thus $\|\mtx Q-\breve{\mtx Q}\|_F^2$ is bounded away from zero with probability tending to one whenever $\nu_2>0$.

Let's now establish an upper bound of $\|\mtx{Q} - \breve{\mtx{Q}}\|$. First, we have 
\[
\|\mtx{Q} - \breve{\mtx{Q}}\| \leq \|\mtx{Q}\| + \|\breve{\mtx{Q}}\| \leq 2 \|\mtx{Q}\|,
\]
where the last inequality is due to the fact that all diagonal entries are bounded by the operator norm in magnitude. On the other hand,
\begin{align*}
    \|\mtx{Q}\| &\leq \left\Vert \begin{pmatrix} \mtx{\Lambda}_\beta^2&0\\0&\frac{1}{n}\mtx{\Lambda}_\varepsilon^2  \end{pmatrix}  \right\Vert \left\Vert \begin{pmatrix}\mtx{Z}^\top \\\sqrt{n}
\mtx{I}_n\end{pmatrix}  \right\Vert^2 \left\Vert \sqrt{n} \mtx{B}_{\gamma_0} \right\Vert\\
& \leq \max_{ k \in [p],~ i \in [n]} \left\{\beta_k^2, \frac{1}{n} \varepsilon_i^2 \right\} \cdot \left(  \| \mtx{Z} \|  + \sqrt{n} \right)^2 \cdot \left\Vert  \sqrt{n}\mtx{B}_{\gamma_0} \right\Vert.
\end{align*}
By the deterministic assumption $\| \beta\|^2_{\infty} = o(n^{-1/2})$ and the bound given in \eqref{eq:max_esp_i}, we have
\[
\max_{ k \in [p],~ i \in [n]} \left\{\beta_k^2, \frac{1}{n} \varepsilon_i^2 \right\} = o_P(n^{-1/2}).
\]
Note that we have $\| \mtx{Z} \| = O_P(\sqrt{n})$ by Theorem \ref{thm:eig_subG}. Also, Lemma \ref{thm:diff_tr_h} implies $\left\Vert \sqrt{n} \mtx{B}_{\gamma_0} \right\Vert = O_P(n^{-1/2})$. Combining the above, we have $\| \mtx{Q}\| = o_P(1)$. This completes the proof of \eqref{eq:oper_frobe_ratio}.

Therefore, by Theorem \ref{thm:normality}, we have
\begin{align}
\label{eq:cond_dist}
    \P\left\{ \frac{\sqrt{n}(\Delta(\gamma_0) - \widetilde{\Delta}_*(\gamma_0))}{\sqrt{\Var\left[\sqrt{n}\Delta(\gamma_0)\vert \mtx{Z}, \vct{\varepsilon} \right]}} \leq t \bigg\vert \mtx{Z}, \vct{\varepsilon}\right\} \stackrel{P}{\longrightarrow} \Phi(t),
\end{align}
where $\Phi(t)$ is the c.d.f of standard normal distribution.

\subsubsection{Asymptotic Distribution of \texorpdfstring{$\sqrt{n}\widetilde{\Delta}_*(\gamma_0)$}{sqrt(n) Delta star(gamma0)}}
This subsection is intended to show the following result that characterizes the asymptotic distribution of $\widetilde{\Delta}_*(\gamma_0)$ defined in \eqref{eq:tilde_delta_star}.
Note that by the definition of $\mtx{B}_\gamma$, we always have $\trace(\mtx{B}_\gamma \mtx{V}_\gamma) = 0$. Therefore,
\begin{align*}
\| \vct{\beta} \|^2\trace\left(\mtx{B}_{\gamma_0} \frac{1}{p} \mtx{Z} \mtx{Z}^\top
 \right)  +  \sigma_0^2 \trace\left( \mtx{B}_{\gamma_0}\right) &=  \sigma_0^2 \trace\left(\mtx{B}_{\gamma_0} \mtx{V}_{\gamma_0}\right) = 0.
\end{align*}
Then, we can represent $\widetilde{\Delta}_*(\gamma_0)$ as
\begin{align}
\label{eq:diff_tiDeltas_Deltass}
\widetilde{\Delta}_*(\gamma_0) &= \sum_{k=1}^p \beta_k^2 \vct{z}_k^\top \mtx{B}_{\gamma_0}\vct{z}_k - \| \vct{\beta} \|^2\trace\left(\mtx{B}_{\gamma_0} \frac{1}{p} \mtx{Z} \mtx{Z}^\top \right) \notag
\\
&~~~~~~+ \trace\left( \mtx{\Lambda}_\varepsilon^2 \mtx{B}_{\gamma_0}\right) - \sigma_0^2 \trace\left(\mtx{B}_{\gamma_0}\right)\notag 
\\
& = \trace\left( (\mtx{\Lambda}_\varepsilon^2 - \sigma_0^2 \mtx{I}_n)\mtx{B}_{\gamma_0}\right) \notag
\\
&~~~~+ \sum_{k=1}^p \beta_k^2 \frac{1}{n}\left(\vct{z}_k^\top \mtx{V}_{\gamma_0}^{-1}\vct{z}_k - \trace\left(\mtx{V}_{\gamma_0}^{-1} \frac{1}{p} \mtx{Z} \mtx{Z}^\top \right)\right) \notag
\\
& ~~~~- \frac{n}{\trace(\mtx{V}_{\gamma_0}^{-1})}\sum_{k=1}^p \beta_k^2 \frac{1}{n}\left(\vct{z}_k^\top \mtx{V}_{\gamma_0}^{-2}\vct{z}_k - \trace\left(\mtx{V}_{\gamma_0}^{-2} \frac{1}{p} \mtx{Z} \mtx{Z}^\top \right)\right).
\end{align}

By Lemma \ref{lem:diag_homo_taylor} and Lemma \ref{lem:leave_one_out}, we have
\begin{align}
\label{eq:diag_fine_approx_1}
&\sum_{k=1}^p \beta_k^2 \frac{1}{n}\left(\vct{z}_k^\top \mtx{V}_{\gamma_0}^{-1}\vct{z}_k - \trace\left(\mtx{V}_{\gamma_0}^{-1} \frac{1}{p} \mtx{Z} \mtx{Z}^\top \right)\right) \notag
\\
&=\frac{1}{\left(1+ \frac{\gamma}{p}\trace(\mtx{V}_{\gamma}^{-1})\right)^2}
\sum_{k=1}^p \beta_k^2\left(\frac{1}{n}\eta_{kk,k}^{(1)} - \frac{1}{n}\trace\left(\mtx{V}_{\gamma,-k}^{-1}\right)\right) + O_P\left(\frac{\log n}{n}\right)
\end{align}
and
\begin{align}
\label{eq:diag_fine_approx_2}
&\sum_{k=1}^p \beta_k^2 \frac{1}{n}\left(\vct{z}_k^\top \mtx{V}_{\gamma_0}^{-2}\vct{z}_k - \trace\left(\mtx{V}_{\gamma_0}^{-2} \frac{1}{p} \mtx{Z} \mtx{Z}^\top \right)\right) \notag
\\
&=-\frac{2\gamma \trace(\mtx{V}_\gamma^{-2})}{p\left(1+ \frac{\gamma}{p}\trace(\mtx{V}_{\gamma}^{-1})\right)^3}
\sum_{k=1}^p \beta_k^2\left(\frac{1}{n}\eta_{kk,k}^{(1)} - \frac{1}{n}\trace\left(\mtx{V}_{\gamma,-k}^{-1}\right)\right) \notag
\\
&~~~~~~~+\frac{1}{\left(1+ \frac{\gamma}{p}\trace(\mtx{V}_\gamma^{-1})\right)^2}
\sum_{k=1}^p \beta_k^2\left(\frac{1}{n}\eta_{kk,k}^{(2)} - \frac{1}{n}\trace\left(\mtx{V}_{\gamma,-k}^{-2}\right)\right) + O_P\left(\frac{\log n}{n}\right).
\end{align}
Then, for $l=1,2$,
\begin{align}
\label{eq:diag_fine_approx_3}
&\E \left(
\sum_{k=1}^p \beta_k^2\left(\frac{1}{n}\eta_{kk,k}^{(l)} - \frac{1}{n}
\trace\left(\mtx{V}_{\gamma,-k}^{-l}\right)\right) \notag
\right)^2
\\
&= \sum_{k=1}^p \beta_k^4 
\E \left[\left(\frac{1}{n}\eta_{kk,k}^{(l)} - \frac{1}{n}
\trace\left(\mtx{V}_{\gamma,-k}^{-l}\right)\right)^2\right] \notag
\\
&~~~~+ \sum_{i \neq j} \beta_i^2 \beta_j^2 
\E \Bigg[\left(\frac{1}{n}\eta_{ii,i}^{(l)} 
- \frac{1}{n}\trace\left(\mtx{V}_{\gamma,-i}^{-l}\right)\right) 
\left(\frac{1}{n}\eta_{jj,j}^{(l)} - \frac{1}{n}
\trace\left(\mtx{V}_{\gamma,-j}^{-l}\right)\right) \Bigg] \notag
\\
&\leq
\frac{C}{n}\|\vct{\beta}\|_4^4
+
\frac{C}{n^{3/2}}\|\vct{\beta}\|_2^4
=
o\left(\frac{1}{n}\right),
\end{align}
for which we have used Lemma \ref{lem:diag_homo_moment} as well as the fact $\|\vct{\beta}\|_4 = o(1)$ (by the assumption $\| \vct{\beta} \|_\infty = o(p^{-1/4})$). Then by \eqref{eq:diag_fine_approx_1}, \eqref{eq:diag_fine_approx_2}, and \eqref{eq:diag_fine_approx_3}, in connection with Lemma \ref{thm:diff_tr_h}, there holds
\begin{equation}
\label{eq:diag_fine_approx_4}
\sum_{k=1}^p \beta_k^2 \frac{1}{n}\left(\vct{z}_k^\top \mtx{V}_{\gamma_0}^{-l}\vct{z}_k - \trace\left(\mtx{V}_{\gamma_0}^{-l} \frac{1}{p} \mtx{Z} \mtx{Z}^\top \right)\right) = o_P\left(\frac{1}{\sqrt{n}}\right), \quad l=1,2.
\end{equation}
Then, equation \eqref{eq:diff_tiDeltas_Deltass} implies
\begin{align}
\label{eq:cond_mean_simp}
\widetilde{\Delta}_*(\gamma_0) = \trace\left( (\mtx{\Lambda}_\varepsilon^2 - \sigma_0^2 \mtx{I}_n)\mtx{B}_{\gamma_0}\right) + o_P\left(\frac{1}{\sqrt{n}}\right).
\end{align}

Before deriving the asymptotic distribution of $\widetilde{\Delta}_*(\gamma_0)$, we first introduce the following row leave-out approximation:

\begin{lemma}
\label{lem:noise_diag_approx_fine}
Under the conditions of Theorem \ref{thm:asymp_normality} or Theorem \ref{thm:asymp_normality_correlated}, for any fixed $\gamma>0$, there holds
\begin{align}\label{eq:noise_diag_approx_fine}
\sum_{i=1}^n (\varepsilon_i^2 - \sigma_0^2) \left( (\mtx{B}_\gamma)_{ii} - \frac{1}{n}\trace\left(\mtx{B}_\gamma\right)  \right) = o_P\left(\frac{1}{\sqrt{n}}\right).
\end{align}
\end{lemma}

With this lemma,  equation \eqref{eq:cond_mean_simp} gives
\begin{align*}
\label{eq:diff_tiDeltas_Deltass2}
\widetilde{\Delta}_*(\gamma_0)&=\sum_{i =1}^n (\varepsilon_i^2 - \sigma_0^2) \left(\mtx{B}_{\gamma_0}\right)_{ii}+ o_P\left(\frac{1}{\sqrt{n}}\right) \notag
\\
&= \frac{1}{n} \left(\sum_{i =1}^n (\varepsilon_i^2 - \sigma_0^2) \right)
\trace\left(\mtx{B}_{\gamma_0}\right) + o_P\left(\frac{1}{\sqrt{n}}\right).
\end{align*}
Recall that Lemma \ref{thm:diff_tr_h} implies $\trace\left(\mtx{B}_{\gamma_0}\right) \stackrel{P}{\longrightarrow} h_1 - \frac{h_2}{h_1}$. Moreover, \eqref{eq:noise_CLT} in Lemma \ref{lem:noise_var_bound} gives
\[
\frac{1}{\sqrt{n}} \left(\left(\sum_{i=1}^n \varepsilon_i^2\right) - n\sigma_0^2\right) \Longrightarrow \mathcal{N}(0, 2\kappa_\varepsilon \sigma_0^4).
\]
Then, by Slutsky's theorem, we get
\begin{align}
\label{eq:con_mean}
\sqrt{n}\widetilde{\Delta}_*(\gamma_0) = \sqrt{n}\E[\Delta(\gamma_0)\vert\mtx{Z},\vct{\varepsilon}]\Longrightarrow
\mathcal{N}\left(0, 2 \kappa_\varepsilon \sigma_0^4 \left(\frac{ h_2}{h_1} - h_1\right)^2\right).
\end{align}

\subsubsection{Asymptotic Distribution of \texorpdfstring{$\hat{\gamma}$}{gamma-hat}}
Denote 
\begin{equation}
\label{eq:nu1_nu2}
\begin{cases}
\nu_1 = 2 \kappa_\varepsilon \sigma_0^4 \left(\frac{ h_2}{h_1} - h_1\right)^2 
\\
\nu_2 = 2\sigma_0^4 \left(\frac{h_2 -  h_1^2}{ h_1^2} - (\tau + 1)\left(h_{1} - \frac{h_{2}}{h_{1}}\right)^2 \right),
\end{cases}
\end{equation}
which are the asymptotic variances given in \eqref{eq:con_mean} and Lemma \ref{lem:conditional_variance_limit}.
If $\nu_1=0$, the following joint convergence argument is interpreted with the first standardized component omitted; the displayed conclusion is unchanged.

To establish the asymptotic distribution of $\hat{\gamma}$, we only need to find that of $\sqrt{n}\Delta(\gamma_0)$ by
Lemma \ref{lem:Taylor_approx}. Furthermore, it suffices to find the asymptotic joint distribution of 
\[
\left(\sqrt{n}(\Delta(\gamma_0) - \widetilde{\Delta}_*(\gamma_0)), \sqrt{n}\widetilde{\Delta}_*(\gamma_0)\right). 
\]
For any $(t,s) \in \mathbb{R}^2$, we have
\begin{align*}
    &~~~\P \left\{  \frac{ \sqrt{n}\widetilde{\Delta}_*(\gamma_0)}{\sqrt{\nu_1}} \leq t, \frac{\sqrt{n}\Delta(\gamma_0) -  \sqrt{n}\widetilde{\Delta}_*(\gamma_0)}{\sqrt{\Var\left[\sqrt{n}\Delta(\gamma_0)\vert \mtx{Z}, \vct{\varepsilon} \right]}} \leq s \right\}
    \\
    & = \E\left[\P \left\{\frac{ \sqrt{n}\widetilde{\Delta}_*(\gamma_0)}{\sqrt{\nu_1}} \leq t, \frac{\sqrt{n}\Delta(\gamma_0) -  \sqrt{n}\widetilde{\Delta}_*(\gamma_0)}{\sqrt{\Var\left[\sqrt{n}\Delta(\gamma_0)\vert \mtx{Z}, \vct{\varepsilon} \right]}} \leq s \bigg\vert \mtx{Z},\vct{\varepsilon}  \right\}\right] 
    \\
    & = \E\left[ \mathbbm{1}_{\{ \sqrt{n}\widetilde{\Delta}_*(\gamma_0)/ \sqrt{\nu_1}\leq t \}}\P \left\{ \frac{\sqrt{n}\Delta(\gamma_0) -  \sqrt{n}\widetilde{\Delta}_*(\gamma_0)}{\sqrt{\Var\left[\sqrt{n}\Delta(\gamma_0)\vert \mtx{Z}, \vct{\varepsilon} \right]}} \leq s \bigg\vert \mtx{Z},\vct{\varepsilon} \right\}\right].
\end{align*}
Note that
\begin{align*}
&\left|\E\left[ \mathbbm{1}_{ \{\sqrt{n}\widetilde{\Delta}_*(\gamma_0)/ \sqrt{\nu_1} \leq t \}}\P \left\{ \frac{\sqrt{n}\Delta(\gamma_0) -  \sqrt{n}\widetilde{\Delta}_*(\gamma_0)}{\sqrt{\Var\left[\sqrt{n}\Delta(\gamma_0)\vert \mtx{Z}, \vct{\varepsilon} \right]}} \leq s \bigg\vert \mtx{Z},\vct{\varepsilon} \right\}\right] - \E\left[ \mathbbm{1}_{ \{\sqrt{n}\widetilde{\Delta}_*(\gamma_0)/ \sqrt{\nu_1} \leq t \}}\Phi(s)\right]\right|
\\
& \leq \E\left[ \mathbbm{1}_{ \{\sqrt{n}\widetilde{\Delta}_*(\gamma_0)/ \sqrt{\nu_1} \leq t\}} \left|\P \left\{ \frac{\sqrt{n}\Delta(\gamma_0) -  \sqrt{n}\widetilde{\Delta}_*(\gamma_0)}{\sqrt{\Var\left[\sqrt{n}\Delta(\gamma_0)\vert \mtx{Z}, \vct{\varepsilon} \right]}} \leq s \bigg\vert \mtx{Z},\vct{\varepsilon} \right\}-\Phi(s)\right|\right] 
\\
& \leq \E\left[ \left|\P \left\{ \frac{\sqrt{n}\Delta(\gamma_0) -  \sqrt{n}\widetilde{\Delta}_*(\gamma_0)}{\sqrt{\Var\left[\sqrt{n}\Delta(\gamma_0)\vert \mtx{Z}, \vct{\varepsilon} \right]}} \leq s \bigg\vert \mtx{Z},\vct{\varepsilon} \right\}-\Phi(s)\right|\right] \rightarrow 0,
\end{align*}
where the last inequality is due to \eqref{eq:cond_dist}. By \eqref{eq:con_mean}, we have
\begin{align*}
     \E\left[ \mathbbm{1}_{\{ \sqrt{n}\widetilde{\Delta}_*(\gamma_0)/ \sqrt{\nu_1} \leq t\}}\Phi(s)\right] 
     = \P \left\{\frac{\sqrt{n}\widetilde{\Delta}_*(\gamma_0)}{\sqrt{\nu_1}} \leq t \right\} \Phi(s) \rightarrow \Phi(t)\Phi(s).
\end{align*}
Thus we can have
\begin{align*}
\P \left\{ \frac{ \sqrt{n}\widetilde{\Delta}_*(\gamma_0)}{\sqrt{\nu_1}} \leq t,~\frac{\sqrt{n}\Delta(\gamma_0) -  \sqrt{n}\widetilde{\Delta}_*(\gamma_0)}{\sqrt{\Var\left[\sqrt{n}\Delta(\gamma_0)\vert \mtx{Z}, \vct{\varepsilon} \right]}} \leq s \right\} \rightarrow \Phi(t)\Phi(s),
\end{align*}
which implies that
\[
\left(\frac{ \sqrt{n}\widetilde{\Delta}_*(\gamma_0)}{\sqrt{\nu_1}},  \frac{\sqrt{n}\Delta(\gamma_0) -  \sqrt{n}\widetilde{\Delta}_*(\gamma_0)}{\sqrt{\Var\left[\sqrt{n}\Delta(\gamma_0)\vert \mtx{Z}, ~
\vct{\varepsilon} \right]}} \right) \Longrightarrow \left( X_1 , X_2\right),
\]
where $[X_1, X_2] \sim \mathcal{N}_2(\vct{0}, \mtx{I}_2)$. By Lemma \ref{lem:conditional_variance_limit}, we have
\[
\Var\left[\sqrt{n}\Delta(\gamma_0)\vert \mtx{Z}, \vct{\varepsilon} \right] \stackrel{P}{\longrightarrow} \nu_2.
\]
Then, Slutsky's theorem implies
\[
\left( \frac{ \sqrt{n}\widetilde{\Delta}_*(\gamma_0)}{\sqrt{\nu_1}},  \frac{\sqrt{n}\Delta(\gamma_0) -  \sqrt{n}\widetilde{\Delta}_*(\gamma_0)}{\sqrt{\Var\left[\sqrt{n}\Delta(\gamma_0)\vert \mtx{Z}, \vct{\varepsilon} \right]}}, \sqrt{\Var\left[\sqrt{n}\Delta(\gamma_0)\vert \mtx{Z}, \vct{\varepsilon} \right]} \right) \Longrightarrow \left( X_1 , X_2, \sqrt{\nu_2}\right).
\]
Letting $g(x,y,z) =  \sqrt{\nu_1} x + yz$, by the continuous mapping theorem, we can have
\begin{align*}
    \sqrt{n}\Delta(\gamma_0) &=  \sqrt{\nu_1} \frac{ \sqrt{n}\widetilde{\Delta}_*(\gamma_0)}{\sqrt{\nu_1}}+ \frac{\sqrt{n}\Delta(\gamma_0) -  \sqrt{n}\widetilde{\Delta}_*(\gamma_0)}{\sqrt{\Var\left[\sqrt{n}\Delta(\gamma_0)\vert \mtx{Z}, \vct{\varepsilon} \right]}} \cdot \sqrt{\Var\left[\sqrt{n}\Delta(\gamma_0)\vert \mtx{Z}, \vct{\varepsilon} \right]} \\
    & \Longrightarrow \sqrt{\nu_1} X_1 + \sqrt{\nu_2} X_2.
\end{align*}
Finally, by Lemma \ref{lem:Taylor_approx} and Slutsky's theorem, we have 
\begin{align*}
    \sqrt{n}\left( \hat{\gamma} - \gamma_0 \right) \Longrightarrow \mathcal{N}\left(0, \frac{\nu_1 + \nu_2}{(\Delta_{\infty}'(\gamma_0))^2  }\right)
\end{align*}
Let $D=h_2-h_1^2>0$. By \eqref{eq:nu1_nu2},
\[
\nu_1+\nu_2
=
2\sigma_0^4\left\{
\frac{D}{h_1^2}+(\kappa_\varepsilon-\tau-1)\frac{D^2}{h_1^2}
\right\}.
\]
Since
\[
\left\{\Delta'_\infty(\gamma_0)\right\}^2
=\frac{\sigma_0^4}{\gamma_0^2}\frac{D^2}{h_1^2},
\]
we have
\begin{align*}
\sqrt{n}\left( \hat{\gamma} - \gamma_0 \right) \Longrightarrow \mathcal{N}\left(0, 2\gamma_0^2\left(\frac{1}{ h_2 - h_1^2}  +\kappa_\varepsilon - \tau - 1\right) \right).
\end{align*}

This proves Theorem \ref{thm:asymp_normality}.

\subsection{Proof of Theorem \ref{thm:asymp_normality_correlated}}

We spell out the modifications relative to the proof of
Theorem \ref{thm:asymp_normality}. The Rademacher representation
\eqref{eq:fixed_effects_rad_double_delta} is unchanged, because it relies only
on the sign-invariance of the design matrix. Thus, conditional on
$(\mtx Z,\vct\varepsilon)$, the same quadratic-form representation and the
same conditional variance formula in Lemma
\ref{lem:conditional_variance_formula} hold.

First, the conditional normal approximation of
\[
\sqrt n\{\Delta(\gamma_0)-\widetilde\Delta_*(\gamma_0)\}
\]
is unchanged. Indeed, the proof of \eqref{eq:oper_frobe_ratio} only uses
\[
\max_i\varepsilon_i^2=O_P(\log n),
\qquad
\|\vct\beta\|_\infty=o(p^{-1/4}),
\qquad
\|\mtx Z\|=O_P(\sqrt n),
\qquad
\|\sqrt n\,\mtx B_{\gamma_0}\|=O_P(n^{-1/2}),
\]
and these remain true under the correlated Gaussian assumptions. Hence
\[
\P\left\{
\frac{
\sqrt n\{\Delta(\gamma_0)-\widetilde\Delta_*(\gamma_0)\}
}{
\sqrt{\Var[\sqrt n\Delta(\gamma_0)\mid \mtx Z,\vct\varepsilon]}
}
\leq t
\Bigm| \mtx Z,\vct\varepsilon
\right\}
\stackrel{P}{\longrightarrow}
\Phi(t).
\]

Second, the limiting conditional variance is the same as in
Lemma \ref{lem:conditional_variance_limit}. The terms $V_1$ and $V_2$ are
unchanged. For $V_3$, write
\[
V_3
=
2n\sum_{i\neq j}
\varepsilon_i^2\varepsilon_j^2(\mtx B_{\gamma_0})_{ij}^2.
\]
Since $\vct\varepsilon$ is Gaussian,
\[
\E[\varepsilon_i^2\varepsilon_j^2]
=
\sigma_i^2\sigma_j^2+2(\mtx\Sigma_\varepsilon)_{ij}^2,
\qquad i\neq j.
\]
The second term is negligible because
\[
0\leq
2n\sum_{i\neq j}
(\mtx\Sigma_\varepsilon)_{ij}^2(\mtx B_{\gamma_0})_{ij}^2
\leq
2n\max_{i\neq j}(\mtx B_{\gamma_0})_{ij}^2
\|\mtx\Sigma_\varepsilon\|_F^2
=
o_P(1),
\]
where we used
$\max_{i\neq j}(\mtx B_{\gamma_0})_{ij}^2
=O_P(\log n/n^3)$ and
$\|\mtx\Sigma_\varepsilon\|_F^2=O(n)$.
The remaining part satisfies
\[
2n\sum_{i\neq j}
\sigma_i^2\sigma_j^2(\mtx B_{\gamma_0})_{ij}^2
=
2\sigma_0^4 n\sum_{i\neq j}
(\mtx B_{\gamma_0})_{ij}^2+o_P(1),
\]
by the row-homogeneity bound in Lemma \ref{lem:B_off_diag},
$\max_i\sigma_i^2=O(1)$, and
$n^{-1}\sum_i\sigma_i^2=\sigma_0^2$.
Finally, the fluctuation around the mean is negligible. One way to see this
is to use the Gaussian Poincare inequality conditionally on $\mtx Z$ for
\[
F(\vct\varepsilon)
=
n\sum_{i\neq j}
(\mtx B_{\gamma_0})_{ij}^2\varepsilon_i^2\varepsilon_j^2.
\]
Since
\[
\frac{\partial F}{\partial \varepsilon_i}
=
4n\varepsilon_i\sum_{j\neq i}
(\mtx B_{\gamma_0})_{ij}^2\varepsilon_j^2,
\]
we have
\[
\|\nabla F(\vct\varepsilon)\|^2
\leq
C n^2
\sum_{i=1}^n
\varepsilon_i^2
\left\{
\sum_{j\neq i}(\mtx B_{\gamma_0})_{ij}^2\varepsilon_j^2
\right\}^2.
\]
Using $\|\mtx\Sigma_\varepsilon\|=O(1)$ and Gaussian moment bounds,
together with
\[
\max_i\sum_{j\neq i}(\mtx B_{\gamma_0})_{ij}^2
=
O_P\!\left(\frac{\log n}{n^2}\right),
\qquad
n\sum_{i\neq j}(\mtx B_{\gamma_0})_{ij}^2=O_P(1),
\]
gives
\[
\Var(F\mid \mtx Z)=o_P(1).
\]
Thus $V_3$ has the same limit as in Lemma
\ref{lem:conditional_variance_limit}. Consequently,
\[
\Var[\sqrt n\Delta(\gamma_0)\mid \mtx Z,\vct\varepsilon]
\stackrel{P}{\longrightarrow}
\nu_2,
\]
where
\[
\nu_2
=
2\sigma_0^4
\left\{
\frac{h_2-h_1^2}{h_1^2}
-
(\tau+1)
\left(h_1-\frac{h_2}{h_1}\right)^2
\right\}.
\]

Third, the conditional mean has the same deterministic reduction:
by Lemma \ref{lem:noise_diag_approx_fine},
\[
\widetilde\Delta_*(\gamma_0)
=
\frac{1}{n}
\left(\vct\varepsilon^\top\vct\varepsilon-n\sigma_0^2\right)
\trace(\mtx B_{\gamma_0})
+
o_P(n^{-1/2}).
\]
Moreover,
\[
\trace(\mtx B_{\gamma_0})
\stackrel{P}{\longrightarrow}
h_1-\frac{h_2}{h_1}.
\]
It remains to identify the limiting law of
$\vct\varepsilon^\top\vct\varepsilon$. Write
$\vct\varepsilon=\mtx\Sigma_\varepsilon^{1/2}\vct g$, where
$\vct g\sim N(\vct 0,\mtx I_n)$. If
$\lambda_1,\ldots,\lambda_n$ are the eigenvalues of
$\mtx\Sigma_\varepsilon$, then
\[
\vct\varepsilon^\top\vct\varepsilon-\trace(\mtx\Sigma_\varepsilon)
=
\sum_{r=1}^n\lambda_r(g_r^2-1).
\]
Since $\|\mtx\Sigma_\varepsilon\|=O(1)$ and
\[
\sum_{r=1}^n\lambda_r^2
=
\|\mtx\Sigma_\varepsilon\|_F^2
=
n\kappa_\Sigma\sigma_0^4,
\]
the Lindeberg condition holds for the triangular array
$\{\lambda_r(g_r^2-1)\}_{r=1}^n$. Hence
\[
\frac{1}{\sqrt n}
\left(
\vct\varepsilon^\top\vct\varepsilon-n\sigma_0^2
\right)
\Longrightarrow
N(0,2\kappa_\Sigma\sigma_0^4).
\]
Therefore
\[
\sqrt n\,\widetilde\Delta_*(\gamma_0)
\Longrightarrow
N\left(
0,
2\kappa_\Sigma\sigma_0^4
\left(h_1-\frac{h_2}{h_1}\right)^2
\right).
\]
Denote this variance by $\nu_1^\Sigma$.

The same conditioning argument used in the proof of
Theorem \ref{thm:asymp_normality} now gives the joint convergence
\[
\sqrt n\Delta(\gamma_0)
\Longrightarrow
N(0,\nu_1^\Sigma+\nu_2).
\]
Finally, Lemma \ref{lem:Taylor_approx} gives
\[
\sqrt n(\hat\gamma-\gamma_0)
=
-
\frac{\sqrt n\Delta(\gamma_0)}
{\Delta_\infty'(\gamma_0)}
+o_P(1),
\]
where
\[
\Delta_\infty'(\gamma_0)
=
\frac{\sigma_0^2}{\gamma_0}
\frac{h_1^2-h_2}{h_1}.
\]
Let $D=h_2-h_1^2>0$. Since
\[
\nu_1^\Sigma+\nu_2
=
2\sigma_0^4
\left\{
\frac{D}{h_1^2}
+
(\kappa_\Sigma-\tau-1)\frac{D^2}{h_1^2}
\right\},
\]
we obtain
\[
\sqrt n(\hat\gamma-\gamma_0)
\Longrightarrow
N\left(
0,
2\gamma_0^2
\left\{
\frac{1}{h_2-h_1^2}
+
\kappa_\Sigma-\tau-1
\right\}
\right).
\]
This proves Theorem \ref{thm:asymp_normality_correlated}.

\subsection{Proof of Proposition \ref{thm:consistency_kappa}}
Straightforward calculation gives
\begin{align}
\label{eq:y_4_moment}
\E[y_i^4] &= \sum_{j=1}^p \left(\E[z_{ij}^4] - 3\right)\beta_j^4 + 3\|\vct{\beta}\|_2^4 + 6\|\vct{\beta}\|_2^2 \sigma_i^2 + \E[\varepsilon_i^4].
\end{align}
Therefore,
\begin{align*}
    \frac{1}{n\sigma_0^4}\sum_{i=1}^n\E[\varepsilon_i^4]
    =
    \frac{1}{n\sigma_0^4}\sum_{i=1}^n \E[y_i^4]
    -3\gamma_0^2-6\gamma_0
    -\frac{1}{n\sigma_0^4}\sum_{i=1}^n\sum_{j=1}^p
    \left(\E[z_{ij}^4]-3\right)\beta_j^4.
\end{align*}
By the assumption $\| \vct{\beta} \|_\infty = o(p^{-1/4})$ and $z_{ij}$ is sub-Gaussian, it is obvious that
\begin{align*}
    \frac{1}{n\sigma_0^4} \sum_{i=1}^n \sum_{j=1}^p \left(\E[z_{ij}^4] - 3\right)\beta_j^4 \leq \max_{i \in [n]; j \in [p]} \left\vert \E[z_{ij}^4] - 3  \right\vert \frac{1}{\sigma_0^4}\sum_{j=1}^p \beta_j^4 = o(1).
\end{align*}
Under the independent heteroscedastic noise setting considered in this proposition, $y_1,\ldots,y_n$ are independent because the rows of $\mtx{Z}$ are independent and the noise coordinates are independent. Furthermore, similar to \eqref{eq:y_4_moment}, the uniform sub-Gaussian bound on the design entries, $\|\vct\beta\|_2=O(1)$, and the additional condition $\max_i\E|\varepsilon_i|^8=O(1)$ imply $\max\limits_{1 \leq i \leq n} \E[y_i^8] = O(1)$. Therefore,
\begin{align*}
  \Var\left[ \frac{1}{n} \sum_{i=1}^n y_i^4 \right] = \frac{1}{n^2} \sum_{i=1}^n \Var[y_i^4] \leq \frac{1}{n^2} \sum_{i=1}^n \E[y_i^8] = O\left(\frac{1}{n}\right).
\end{align*}
Then we can have $\frac{1}{n}\sum_{i=1}^n \E[y_i^4]=\frac{1}{n} \sum_{i=1}^n y_i^4 + O_P(n^{-1/2})$. Combining the above, we have
\begin{align*}
\frac{1}{n \sigma_0^4} \sum_{i=1}^n y_i^4 -3\gamma_0^2-6\gamma_0
\stackrel{P}{\longrightarrow}
\frac{1}{n\sigma_0^4}\sum_{i=1}^n\E[\varepsilon_i^4].
\end{align*}
In Theorem \ref{thm:consistency} we have already shown that $\hat{\sigma}^2 \stackrel{P}{\longrightarrow} \sigma_0^2$ and $\hat{\gamma}_n  \stackrel{P}{\longrightarrow} \gamma_0$. Thus, by Slutsky's theorem,
\[
\frac{1}{n\hat{\sigma}^4}\sum_{i=1}^n y_i^4
-3\hat{\gamma}_n^2-6\hat{\gamma}_n
\stackrel{P}{\longrightarrow}
\frac{1}{n\sigma_0^4}\sum_{i=1}^n\E[\varepsilon_i^4].
\]
If the noise is Gaussian, then $\E[\varepsilon_i^4]=3\sigma_i^4$ and
\[
\frac{1}{n\sigma_0^4}\sum_{i=1}^n\E[\varepsilon_i^4]
=\frac{3}{n\sigma_0^4}\sum_{i=1}^n\sigma_i^4
=3\kappa_\varepsilon,
\]
and hence
$\hat{\kappa}_{\varepsilon,\mathrm{G}}\stackrel{P}{\longrightarrow}\kappa_\varepsilon$.
If the noise is homogeneous, then $\sigma_i^2=\sigma_0^2$ for all $i$, so
\[
\kappa_\varepsilon
=\frac{1}{2n\sigma_0^4}\sum_{i=1}^n
\left\{\E[\varepsilon_i^4]-\sigma_0^4\right\}
=\frac{1}{2}\left\{
\frac{1}{n\sigma_0^4}\sum_{i=1}^n\E[\varepsilon_i^4]-1
\right\},
\]
and therefore $\hat{\kappa}_{\varepsilon,\mathrm{H}}\stackrel{P}{\longrightarrow}\kappa_\varepsilon$.
This proves the proposition.

\subsection{Proof of Corollary \ref{cor:ci_gamma}}
Since $\tau_n=n/p \to \tau$, continuity of $h_1(\gamma,\tau)$ and $h_2(\gamma,\tau)$ on $(0,\infty)^2$ implies
\[
\mathcal{V}(\hat{\gamma},\bar\kappa_n,\tau_n)
\stackrel{P}{\longrightarrow}
\mathcal{V}(\gamma_0,\kappa_\varepsilon,\tau).
\]
Therefore
\[
n\hat{s}_n^2=\mathcal{V}(\hat{\gamma},\bar\kappa_n,\tau_n)
\stackrel{P}{\longrightarrow}
\mathcal{V}(\gamma_0,\kappa_\varepsilon,\tau),
\]
or equivalently $\hat{s}_n^2=n^{-1}\mathcal{V}(\gamma_0,\kappa_\varepsilon,\tau)+o_P(n^{-1})$. Combining this with Theorem \ref{thm:asymp_normality}, Slutsky's theorem yields
\[
\frac{\hat{\gamma}-\gamma_0}{\hat{s}_n}
\Longrightarrow \mathcal{N}(0,1).
\]
Hence,
\[
\P\left(\left|\frac{\hat{\gamma}-\gamma_0}{\hat{s}_n}\right|\le z_{1-\alpha/2}\right)\longrightarrow 1-\alpha,
\]
which is equivalent to the asserted asymptotic coverage of \eqref{eq:ci_gamma}.

\section{Discussion}
\label{sec:discussion}
This paper studies Gaussian random-effects MLEs for SNR estimation under coefficient- and noise-model misspecification. We establish consistency and asymptotic normality for fixed dense coefficients and independent finite-moment noise, with a parallel correlated Gaussian benchmark. The limiting variance identifies a scalar noise-square fluctuation parameter, leading to plug-in confidence intervals for heterogeneous Gaussian noise and homogeneous non-Gaussian noise. The simulations show that this calibration is important when the homogeneous Gaussian approximation understates uncertainty.

The symmetry assumption in the independent-noise theory is mainly tied to the proof strategy. We use row and column Rademacher sign flips to condition the estimating function $\Delta(\gamma)$, reveal the variance contributions of heterogeneity and correlation, and apply normal approximation tools for Rademacher quadratic forms \citep{Chat2008}. The variance calculation also extends the ``leave-$k$-column-out'' argument of \cite{jiang2016high} to both rows and columns.

Several questions remain. Relaxing the symmetry condition would broaden applicability to skewed covariates; the supplementary design-distribution checks in Appendix \ref{app:additional_simulation_figures} give preliminary empirical evidence in this direction. It would also be useful to extend the analysis to multivariate outcomes and grouped features, as in heritability components or group ridge regression \citep{yang2011gcta,ignatiadis2020group}. In these settings, estimating standard errors under unknown dependence remains a central challenge.

\section*{Acknowledgment}
X. Li and X. Hu are partially supported by the NSF via the Career Award DMS-1848575. We would like to thank Debashis Paul for inspiring discussions.

\newpage
\begin{appendices}
\section{Preliminaries}
Let's first recall the famous Mar$\check{c}$enko-Pastur law in random matrix theory.
\begin{theorem}[Mar$\check{c}$enko-Pastur law, \cite{M-PLaw}]
\label{thm:M-P_Law} 
Let $\mtx{Z}$ be an $n \times p$ random matrix whose entries are i.i.d. random variables with mean $0$ and variance $1$ in which $n/p \rightarrow \tau \in (0,\infty)$ as $n,p \rightarrow \infty$. Then the empirical spectral distribution (ESD) of $S = p^{-1} \mtx{Z} \mtx{Z}^\top$, which is defined as $F^{\mtx{S}}$, converges almost surely (a.s.) in distribution to $F_{\tau}$, whose p.d.f. is given by
\begin{align*}
    f_{\tau}(x) = \left\{ \begin{aligned}
        &\max\{\tau - 1, 0\}\delta_0(x) + \frac{1}{2 \pi \tau x} \sqrt{\left( b_{+}(\tau) - x\right) \left(x - b_{-}(\tau) \right) }& & b_{-}(\tau) \leq x \leq b_{+}(\tau)\\
        &0 & & elsewhere
    \end{aligned}
    \right.
\end{align*}
where $b_{\pm}(\tau) = (1 \pm \sqrt{\tau})^2$ and $\delta_0(x)$ is a point mass $\tau^{-1}$ at the origin.
\end{theorem}
Note that in our settings, the entries of the design matrix are not necessarily identically distributed. To this end, we consider the following extension of Mar$\check{c}$enko-Pastur law.
\begin{theorem}[\cite{Bai1999}, Theorem 2.8]
\label{thm:MP_not_iid}
Let $\mtx{Z}$ be an $n \times p$ random matrix whose entries are independent random variables with mean $0$ and variance $1$. Assume that $n/p \rightarrow \tau \in (0,\infty)$ and that for any $\delta > 0$,
\[
\frac{1}{\delta^2 np} \sum_{i,j} \E\left[ | z_{ij}^{(n)} |^2 I_{(| z_{ij}^{(n)} | \geq \delta \sqrt{n} )}\right] \rightarrow 0.
\]
Then $F^{\mtx{S}}$, defined as in Theorem \ref{thm:M-P_Law}, tends almost surely to the Mar$\check{c}$enko-Pastur law with ratio index $\tau$.
\end{theorem}

\begin{corollary} \label{thm:MP_trace} Under the assumption of Theorem \ref{thm:M-P_Law} or \ref{thm:MP_not_iid}, for any integer $l$, we have 
\[
\frac{1}{n} \trace(\mtx{S}^l)  \stackrel{a.s.}{\longrightarrow} \int_{b_{-}(\tau)}^{b_{+}(\tau)} x^l f_{\tau}(x)  \mathrm{d}x \quad \text{as}\quad n,p \rightarrow \infty.
\]

\end{corollary}

Define the sub-Gaussian norm of a random variable $\zeta$ as
$$
\|\zeta\|_{\psi_2} \equiv \sup\limits_{q \geq 1} \{ q^{-1/2} (\mathbb{E}|\zeta|^q)^{1/q}\}. 
$$
A random variable $\zeta$ is sub-Gaussian if and only if its sub-Gaussian norm $\|\zeta\|_{\psi_2} < \infty$. We have the following equivalent characterizations on the sub-Gaussianity of a random variable:
\begin{lemma}[\cite{vershynin2010introduction}, Lemma 5.5] \label{sub-Gaussian}
A random variable $\zeta$ is sub-Gaussian if and only if
\begin{itemize}
\item [1)]  $\|\zeta\|_{\psi_2} < \infty$;\quad \text{or}
\item [2)] $\P\{\vert \zeta \vert > t\} \leq \exp(1 - t^2/K^2)$ for some parameter $K>0$ and all $t>0$.
\end{itemize}
\end{lemma}

Part 2) implies that the design matrix under the setting of Theorem \ref{thm:consistency}, in which the entries have sub-Gaussian norms that are uniformly upper bounded, satisfies the conditions in Theorem \ref{thm:MP_not_iid}. Indeed, if $\zeta$ is a sub-Gaussian random variable, then by the identity $\E[X] = \int_{0}^\infty \P(X > t) \mathrm{d}t$ for any nonnegative random variable $X$, we have
\begin{align*}
    \E\left[ | \zeta |^2 I_{(| \zeta | \geq \delta \sqrt{n} )}\right] 
    &= \int_{\delta \sqrt{n}}^{\infty} \P\{\vert \zeta \vert > t\} 2 t \mathrm{d}t + \delta^2 n\P\{\vert \zeta \vert > \delta \sqrt{n}\}\\
    &\leq 2\int_{\delta \sqrt{n}}^{\infty} e^{1 - \frac{t^2}{K^2}} t \mathrm{d}t + \delta^2 n e^{1 - \frac{\delta^2 n}{K^2}}\\
    & = (K^2 + \delta^2 n) e^{1 - \frac{\delta^2 n}{K^2}}.
\end{align*}
This implies that for $n \times p$ random matrices $\mtx{Z}$ whose entries have uniformly upper bounded sub-Gaussian norms,
\[
\frac{1}{\delta^2 np} \sum_{i,j} \E\left[ | z_{ij}^{(n)} |^2 I_{(| z_{ij}^{(n)} | \geq \delta \sqrt{n} )}\right] \rightarrow 0,
\]
as $n, p \rightarrow \infty$, for any $\delta > 0$.

Our proof also relies crucially on the following fundamental concentration inequalities.
\begin{proposition}[Hanson–Wright inequality, \cite{10.1214/ECP.v18-2865}]\label{ineq:Hanson-Wright}
Let $\vct{\zeta} = (\zeta_1, \cdots, \zeta_n)^\top$, where the $\zeta_i$'s are independent random variables satisfying $\mathbb{E}(\zeta_i) = 0$ and $\|\zeta_i\|_{\psi_2} \leq K < \infty$. Let $\mtx{A}$ be an $n \times n$ deterministic matrix. Then we have for any $t>0$,
$$
\P \{ |\vct{\zeta}^\top \mtx{A} \vct{\zeta} - \mathbb{E}(\vct{\zeta}^\top \mtx{A} \vct{\zeta})| > t  \} \leq 2 \exp \left\{ -c \min \left( \frac{t^2}{K^4 \|\mtx{A}\|_F^2}, \frac{t}{K^2 \|\mtx{A}\|}\right) \right\},
$$
where $c>0$ is an absolute constant. Here $\|\mtx{A}\|$ and $\|\mtx{A}\|_F$ denote the operator and Frobenius norms of $\mtx{A}$, respectively.
\end{proposition}

\begin{proposition}[Hoeffding-type inequality, \cite{vershynin2010introduction}, Proposition 5.10]\label{ineq:Hoeffding}
Let $\vct{\zeta} = (\zeta_1, \cdots, \zeta_n)^\top$, where the $\zeta_i$'s are independent centered sub-Gaussian random variables. Let $K = \max_{1\leq i \leq n} \|\zeta_i\|_{\psi_2} $ and $\vct{a} = (a_1,\cdots,a_N)^\top \in \mathbb{R}^N$. Then we have for any $t \geq 0$,
$$
\P \{ |\vct{a}^\top \vct{\zeta}| > t  \} \leq e \exp \left\{ -c \frac{t^2}{K^2 \| \vct{a} \|_2^2}\right\},
$$
where $c>0$ is an absolute constant.
\end{proposition}

\begin{proposition}[Bernstein-type inequality, \cite{vershynin2010introduction}, Proposition 5.16] \label{ineq:Bernstein}
Let $\vct{\zeta} = (\zeta_1, \cdots, \zeta_n)^\top$, where the $\zeta_i$'s are independent centered sub-exponential random variables. Let $K = \max_{1\leq i \leq n} \|\zeta_i\|_{\psi_2} $ and $\vct{a} = (a_1,\cdots,a_N)^\top \in \mathbb{R}^N$. Then we have for any $t \geq 0$,
$$
\P \{ |\vct{a}^\top \vct{\zeta}| > t  \} \leq 2 \exp \left\{ -c \min \left( \frac{t^2}{K^2 \| \vct{a} \|_2^2}, \frac{t}{K \| \vct{a} \|_\infty} \right) \right\},
$$
where $c>0$ is an absolute constant.
\end{proposition}

The next result, the famous Sherman-Morrison-Woodbury formula in matrix analysis is repeatedly used in our proofs, as the corner stone of leave-one-out analysis.
\begin{theorem}[Sherman-Morrison-Woodbury formula, \cite{matrix1990}, Page 19]
\label{thm:matrix_inv}
	Let $\mtx{P}$ and $\mtx{Q}$ be n-dimensional non-singular matrices such that $\mtx{Q} = \mtx{P}+\mtx{U}\mtx{V}^\top$, where $\mtx{U}, \mtx{V} \in \mathbb{R}^{n \times q}$. Then
	$$
	\mtx{Q}^{-1}  = (\mtx{P} + \mtx{U}\mtx{V}^\top)^{-1} = \mtx{P}^{-1} - \mtx{P}^{-1} \mtx{U} (\mtx{I}_q + \mtx{V}^\top \mtx{P}^{-1}\mtx{U})^{-1} \mtx{V}^\top \mtx{P}^{-1}.
	$$
\end{theorem}

The following results, implied by \cite{Chat2008} and \cite{chatterjee2009fluctuations}, are conditions for the normality of quadratic forms.
\begin{theorem}[\cite{Chat2008}, Proposition 3.1]\label{thm:normality}
Let $X = (X_1,\ldots,X_n)$ be i.i.d. Rademacher random variables and $A = (a_{ij} )1\leq i,j\leq n$ be a real symmetric matrix. Let $W = X^\top \mtx{A} X$ and 
\[
\sigma^2 = \Var(W) = \frac{1}{2} \trace(\mtx{A}^2).
\]
Let $\mu$ be the law of $(W-\E(W))/\sqrt{\Var(W)}$ and let $\nu$ be the standard Gaussian law. We define
$$d_W \coloneqq \mathcal{W}(\mu,\nu),$$
where $\mathcal{W}$ is the Kantorovich–Wasserstein distance between two probability measures with
\[
\mathcal{W}(\mu,\nu) = \sup \left\{\left|\int h d\mu - \int h d \nu \right|: h \text{ Lipschitz}, \text{ with } \|h\|_{\text{Lip}} \leq 1  \right\}
\]
Then,
\[
d_W \leq \left(\frac{\trace(\mtx{A}^4)}{2 \sigma^4} \right)^{1/2} + \frac{5}{2 \sigma^3} \sum_{i=1}^n\left(\sum_{j=1}^n a_{ij}^2 \right)^{3/2} \leq 6\sqrt{2}\frac{\|\mtx{A}\|^2}{\|\mtx{A}\|_F^2}.  
\]
\end{theorem}

\begin{theorem}[\cite{chatterjee2009fluctuations}] \label{thm: nor_xtx}
Suppose $\vct{x}$ is a gaussian random vector with mean 0 and covariance matrix $\mtx{\Sigma}$. Take any $g \in C^2(\mathbb{R})$ and let $\nabla g$ and $\nabla^2 g$  denote the gradient and Hessian of $g$. Let 
\begin{align*}
    \varsigma_1  = \left( \E \| \nabla g(\vct{x}) \|^4 \right)^\frac{1}{4}, \quad  \varsigma_2  = \left( \E \| \nabla^2 g(\vct{x}) \|^4 \right)^\frac{1}{4}.
\end{align*}
Then let $W =  g(\vct{x})$ have a finite fourth moment and $U$ be a normal random variable having the same mean and variance as $W$,
\begin{align*}
    d_{TV}(W,U) \leq \frac{2 \sqrt{5} \|\mtx{\Sigma} \|^{\frac{3}{2}}\varsigma_1 \varsigma_2 }{\Var \left[ W\right]}.
\end{align*}
Here $d_{TV}$ is the total variation distance between random variables $u$ and $v$, 
\begin{align*}
    d_{TV}(u,v) = \sup_{B \in \mathcal{B}(\mathbb{R})} | \P(u \in B) - \P (v \in B)|, 
\end{align*}
where $\mathcal{B}(\mathbb{R})$ denotes the collection of Borel sets in $\mathbb{R}$.
\end{theorem}

Next, there is a famous result for the bounds of eigenvalues of the sub-gaussian random matrix.
\begin{theorem}[Theorem 5.39, \cite{vershynin2010introduction}]\label{thm:eig_subG}
Let $\mtx{Z}$ be an $n\times p$ matrix whose rows are independent sub-gaussian isotropic random vectors. Then for every $t \geq 0$, with probability at least $1- 2\exp(-ct^2)$ one has
\[
\sqrt{n} - C\sqrt{p} - t \leq \lambda_{\min} (\mtx{Z}) \leq \lambda_{\max} (\mtx{Z}) \leq \sqrt{n} + C\sqrt{p} + t
\]
Here $C = C_K$, $c = c_K > 0$ depend only on the subgaussian norm $K$ of the rows.
\end{theorem}

\section{Proofs of Lemmas in Section \ref{sec:proofs}}
\label{appdenixB}

In this section we give detailed proofs of the technical lemmas that appear in Section \ref{sec:proofs}. As mentioned earlier, the proofs of Lemmas \ref{lem:leave_one_out}, \ref{lem:diag_homo}, \ref{lem:diag_homo_taylor}, \ref{lem:diag_homo_moment}, \ref{lem:leave_two_out}, and \ref{lem:Delta_starstar} follow the proof ideas in \cite{jiang2016high}, with modifications for the present setting. We provide self-contained proofs here for completeness.

\subsection{Proof of Lemma \ref{lem:noise_var_bound}}
Under the assumptions of Theorem \ref{thm:consistency}, Markov's inequality gives, for any $M>0$,
\[
\P\left(\max_{i\in[n]}\varepsilon_i^2>Mn^{2/(4+\delta)}\right)
\leq
\frac{nC_\varepsilon}{M^{(4+\delta)/2}n}
=C_\varepsilon M^{-(4+\delta)/2},
\]
which proves \eqref{eq:max_esp_i} by taking $M\to\infty$. Moreover,
\[
\Var\left(\frac{1}{n}\sum_{i=1}^n(\varepsilon_i^2-\sigma_i^2)\right)
\leq \frac{1}{n^2}\sum_{i=1}^n \E \varepsilon_i^4
=O(n^{-1}),
\]
and hence
\[
\frac{1}{n}\sum_{i=1}^n(\varepsilon_i^2-\sigma_i^2)=o_P(1),
\]
which proves \eqref{eq:noise_var_concent}. Under the assumptions of Theorem \ref{thm:asymp_normality}, set
\[
s_n^2=\sum_{i=1}^n\Var(\varepsilon_i^2)=2n\kappa_\varepsilon\sigma_0^4.
\]
If $\kappa_\varepsilon>0$, then $s_n^2\asymp n$. The moment assumption gives
$\max_i\E|\varepsilon_i^2-\sigma_i^2|^{2+\delta/2}=O(1)$, and hence, for any fixed $\eta>0$,
\[
\frac{1}{s_n^2}\sum_{i=1}^n
\E\left[(\varepsilon_i^2-\sigma_i^2)^2
1_{\{|\varepsilon_i^2-\sigma_i^2|>\eta s_n\}}\right]
\leq
\frac{1}{\eta^{\delta/2}s_n^{2+\delta/2}}
\sum_{i=1}^n \E|\varepsilon_i^2-\sigma_i^2|^{2+\delta/2}
=o(1).
\]
Thus the Lindeberg condition for the triangular array
$\{\varepsilon_i^2-\sigma_i^2\}_{i=1}^n$ holds, and the Lindeberg-Feller central limit theorem yields \eqref{eq:noise_CLT}. If $\kappa_\varepsilon=0$, then $s_n^2=0$ and the left-hand side of \eqref{eq:noise_CLT} is degenerate, so the same display holds with zero limiting variance.

It remains to verify the corresponding statements under Theorem \ref{thm:consistency_correlated} and Theorem \ref{thm:asymp_normality_correlated}. The bound $\max_{i\in[n]}\varepsilon_i^2=O_P(\log n)$ follows from a union bound applied to the marginal Gaussian tails and $\max_i\sigma_i^2=O(1)$. Also,
\[
\E\left[\vct{\varepsilon}^\top\vct{\varepsilon}\right]=\trace(\mtx{\Sigma}_\varepsilon)=n\sigma_0^2,\qquad
\Var\left(\vct{\varepsilon}^\top\vct{\varepsilon}\right)=2\|\mtx{\Sigma}_\varepsilon\|_F^2,
\]
so \eqref{eq:noise_var_concent} follows from $\|\mtx{\Sigma}_\varepsilon\|_F=o(n)$. Finally, applying Theorem \ref{thm: nor_xtx} to
\[
g(\vct{x})=\frac{1}{\sqrt n}\sum_{i=1}^n x_i^2
\]
and using $\|\mtx{\Sigma}_\varepsilon\|=O(1)$ and
$\|\mtx{\Sigma}_\varepsilon\|_F^2=n\kappa_\Sigma\sigma_0^4$ gives
\[
\frac{1}{\sqrt{n}}\left(\vct{\varepsilon}^\top\vct{\varepsilon}-n\sigma_0^2\right)
\Longrightarrow
\mathcal{N}(0,2\kappa_\Sigma\sigma_0^4),
\]
which is the Gaussian correlated version of \eqref{eq:noise_CLT}.

\subsection{Proof of Lemma \ref{lem:leave_one_out}}
For convenience, define
\begin{equation}
\label{eq:rho_phi_psi}
\begin{cases}
\rho_k \coloneqq \eta_{kk,k}^{(1)} \coloneqq \vct{z}_k^\top \mtx{V}_{\gamma, -k}^{-1} \vct{z}_k, \\
\phi_k \coloneqq \eta_{kk,k}^{(2)} \coloneqq \vct{z}_k^\top \mtx{V}_{\gamma, -k}^{-2} \vct{z}_k, \\
\psi_k \coloneqq \eta_{kk,k}^{(3)} \coloneqq \vct{z}_k^\top \mtx{V}_{\gamma, -k}^{-3} \vct{z}_k.
\end{cases}
\end{equation}
First, there is a simple relationship: $\psi_k \leq \phi_k \leq \rho_k$. Indeed, since $\mtx{I}_n-\mtx{V}_{\gamma, -k}^{-1} \succeq \mtx{0}$, we know that
\[
\rho_k - \phi_k = \vct{z}_k^\top \mtx{V}_{\gamma, -k}^{-1/2}(I-\mtx{V}_{\gamma, -k}^{-1})\mtx{V}_{\gamma, -k}^{-1/2}\vct{z}_k \geq 0.
\]
i.e., $\phi_k \leq \rho_k$. We can similarly obtain $\psi_k \leq \phi_k$.

Using Sherman-Morrison-Woodbury formula (Theorem \ref{thm:matrix_inv}), we have
\begin{equation}
\label{eq:V_inv_sherman}
\mtx{V}_\gamma^{-1} = \mtx{V}_{\gamma, -k}^{-1} - \frac{\gamma}{p}(1+ \frac{\gamma}{p}\rho_k)^{-1} \mtx{V}_{\gamma, -k}^{-1} \vct{z}_k \vct{z}_k^\top \mtx{V}_{\gamma, -k}^{-1},
\end{equation}
and
\begin{align}
\mtx{V}_\gamma^{-2} &= \left(\mtx{V}_{\gamma, -k}^{-1} - \frac{\gamma}{p}(1+ \frac{\gamma}{p}\rho_k)^{-1} \mtx{V}_{\gamma, -k}^{-1} \vct{z}_k \vct{z}_k^\top \mtx{V}_{\gamma, -k}^{-1}\right)^2    \nonumber
\\
&= \mtx{V}_{\gamma, -k}^{-2} 
- \frac{\gamma}{p}(1+ \frac{\gamma}{p}\rho_k)^{-1}\mtx{V}_{\gamma, -k}^{-2} \vct{z}_k \vct{z}_k^\top \mtx{V}_{\gamma, -k}^{-1} 
- \frac{\gamma}{p}(1+ \frac{\gamma}{p}\rho_k)^{-1}\mtx{V}_{\gamma, -k}^{-1} \vct{z}_k \vct{z}_k^\top \mtx{V}_{\gamma, -k}^{-2} \nonumber
\\
&~~+ \left(\frac{\gamma}{p}\right)^2(1+ \frac{\gamma}{p}\rho_k)^{-2} \phi_k \mtx{V}_{\gamma, -k}^{-1} \vct{z}_k \vct{z}_k^\top \mtx{V}_{\gamma, -k}^{-1}. 
\label{eq:V_inv_sq_sherman}
\end{align}
By \eqref{eq:V_inv_sherman} and \eqref{eq:V_inv_sq_sherman}, we can also have 
\[
\trace(\mtx{V}_\gamma^{-1}) = \trace(\mtx{V}_{\gamma, -k}^{-1}) - \frac{\gamma}{p}(1+ \frac{\gamma}{p}\rho_k)^{-1} \phi_k,
\]
and 
\[
\trace(\mtx{V}_\gamma^{-2}) = \trace(\mtx{V}_{\gamma, -k}^{-2}) - \frac{2 \gamma}{p}( 1 + \frac{\gamma}{p}\rho_k)^{-1} \psi_k + (\frac{\gamma}{p})^2(1+ \frac{\gamma}{p}\rho_k)^{-2}\phi_k^2.
\]
Then
\begin{align}
\left\vert \trace(\mtx{V}_\gamma^{-1}) - \trace(\mtx{V}_{\gamma, -k}^{-1}) \right\vert = \frac{\gamma}{p}(1+ \frac{\gamma}{p}\rho_k)^{-1} \phi_k \leq \frac{\gamma}{p}(1+ \frac{\gamma}{p}\rho_k)^{-1} \rho_k <1, \label{eq:diff_Vgamma_Vgamma_k_1}
\end{align}
and
\begin{align}
\left\vert \trace(\mtx{V}_\gamma^{-2}) - \trace(\mtx{V}_{\gamma, -k}^{-2}) \right\vert \leq \frac{2 \gamma}{p}( 1 + \frac{\gamma}{p}\rho_k)^{-1} \rho_k + (\frac{\gamma}{p})^2(1+ \frac{\gamma}{p}\rho_k)^{-2}\rho_k^2 < 3 \label{eq:diff_Vgamma_Vgamma_k_2}.
\end{align}
Similarly, we can also prove that
\[
\left\vert \trace(\mtx{V}_\gamma^{-3}) - \trace(\mtx{V}_{\gamma, -k}^{-3}) \right\vert \leq 7 \text{ and } \left\vert \trace(\mtx{V}_\gamma^{-4}) - \trace(\mtx{V}_{\gamma, -k}^{-4}) \right\vert \leq 15.
\]

Since the entries of $\mtx{Z}$ are independent sub-Gaussian and $\mathbb{E}(z_{ik}) = 0$, using Proposition \ref{ineq:Hanson-Wright}, we have, for any $1\leq k \leq p$ and $t> 0$:
\begin{align*}
    \P{\left\{ |\rho_k - \trace (\mtx{V}_{\gamma, -k}^{-1})| > t \vert \mtx{V}_{\gamma, -k}  \right\}} \leq 2\exp\left\{ -c \min \left( \frac{t^2}{K^4 \|\mtx{V}_{\gamma, -k}^{-1} \|_F^2}, \frac{t}{K^2\|\mtx{V}_{\gamma, -k}^{-1} \|} \right) \right\},
\end{align*}
where $c$ and $K$ are positive constants. If we set
\begin{align*}
    t = t_{k} = K^2 \max \left( \sqrt{\frac{2 \log p  }{c}}\|\mtx{V}_{\gamma, -k}^{-1} \|_F, \frac{2 \log p  }{c} \|\mtx{V}_{\gamma, -k}^{-1} \|  \ \right),
\end{align*}
it follows that $
    \P{ \left\{ |\rho_k - \trace (\mtx{V}_{\gamma, -k}^{-1}) | >  t_k |\mtx{V}_{\gamma, -k}^{-1}    \right\} }\leq 2/p^2.
$ Thus
\begin{align*}
    \P{ \left\{\max_{1\leq k \leq p} t_{k}^{-1} |\rho_k - \trace(\mtx{V}_{\gamma, -k}^{-1}) | >  1 \right\}} \leq  \frac{2}{p}.
\end{align*}
By Lemma \ref{thm:diff_tr_h}, $\|\mtx{V}_{\gamma, -k}^{-1}\| \leq 1$, and $\|\mtx{V}_{\gamma, -k}^{-1}\|_F \leq \sqrt{n} \|\mtx{V}_{\gamma, -k}^{-1}\| \leq \sqrt{n}$, we can obtain that
$$
t_{k} \leq K^2 \max\left( \sqrt{\frac{2}{c}} \sqrt{n \log p}, \frac{2}{c}\log p \right),
$$
which implies 
\[
\P{ \left\{\max_{1\leq k \leq p}  |\rho_k - \trace(\mtx{V}_{\gamma, -k}^{-1}) | >  C \sqrt{n \log p} \right\}} \leq 2/p 
\]
for some constant $C > 0$. Then, it follows that
\begin{align}
    \max_{1\leq k \leq p} |\rho_k -\trace( \mtx{V}_{\gamma, -k}^{-1} )| & = O_P(\sqrt{n \log n} ) \label{eq:diff_rhok_Vgamma_k}.
\end{align}
By a similar argument, we have
\begin{align}
\max_{1\leq k \leq p} |\phi_k - \trace(\mtx{V}_{\gamma, -k}^{-2} )| & = O_P(\sqrt{n \log n} ) \label{eq:diff_phik_Vgamma_k}.
\end{align}
Combining \eqref{eq:diff_Vgamma_Vgamma_k_1}, \eqref{eq:diff_rhok_Vgamma_k},  \eqref{eq:diff_Vgamma_Vgamma_k_2} and \eqref{eq:diff_phik_Vgamma_k}, we have
\begin{align}
    \max_{1\leq k \leq p} |\rho_k -\trace( \mtx{V}_{\gamma}^{-1} )| & = O_P(\sqrt{n \log n} ), \quad \text{and}\label{eq:diff_rhok_Vgamma}  \\
    \max_{1\leq k \leq p} |\phi_k - \trace(\mtx{V}_{\gamma}^{-2} )| & = O_P(\sqrt{n \log n} ) \label{eq:diff_phik_Vgamma}.   
\end{align}

\subsection{Proof of Lemma \ref{lem:diag_homo}}
Based on \eqref{eq:V_inv_sherman} and \eqref{eq:V_inv_sq_sherman}, there holds
\begin{align}
\vct{z}_k^\top \mtx{V}_\gamma^{-1} \vct{z}_k = (1+ \frac{\gamma}{p} \rho_k)^{-1} \rho_k,
\label{eq:mean_A1_rhok}
\\
\vct{z}_k^\top \mtx{V}_\gamma^{-2} \vct{z}_k = (1+ \frac{\gamma}{p} \rho_k)^{-2} \phi_k,
\label{eq:mean_A2_phik}
\end{align}
where $\rho_k$ and $\phi_k$ are defined in \eqref{eq:rho_phi_psi}. 

Let's now come back to find approximations of $\E[A_1 \vert \mtx{Z}]$ and $\E[A_2 \vert \mtx{Z}]$. We define the following intermediate quantities
\begin{align}
\label{eq:theta_12}
\theta_1 = \frac{1}{n}\frac{\trace(\mtx{V}_{\gamma}^{-1})}{1 +  \frac{\gamma}{p} \trace(V_{\gamma}^{-1}) } 
\text{~~and~~} 
\theta_2 = \frac{1}{n}\frac{\trace(\mtx{V}_{\gamma}^{-2})}{ \left(1 +  \frac{\gamma}{p} \trace(V_{\gamma}^{-1}) \right)^2 }.
\end{align}
Then by \eqref{eq:mean_A1_rhok} and \eqref{eq:diff_rhok_Vgamma}, we can have
\begin{align} 
\label{eq:diff_theta1_mean_A1}
\max_{1\leq k \leq p} \left|\theta_1 - \frac{\vct{z}_k^\top \mtx{V}_\gamma^{-1} \vct{z}_k}{n}  \right|
      \leq \max_{1\leq k \leq p} \left| \frac{\trace(\mtx{V}_{\gamma}^{-1}) - \rho_k}{n}  \right|= O_P\left(\sqrt{\frac{\log n}{n}}  \right),
\end{align}
which implies \eqref{eq:diag_homo_1}. Similarly, by \eqref{eq:mean_A2_phik} and \eqref{eq:diff_phik_Vgamma}, there holds
\begin{align}
\label{eq:diff_theta2_mean_A2_detail}
\max_{1\leq k \leq p} \left|\theta_2 - \frac{\vct{z}_k^\top \mtx{V}_\gamma^{-2} \vct{z}_k}{n}  \right| & \leq  \max_{1\leq k \leq p} \frac{1}{n} \left| \trace(\mtx{V}_{\gamma}^{-2}) \left(1 +  \frac{\gamma}{p} \rho_k\right)^2 - \left(1 +  \frac{\gamma}{p} \trace(V_{\gamma}^{-1}) \right)^2 \phi_k  \right| \nonumber
\notag\\
& \leq \frac{1}{n} \left[ \max_{1\leq k \leq p} \left|  \trace(\mtx{V}_{\gamma}^{-2}) - \phi_k\right| + \max_{1\leq k \leq p}  \frac{2 \gamma}{p} \rho_k \left| \trace(\mtx{V}_{\gamma}^{-2})  -  \phi_k \right|    \right. \nonumber
\notag\\
& \quad \left. + \max_{1\leq k \leq p}\frac{2 \gamma}{p} \phi_k  \left| \rho_k - \trace(\mtx{V}_{\gamma}^{-1}) \right|   + \max_{1\leq k \leq p} \frac{\gamma^2}{p^2} \rho_k^2 \left|  \trace(\mtx{V}_{\gamma}^{-2}) -  \phi_k \right|  \right.\nonumber
\notag\\
& \quad  \left.  + \max_{1\leq k \leq p}\frac{\gamma^2}{p^2}\phi_k \left(\left| \rho_k -  \trace(\mtx{V}_{\gamma}^{-1})  \right| \left| \rho_k + \trace(\mtx{V}_{\gamma}^{-1})  \right| \right) \right].
\end{align}
It follows, by the facts $\trace(\mtx{V}_{\gamma}^{-1}) = O_P(n)$, $\trace(\mtx{V}_{\gamma}^{-2}) = O_P(n)$, $\rho_k = O_P(n)$ and $\phi_k = O_P(n)$ (Lemma \ref{thm:diff_tr_h} and Lemma \ref{lem:leave_one_out}),\eqref{eq:diag_homo_2} is true.\\
Then we can have for $l = 1,2$,
\begin{align}\label{eq:thm_D1}
    &\left|\frac{\vct{z}_k^\top \mtx{V}_\gamma^{-l} \vct{z}_k}{n} - \frac{1}{np} \trace\left( \mtx{V}_{\gamma}^{-1} \mtx{Z} \mtx{Z}^\top  \right)  \right| \notag\\
    & = \left|\frac{\vct{z}_k^\top \mtx{V}_\gamma^{-l} \vct{z}_k}{n} - \theta_1 +  \frac{1}{p} \sum_{i=1}^p \left( \theta_1 - \sum_{i=1}^p \frac{1}{n}\vct{z}_i^\top \mtx{V}_\gamma^{-1} \vct{z}_i \right)\right| = O_P\left(\sqrt{\frac{\log n}{n}}  \right),
\end{align}
which implies \eqref{eq:diag_homo_3}. Then by the definition of $\mtx{B}_\gamma$, when $l=1$ we can get
\eqref{eq:diag_homo_4} from \eqref{eq:diag_homo_3}.

When $l=2$, for $\left(\vct{z}_k^\top \mtx{B}_{\gamma_0} \vct{z}_k \right)^2$,
\begin{align}\label{eq:z1_B_z1}
    \left(\vct{z}_k^\top \mtx{B}_{\gamma} \vct{z}_k \right)^2 &=  \left(\frac{1}{n}\vct{z}_k^\top\mtx{V}_\gamma^{-1}\vct{z}_k  - \frac{\frac{1}{n} \vct{z}_k^\top \mtx{V}_\gamma^{-2}\vct{z}_k }{\frac{1}{n}\trace(\mtx{V}_\gamma^{-1})}  \right)^2\notag \\
    & = \left(\frac{1}{n} \eta_{kk}^{(1)} \right)^2-2 \frac{ \frac{1}{n} \eta_{kk}^{(1)}\frac{1}{n} \eta_{kk}^{(2)}}{\frac{1}{n}\trace(\mtx{V}_\gamma^{-1}) }  + \frac{\left(\frac{1}{n} \eta_{kk}^{(2)} \right)^2}{\left( \frac{1}{n}\trace(\mtx{V}_\gamma^{-1})\right)^2 }.
\end{align}
Then by triangle inequality, for any $l,m = 1,2$,
\begin{align*}
     &~~~\frac{1}{n} \eta_{kk}^{(l)}\frac{1}{n} \eta_{kk}^{(m)} - \frac{1}{n} \trace \left(\mtx{V}_{\gamma}^{-l} \frac{1}{p}\mtx{Z} \mtx{Z}^\top \right)\frac{1}{n} \trace\left(\mtx{V}_{\gamma}^{-m} \frac{1}{p}\mtx{Z} \mtx{Z}^\top \right) \\
     & \leq \frac{1}{n} \trace \left(\mtx{V}_{\gamma}^{-l}\frac{1}{p}\mtx{Z} \mtx{Z}^\top  \right)\left| \frac{1}{n} \eta_{kk}^{(m)}  - \frac{1}{n} \trace \left(\mtx{V}_{\gamma}^{-m}\frac{1}{p}\mtx{Z} \mtx{Z}^\top 
 \right)\right| \\
 &~~~ + \frac{1}{n} \trace \left(\mtx{V}_{\gamma}^{-m}\frac{1}{p}\mtx{Z} \mtx{Z}^\top 
 \right) \left|\frac{1}{n} \eta_{kk}^{(l)} - \frac{1}{n} \trace \left(\mtx{V}_{\gamma}^{-l} \frac{1}{p}\mtx{Z} \mtx{Z}^\top \right)\right|\\
 &~~~ + \left| \frac{1}{n} \eta_{kk}^{(m)}  - \frac{1}{n} \trace \left(\mtx{V}_{\gamma}^{-m}\frac{1}{p}\mtx{Z} \mtx{Z}^\top 
 \right)\right|\left|\frac{1}{n} \eta_{kk}^{(l)} - \frac{1}{n} \trace \left(\mtx{V}_{\gamma}^{-l}\frac{1}{p}\mtx{Z} \mtx{Z}^\top  \right)\right|.
\end{align*}
By \eqref{eq:thm_D1} and the fact that
\[
\frac{1}{n} \trace \left(\mtx{V}_{\gamma}^{-l} \frac{1}{p}\mtx{Z} \mtx{Z}^\top \right) = O_P(1), \quad \frac{1}{n} \trace \left(\mtx{V}_{\gamma}^{-m}\frac{1}{p}\mtx{Z} \mtx{Z}^\top  \right) = O_P(1),
\]
we can have
\begin{align}
\label{diff:eta_tr}
\max_{1 \leq k \leq p}\left|\frac{1}{n} \eta_{kk}^{(l)}\frac{1}{n} \eta_{kk}^{(m)} - \frac{1}{n} \trace \left(\mtx{V}_{\gamma}^{-l} \frac{1}{p}\mtx{Z} \mtx{Z}^\top \right)\frac{1}{n} \trace \left(\mtx{V}_{\gamma}^{-m}\frac{1}{p}\mtx{Z} \mtx{Z}^\top  \right) \right|=  O_P\left( \sqrt{\frac{\log n}{n}}  \right).
\end{align}
Since\begin{align*}
    &\left( \trace\left( \mtx{B}_{\gamma} \frac{1}{p}\mtx{Z} \mtx{Z}^\top   \right) \right)^2\\
    &= \left(\frac{1}{n} \trace\left( \mtx{V}_\gamma^{-1} \frac{1}{p}\mtx{Z} \mtx{Z}^\top   \right) - \frac{ \frac{1}{n}\trace\left( \mtx{V}_\gamma^{-1} \frac{1}{p}\mtx{Z} \mtx{Z}^\top   \right) }{\frac{1}{n}\trace(\mtx{V}_\gamma^{-1})}  \right)^2\\
    & = \left(\frac{1}{n} \trace\left( \mtx{V}_\gamma^{-1} \frac{1}{p}\mtx{Z} \mtx{Z}^\top   \right) \right)^2 -2 \frac{ \frac{1}{n} \trace\left( \mtx{V}_\gamma^{-1} \frac{1}{p}\mtx{Z} \mtx{Z}^\top   \right)\frac{1}{n} \trace\left( \mtx{V}_\gamma^{-2} \frac{1}{p}\mtx{Z} \mtx{Z}^\top   \right)}{\frac{1}{n}\trace(\mtx{V}_\gamma^{-1}) } \notag\\
    &~~~ + \frac{\left(\frac{1}{n} \trace\left( \mtx{V}_\gamma^{-2} \frac{1}{p}\mtx{Z} \mtx{Z}^\top   \right) \right)^2}{\left(\frac{1}{n}\trace(\mtx{V}_\gamma^{-1})\right)^2 },
\end{align*}
then by \eqref{eq:z1_B_z1} and \eqref{diff:eta_tr}, there holds that
\begin{align}  \label{eq:diff_zkBzk_tr}
 \max_{1 \leq k \leq p}\left| \left(\vct{z}_k^\top \mtx{B}_{\gamma} \vct{z}_k \right)^2 -  \left( \trace\left( \mtx{B}_{\gamma} \frac{1}{p}\mtx{Z} \mtx{Z}^\top   \right) \right)^2 \right| = O_P\left( \sqrt{\frac{\log n}{n}}  \right).
\end{align}
Finally, as we can know from \eqref{eq:diag_homo_1}, \eqref{eq:diag_homo_3} and Lemma \ref{thm:diff_tr_h} that $ (np)^{-1} \trace\left(\mtx{V}_{\gamma}^{-1} \mtx{Z} \mtx{Z}^\top\right) $ converges to the same limit as 
\[
\frac{ n^{-1} \trace(\mtx{V}_{\gamma}^{-1})}{1 +  \gamma p^{-1} \trace(V_{\gamma}^{-1})}, 
\]
which means, with $\tau_n=n/p$,
\[
\frac{1}{\gamma}(1 - h_1(\gamma,\tau_n)) = \frac{h_1(\gamma,\tau_n)}{ 1 + \gamma \tau_n h_1(\gamma,\tau_n)}.
\]
And Lemma \ref{thm:diff_tr_h} shows that
\begin{align*}
    \left| \frac{1}{np} \trace\left(\mtx{V}_{\gamma}^{-1} \mtx{Z} \mtx{Z}^\top\right) - \frac{1}{\gamma}(1 - h_1(\gamma,\tau_n))  \right| = O_P\left( \frac{1}{n}  \right),
\end{align*}
and
\begin{align*}
    \left|  \frac{1}{n}\frac{\trace(\mtx{V}_{\gamma}^{-1})}{1 +  \frac{\gamma}{p} \trace(V_{\gamma}^{-1})} - \frac{h_1(\gamma,\tau_n)}{ 1 + \gamma \tau_n h_1(\gamma,\tau_n)}  \right|=O_P\left( \frac{1}{n}  \right).
\end{align*}
Combine the above two inequalities, we can get \eqref{eq:diag_homo_5}. Similarly, by \eqref{eq:diag_homo_2}, \eqref{eq:diag_homo_3} and Lemma \ref{thm:diff_tr_h} we can get \eqref{eq:diag_homo_6}.

\subsection{Proof of Lemma \ref{lem:diag_homo_taylor}}
By \eqref{eq:diag_homo_5} and \eqref{eq:mean_A1_rhok}, we can know that
\begin{align}\label{eq:eq:diag_homo_5_01}
&\max_{1\leq k \leq p} \Bigg\vert\frac{1}{n}\vct{z}_k^\top \mtx{V}_\gamma^{-1} \vct{z}_k - \frac{1}{np} \trace\left(\mtx{V}_{\gamma}^{-1} \mtx{Z} \mtx{Z}^\top\right) \notag
\\
&~~~~~-\frac{1}{\left(1 +  \frac{\gamma}{p} \trace(V_{\gamma}^{-1}) \right)^2}\left(\frac{1}{n}\eta_{kk,k}^{(1)} - \frac{1}{n}\trace(V_{\gamma}^{-1})  \right)\Bigg\vert  \notag\\
&=  \max_{1\leq k \leq p} \Bigg\vert\frac{1}{n} \frac{ \eta_{kk,k}^{(1)}}{1 + \frac{\gamma}{p}\eta_{kk,k}^{(1)} } -   \frac{1}{n}\frac{\trace(\mtx{V}_{\gamma}^{-1})}{1 +  \frac{\gamma}{p} \trace(V_{\gamma}^{-1})} \notag
\\
&~~~~~-\frac{1}{\left(1 +  \frac{\gamma}{p} \trace(V_{\gamma}^{-1}) \right)^2}\left(\frac{1}{n}\eta_{kk,k}^{(1)} - \frac{1}{n}\trace(V_{\gamma}^{-1})  \right)\Bigg\vert + O_P\left( \frac{1}{n}  \right),
\end{align}
and similarly by \eqref{eq:diag_homo_6} and \eqref{eq:mean_A2_phik}, there holds that
\begin{align*}
&\max_{1\leq k \leq p} \Bigg\vert\frac{1}{n}\vct{z}_k^\top \mtx{V}_\gamma^{-2}\vct{z}_k - \frac{1}{np} \trace\left(\mtx{V}_{\gamma}^{-2} \mtx{Z} \mtx{Z}^\top\right)
\\
&~~~~~~~~+\frac{\trace(\mtx{V}_{\gamma}^{-2})}{\left(1 +  \frac{\gamma}{p} \trace(V_{\gamma}^{-1}) \right)^3} \frac{2\gamma}{p}\left(\frac{1}{n}\eta_{kk,k}^{(1)} - \frac{1}{n}\trace(V_{\gamma}^{-1})  \right)
\\
&~~~~~~~~-\frac{1}{\left(1 +  \frac{\gamma}{p} \trace(V_{\gamma}^{-1}) \right)^2} \left(\frac{1}{n}\eta_{kk,k}^{(2)} - \frac{1}{n}\trace(V_{\gamma}^{-2})  \right)\Bigg\vert\\
&= \max_{1\leq k \leq p} \Bigg\vert    \frac{  \eta_{kk,k}^{(2)} }{ \left(1 + \frac{\gamma}{p}\eta_{kk,k}^{(1)} \right)^2 } - \frac{  \trace\left(\mtx{V}_{\gamma }^{-2} \right)  }{\left( 1 + \frac{\gamma}{p}\trace\left(\mtx{V}_{\gamma }^{-1}\right)\right)^2} \notag  \\
\\
&~~~~~~~~+\frac{\trace(\mtx{V}_{\gamma}^{-2})}{\left(1 +  \frac{\gamma}{p} \trace(V_{\gamma}^{-1}) \right)^3} \frac{2\gamma}{p}\left(\frac{1}{n}\eta_{kk,k}^{(1)} - \frac{1}{n}\trace(V_{\gamma}^{-1})  \right)
\\
&~~~~~~~~-\frac{1}{\left(1 +  \frac{\gamma}{p} \trace(V_{\gamma}^{-1}) \right)^2} \left(\frac{1}{n}\eta_{kk,k}^{(2)} - \frac{1}{n}\trace(V_{\gamma}^{-2})  \right)\Bigg\vert +  O_P\left( \frac{1}{n}  \right).
\end{align*}
Define
\begin{align*}
    z_\gamma(x) = \frac{x}{1+ \gamma\frac{n}{p}x}, \quad w_\gamma(x,y) = \frac{x}{\left( 1+ \gamma\frac{n}{p}y \right)^2}.
\end{align*}
By the Taylor series expansion, as $\frac{1}{n}\eta_{kk,k}^{(1)}  \rightarrow \frac{1}{n} \trace\left(\mtx{V}_{\gamma }^{-1} \right)$
\begin{align*}
     z_\gamma\left( \frac{1}{n}\eta_{kk,k}^{(1)} \right) &=  z_\gamma\left( \frac{1}{n} \trace\left(\mtx{V}_{\gamma }^{-1} \right)\right) + z_\gamma'\left( \frac{1}{n} \trace\left(\mtx{V}_{\gamma }^{-1} \right)\right) \left( \frac{1}{n}\eta_{kk,k}^{(1)} - \frac{1}{n} \trace\left(\mtx{V}_{\gamma }^{-1} \right)\right)\\
     &~~~ + R_1\left( \frac{1}{n}\eta_{kk,k}^{(1)} , \frac{1}{n} \trace\left(\mtx{V}_{\gamma }^{-1} \right)\right).
\end{align*}
Here $R_1$ is the remainder term \begin{align*}
    R_1 = \frac{1}{2}z_\gamma^{''}\left( c_k \right) \left( \frac{1}{n}\eta_{kk,k}^{(1)} - \frac{1}{n} \trace\left(\mtx{V}_{\gamma }^{-1} \right)\right)^2  
\end{align*}
where $c_k$ is some constant between $\frac{1}{n}\eta_{kk,k}^{(1)} $ and $ \frac{1}{n} \trace\left(\mtx{V}_{\gamma }^{-1} \right)$. Then 
\begin{align*}
    & \frac{  \frac{1}{n}\eta_{kk,k}^{(1)} }{ 1 + \frac{\gamma}{p}\eta_{kk,k}^{(1)}  } - \frac{ \frac{1}{n} \trace\left(\mtx{V}_{\gamma }^{-1} \right)  }{ 1 + \frac{\gamma}{p}\trace\left(\mtx{V}_{\gamma }^{-1}\right)}\\
    & = z_\gamma\left( \frac{1}{n}\eta_{kk,k}^{(1)} \right) - z_\gamma\left( \frac{1}{n} \trace\left(\mtx{V}_{\gamma }^{-1} \right)\right) \\
    & = z_\gamma'\left( \frac{1}{n} \trace\left(\mtx{V}_{\gamma }^{-1} \right)\right) \left( \frac{1}{n}\eta_{kk,k}^{(1)} - \frac{1}{n} \trace\left(\mtx{V}_{\gamma }^{-1} \right)\right) + R_1\left( \frac{1}{n}\eta_{kk,k}^{(1)} , \frac{1}{n} \trace\left(\mtx{V}_{\gamma }^{-1} \right)\right), 
\end{align*}
where
\begin{align*}
     z_\gamma'\left( x\right) = \frac{1}{\left(1+ \frac{\gamma n}{p} x  \right)^2}, \quad z_\gamma''\left( x\right) =\frac{\gamma n}{p}\frac{1 }{ \left(1+ \frac{\gamma n}{p} x \right)^3}.
\end{align*}
By \eqref{eq:eq:diag_homo_5_01}, this implies that
\begin{align*}
&\max_{1\leq k \leq p} \Bigg\vert\frac{1}{n}\vct{z}_k^\top \mtx{V}_\gamma^{-1} \vct{z}_k - \frac{1}{np} \trace\left(\mtx{V}_{\gamma}^{-1} \mtx{Z} \mtx{Z}^\top\right)
\\
&~~~~~-\frac{1}{\left(1 +  \frac{\gamma}{p} \trace(V_{\gamma}^{-1}) \right)^2}\left(\frac{1}{n}\eta_{kk,k}^{(1)} - \frac{1}{n}\trace(V_{\gamma}^{-1})  \right)\Bigg\vert\\
&=  \max_{1\leq k\leq p} \left|R_1\left( \frac{1}{n}\eta_{kk,k}^{(1)} , \frac{1}{n} \trace\left(\mtx{V}_{\gamma }^{-1} \right)\right)  \right| + O_P\left( \frac{1}{n}  \right),
\end{align*}
and by Lemma \ref{lem:leave_one_out},
\begin{align*}
     &\max_{1\leq k\leq p} \left|R_1\left( \frac{1}{n}\eta_{kk,k}^{(1)} , \frac{1}{n} \trace\left(\mtx{V}_{\gamma }^{-1} \right)\right)  \right|\\
     & = \max_{1\leq k\leq p} \left| \frac{1}{2}z_\gamma^{''}\left( c_k \right) \left( \frac{1}{n}\eta_{kk,k}^{(1)} - \frac{1}{n} \trace\left(\mtx{V}_{\gamma }^{-1} \right)\right)^2   \right|\\
     & \leq \frac{1}{2}\max_{1\leq k\leq p} \left|\frac{\gamma n}{p}\frac{1 }{ \left(1+ \frac{\gamma n}{p}c_k  \right)^3} \right|\max_{1\leq k\leq p} \left| \left(\frac{1}{n}\eta_{kk,k}^{(1)} - \frac{1}{n} \trace\left(\mtx{V}_{\gamma }^{-1} \right)\right)^2   \right|\\
     & \leq \frac{\gamma n}{2p}\max_{1\leq k\leq p} \left| \left(\frac{1}{n}\eta_{kk,k}^{(1)} - \frac{1}{n} \trace\left(\mtx{V}_{\gamma }^{-1} \right)\right)^2   \right|\\
     & = O_P\left(\frac{\log n}{n} \right).
\end{align*}
Combine the above two inequalities, we can get \eqref{eq:diag_homo_taylor_1}.\\
Similarly, by the Taylor series expansion, as $\frac{1}{n}\eta_{kk,k}^{(2)}  \rightarrow \frac{1}{n} \trace\left(\mtx{V}_{\gamma }^{-2} \right)$, we can have
\begin{align*}
    &\max_{1\leq k \leq p} \Bigg\vert\frac{1}{n}\vct{z}_k^\top \mtx{V}_\gamma^{-2}\vct{z}_k - \frac{1}{np} \trace\left(\mtx{V}_{\gamma}^{-2} \mtx{Z} \mtx{Z}^\top\right)
\\
&~~~~~~~~+\frac{\trace(\mtx{V}_{\gamma}^{-2})}{\left(1 +  \frac{\gamma}{p} \trace(V_{\gamma}^{-1}) \right)^3} \frac{2\gamma}{p}\left(\frac{1}{n}\eta_{kk,k}^{(1)} - \frac{1}{n}\trace(V_{\gamma}^{-1})  \right)
\\
&~~~~~~~~-\frac{1}{\left(1 +  \frac{\gamma}{p} \trace(V_{\gamma}^{-1}) \right)^2} \left(\frac{1}{n}\eta_{kk,k}^{(2)} - \frac{1}{n}\trace(V_{\gamma}^{-2})  \right)\Bigg\vert\\
&=\max_{1\leq k\leq p} \left|\tilde{R}_1\left( \frac{1}{n}\eta_{kk,k}^{(2)}, \frac{1}{n}\trace\left(\mtx{V}_{\gamma }^{-2}\right) ,\frac{1}{n}\eta_{kk,k}^{(1)},\frac{1}{n} \trace\left(\mtx{V}_{\gamma }^{-1} \right)  \right)   \right| + O_P\left( \frac{1}{n}  \right).
\end{align*}
Since
\begin{align*}
      \frac{\partial^2 w_\gamma}{\partial x^2}\left(x,y \right) & = 0,\\
    \left| \frac{\partial^2 w_\gamma}{\partial x\partial y}\left(x,y \right)  \right|& =  \left| -2 \frac{\gamma n}{p} \frac{1}{\left( 1+ \gamma\frac{n}{p}y \right)^3 }\right| \leq 2 \frac{\gamma n}{p} \quad \text{for } y \geq 0,\\
    \left| \frac{\partial^2 w_\gamma}{\partial y^2}\left(x,y \right)\right| & = 6 \left(\frac{\gamma n}{p}\right)^2 \frac{x}{\left( 1+ \gamma\frac{n}{p}y \right)^3 } \leq 6 \left(\frac{\gamma n}{p}\right)^2  x \quad \text{for } x,y \geq 0,
\end{align*}
by Lemma \ref{lem:leave_one_out} we can have
\begin{align*}
    &\max_{1\leq k \leq q}\left|\tilde{R}_1\left( \frac{1}{n}\eta_{kk,k}^{(2)}, \frac{1}{n}\trace\left(\mtx{V}_{\gamma }^{-2}\right) ,\frac{1}{n}\eta_{kk,k}^{(1)},\frac{1}{n} \trace\left(\mtx{V}_{\gamma }^{-1} \right)  \right)\right|\\
    & \leq \max_{1\leq k \leq q}\left| \frac{\partial^2 w_\gamma}{\partial x\partial y}\left(c_{k1},c_{k2} \right)\left( \frac{1}{n}\eta_{kk,k}^{(2)} - \frac{1}{n}\trace\left(\mtx{V}_{\gamma }^{-2}\right) \right)\left(\frac{1}{n}\eta_{kk,k}^{(1)} - \frac{1}{n} \trace\left(\mtx{V}_{\gamma }^{-1} \right)   \right)\right| \\
    &~~~ + \max_{1\leq k \leq q}\left|\frac{ \frac{\partial^2 w_\gamma}{\partial y^2}\left(c_{k1},c_{k2} \right)}{2}\left(\frac{1}{n}\eta_{kk,k}^{(1)} - \frac{1}{n} \trace\left(\mtx{V}_{\gamma }^{-1} \right)   \right) ^2\right|\\
    &\leq 2 \frac{\gamma n}{p} \max_{1\leq k \leq q}\left|\left( \frac{1}{n}\eta_{kk,k}^{(2)} - \frac{1}{n}\trace\left(\mtx{V}_{\gamma }^{-2}\right) \right)\right|\max_{1\leq k \leq q}\left|\left(\frac{1}{n}\eta_{kk,k}^{(1)} - \frac{1}{n} \trace\left(\mtx{V}_{\gamma }^{-1} \right)   \right)\right|\\
    &~~~ + 3 \left(\frac{\gamma n}{p}\right)^2  c_{k1}\max_{1\leq k \leq q}\left|\left(\frac{1}{n}\eta_{kk,k}^{(1)} - \frac{1}{n} \trace\left(\mtx{V}_{\gamma }^{-1} \right)   \right)^2\right|\\
    & = O_P\left( \frac{\log n}{n}\right),
\end{align*}
where $c_{k1}$ is some constant between $\frac{1}{n} \eta_{kk,k}^{(2)} $ and $ \frac{1}{n} \trace\left(\mtx{V}_{\gamma }^{-2} \right)$ and $c_{k2}$ is some constant between $\frac{1}{n}\eta_{kk,k}^{(1)}$ and $ \frac{1}{n} \trace\left(\mtx{V}_{\gamma }^{-1} \right) $.

\subsection{Proof of Lemma \ref{lem:diag_homo_moment}}

Write
\[
X_{k,l}
=
\frac{1}{n}\left\{
\eta_{kk,k}^{(l)}
-
\trace(\mtx V_{\gamma,-k}^{-l})
\right\},
\qquad l=1,2.
\]
We first prove the variance bound. Conditional on $\mtx Z_{-k}$,
the matrix $\mtx V_{\gamma,-k}^{-l}$ is deterministic and independent of
$\vct z_k$. Since the entries of $\vct z_k$ are independent, centered,
unit-variance, and uniformly sub-Gaussian, the standard quadratic-form
moment bound gives
\[
\E\left[
\left\{
\vct z_k^\top \mtx V_{\gamma,-k}^{-l}\vct z_k
-
\trace(\mtx V_{\gamma,-k}^{-l})
\right\}^2
\Bigm| \mtx Z_{-k}
\right]
\leq
C\trace(\mtx V_{\gamma,-k}^{-2l}).
\]
Because $\|\mtx V_{\gamma,-k}^{-1}\|\leq 1$, we have
$\trace(\mtx V_{\gamma,-k}^{-2l})\leq n$. Hence
\[
\E X_{k,l}^2\leq \frac{C}{n},
\]
uniformly over $k$.

We next prove the cross-moment bound. Fix $i\neq j$ and put
\[
\mtx A=\mtx V_{\gamma,-ij}.
\]
We give the argument for $l=1$; the case $l=2$ follows from the same
resolvent expansion with one additional bounded resolvent factor.
By the Sherman--Morrison formula,
\[
\mtx V_{\gamma,-i}^{-1}
=
\mtx A^{-1}
-
\frac{\gamma}{p}
\frac{\mtx A^{-1}\vct z_j\vct z_j^\top\mtx A^{-1}}
{1+\frac{\gamma}{p}\vct z_j^\top \mtx A^{-1}\vct z_j}.
\]
Therefore
\[
X_{i,1}
=
U_i+R_i,
\]
where
\[
U_i
=
\frac{1}{n}
\left\{
\vct z_i^\top \mtx A^{-1}\vct z_i
-
\trace(\mtx A^{-1})
\right\},
\]
and
\[
R_i
=
-\frac{1}{n}\frac{\gamma}{p}
\frac{
(\vct z_i^\top \mtx A^{-1}\vct z_j)^2
-
\vct z_j^\top \mtx A^{-2}\vct z_j
}{
1+\frac{\gamma}{p}\vct z_j^\top \mtx A^{-1}\vct z_j
}.
\]
Similarly,
\[
X_{j,1}=U_j+R_j,
\qquad
U_j
=
\frac{1}{n}
\left\{
\vct z_j^\top \mtx A^{-1}\vct z_j
-
\trace(\mtx A^{-1})
\right\}.
\]
Conditional on $\mtx A$, the random vectors $\vct z_i$ and $\vct z_j$ are
independent. Thus
\[
\E[U_iU_j\mid \mtx A]=0.
\]
Moreover, by the same quadratic-form moment bound and the fact
$\|\mtx A^{-1}\|\leq 1$,
\[
\E[U_i^2\mid \mtx A]\leq \frac{C}{n},
\qquad
\E[U_j^2\mid \mtx A]\leq \frac{C}{n}.
\]
For the remainder term, using
\[
\frac{\frac{\gamma}{p}\vct z_j^\top \mtx A^{-1}\vct z_j}
{1+\frac{\gamma}{p}\vct z_j^\top \mtx A^{-1}\vct z_j}
\leq 1
\]
and the Hanson--Wright moment bound, we obtain
\[
\E[R_i^2\mid \mtx A]\leq \frac{C}{n^2},
\qquad
\E[R_j^2\mid \mtx A]\leq \frac{C}{n^2}.
\]
Consequently, by Cauchy's inequality,
\[
|\E X_{i,1}X_{j,1}|
\leq
|\E U_iU_j|
+
|\E U_iR_j|
+
|\E U_jR_i|
+
|\E R_iR_j|
\leq
\frac{C}{n^{3/2}}.
\]
The same argument applies to $l=2$. Indeed, applying the resolvent identity
to $\mtx V_{\gamma,-i}^{-2}$ gives a finite linear combination of terms of
the form
\[
\frac{1}{n}
\left\{
\vct z_i^\top \mtx A^{-a}
\vct z_i
-
\trace(\mtx A^{-a})
\right\},
\qquad a=1,2,
\]
plus rank-one remainder terms whose conditional second moments are
$O(n^{-2})$. Since all resolvent powers have operator norm at most one,
the same estimates yield
\[
\max_{i\neq j}|\E X_{i,2}X_{j,2}|
\leq Cn^{-3/2}.
\]
This proves the lemma.

\subsection{Proof of Lemma \ref{lem:leave_two_out} }
Note that
\begin{align*}
    (\vct{z}_k^\top \mtx{V}_{\gamma}^{-1} \vct{z}_j)^2  =  (1+ \frac{\gamma}{p}\rho_k)^{-2}  (\vct{z}_k^\top \mtx{V}_{\gamma,-k}^{-1} \vct{z}_j)^2 \leq  \left(\vct{z}_k^\top \mtx{V}_{\gamma,-k}^{-1} \vct{z}_j\right)^2
\end{align*}
and
\begin{align}
    &(\vct{z}_k^\top \mtx{V}_{\gamma}^{-2} \vct{z}_j)^2 
    \\
    &= \left((1+ \frac{\gamma}{p}\rho_k)^{-1}  \vct{z}_k^\top \mtx{V}_{\gamma,-k}^{-2}\vct{z}_j + \left(-\frac{\gamma}{p} \phi_k (1+ \frac{\gamma}{p}\rho_k)^{-2} \vct{z}_k^\top \mtx{V}_{\gamma, -k}^{-1} \vct{z}_j  \right) \right)^2 \nonumber\\
    & \leq \frac{1}{2} (1+ \frac{\gamma}{p}\rho_k)^{-2}  \left( \vct{z}_k^\top \mtx{V}_{\gamma,-k}^{-2}\vct{z}_j \right)^2  + \frac{1}{2}\left(\frac{\gamma}{p} \phi_k\right)^2 (1+ \frac{\gamma}{p}\rho_k)^{-4}  \left( \vct{z}_k^\top \mtx{V}_{\gamma,-k}^{-1}\vct{z}_j \right)^2 \nonumber\\
     &\leq 2 \left(  \left( \vct{z}_k^\top \mtx{V}_{\gamma,-k}^{-2}\vct{z}_j \right)^2 +   \left( \vct{z}_k^\top \mtx{V}_{\gamma,-k}^{-1}\vct{z}_j \right)^2  \right) \label{ieq: zk_invV2_zj}
\end{align}
where the last inequality is due to $0 < \phi_k \leq \rho_k$.

Denote $\mtx{Z}_{-k} = [\vct{z}_1, \ldots, \vct{z}_{k-1}, \vct{z}_{k+1}, \ldots, \vct{z}_p]$. Note that the components of $\vct{z}_k$ are independent mean-zero sub-Gaussian random variables, conditional on $\mtx{Z}_{-k}$, by Proposition \ref{ineq:Hoeffding}, we have, for any $k \neq j$ and $t \geq  0$,
\[
\P\left\{|\vct{z}_k^\top \mtx{V}_{\gamma, -k}^{-1} \vct{z}_j| \geq t  \Big\vert \mtx{Z}_{-k}\right\} \leq  e \exp \left\{ -c \frac{t^2}{K^2 \|\mtx{V}_{\gamma, -k}^{-1} \vct{z}_j \|_2^2}\right\},
\]
where $c$ and $K$ are some positive constants. By letting $t = K \sqrt{\frac{3 \log p}{c}}  \|\mtx{V}_{\gamma, -k}^{-1} \vct{z}_j \|_2$, it follows 
\[
\P\left\{|\vct{z}_k^\top \mtx{V}_{\gamma, -k}^{-1} \vct{z}_j| \geq C\sqrt{\log p} \|\mtx{V}_{\gamma, -k}^{-1} \vct{z}_j\| \Big\vert \mtx{Z}_{-k}\right\} \leq \frac{e}{p^3},
\]
where $C$ is some positive constant. It further implies the unconditional probability inequality
\begin{align}\label{ieq: zk_V_zj}
\P\left\{|\vct{z}_k^\top \mtx{V}_{\gamma, -k}^{-1} \vct{z}_j| \geq C\sqrt{\log p} \|\mtx{V}_{\gamma, -k}^{-1} \vct{z}_j\| \right\} \leq \frac{e}{p^3}.
\end{align}
By the fact $\|\mtx{V}_{\gamma, -k}^{-1}\| \leq 1$, we have $\|\mtx{V}_{\gamma, -k}^{-1} \vct{z}_j\| \leq \|\vct{z}_j\|$. Note that $z_{1j}^2-1, z_{2j}^2-1, \ldots, z_{nj}^2-1$ are independent centered sub-exponential random variables, by the Proposition \ref{ineq:Bernstein}, we have
\[
\P\left\{\left\vert \sum_{i=1}^n (z_{ij}^2-1) \right\vert \geq t  \right\} \leq  2 \exp \left\{ -c \min \left( \frac{t^2}{K^2 n }, \frac{t}{K} \right) \right\},
\]
where $c$ and $K$ are some positive constants. Take $t = K\sqrt{3n\log p}$, then we get
\begin{align}
\P\left\{\left\vert \|\vct{z}_j\|^2 - n \right\vert \geq C\sqrt{n \log p} \right\} \leq \frac{2}{p^3},
\label{ieq: norm_zk}
\end{align}
for some constant $C$. Combining the above \eqref{ieq: zk_V_zj} and \eqref{ieq: norm_zk} together, with probability at least $1 - (2+e)/p^3$, there holds
\[
(\vct{z}_k^\top \mtx{V}_{\gamma}^{-1} \vct{z}_j)^2 \leq  (\vct{z}_k^\top \mtx{V}_{\gamma,-k}^{-1} \vct{z}_j)^2 \leq C (n + C\sqrt{n \log p}) \log p.
\]
By a similar argument with the fact that  $\|\mtx{V}_{\gamma, -k}^{-2}\| \leq 1$, we can have with probability at least $1 - (2+e)/p^3$
\[
 (\vct{z}_k^\top \mtx{V}_{\gamma}^{-2} \vct{z}_j)^2 \leq   2 \left(  \left( \vct{z}_k^\top \mtx{V}_{\gamma,-k}^{-2}\vct{z}_j \right)^2 +   \left( \vct{z}_k^\top \mtx{V}_{\gamma,-k}^{-1}\vct{z}_j \right)^2  \right) \leq 4 C (n + C\sqrt{n \log p}) \log p,
\]
for some constant $C$.

Above inequalities imply that with probability at least $1 - (2+e)/p$, there hold
\begin{align*}
\max_{k \neq j} \vert \vct{z}_k^\top \mtx{V}_{\gamma}^{-1} \vct{z}_j|^2 &\leq C (n + C\sqrt{n \log p}) \log p, \quad \text{and}
\\
\max_{k \neq j} \vert \vct{z}_k^\top \mtx{V}_{\gamma}^{-2} \vct{z}_j|^2 &\leq 4 C (n + C\sqrt{n \log p}) \log p. 
\end{align*}
Then, there follows that
\begin{align*}
\max_{k \neq j} \vert \vct{z}_k^\top \mtx{V}_{\gamma}^{-1} \vct{z}_j|^2 = O_P (n \log p)
\quad \text{and} \quad 
\max_{k \neq j} \vert \vct{z}_k^\top \mtx{V}_{\gamma}^{-2} \vct{z}_j|^2 = O_P (n \log p).
\end{align*}

Next, define
\begin{align*}
    \bar{\theta}_1(\gamma,\tau) &=  \kappa_{1,1}(\gamma,\tau) -2 \frac{\kappa_{1,2}(\gamma,\tau) }{h_{1}(\gamma,\tau)} + \frac{ \kappa_{2,2}(\gamma,\tau)}{ h^2_{1}(\gamma,\tau) },
\end{align*}
where
\[
\kappa_{m,l}(\gamma,\tau) = \sum_{q_1 = 1}^l\sum_{q_2 = 1}^m \bar{a}_{q_1}^{(l)}(\gamma,\tau) \bar{a}_{q_2}^{(m)}(\gamma,\tau) h_{q_1+q_2}(\gamma,\tau),
\]
and
\[
\bar{a}_{1}^{(1)}(\gamma, \tau)=\frac{1}{\left(1+\tau \gamma h_{1}(\gamma, \tau) \right)^{2}}, \quad \bar{a}_{1}^{(2)}(\gamma, \tau)=\frac{-2 \tau \gamma h_{2}(\gamma, \tau)}{\left(1+\tau \gamma h_{1}(\gamma, \tau) \right)^{3}}, \quad \bar{a}_{2}^{(2)}(\gamma, \tau)=\frac{1}{\left(1+\tau \gamma h_{1}(\gamma, \tau) \right)^{2}}.
\]
Recall that $h_1(\gamma, \tau)$ and $h_2(\gamma, \tau)$ are defined in \eqref{eq:h_MP}.
Now, by the definition of $\eta_{ij}^{(l)}$ in \eqref{def:eta_noC}, we can rewrite $\left(\vct{z}_i^\top \mtx{B}_{\gamma} \vct{z}_j \right)^2$ as
\begin{align}\label{eq:zi_B_zj}
    \left(\vct{z}_i^\top \mtx{B}_{\gamma} \vct{z}_j \right)^2 &=  \left(\frac{1}{n}\vct{z}_i^\top\mtx{V}_\gamma^{-1}\vct{z}_j - \frac{\frac{1}{n} \vct{z}_i^\top \mtx{V}_\gamma^{-2}\vct{z}_j }{\frac{1}{n}\trace(\mtx{V}_\gamma^{-1})}  \right)^2\notag \\
    & = \left(\frac{1}{n} \eta_{ij}^{(1)} \right)^2 -2 \frac{ \frac{1}{n} \eta_{ij}^{(1)}\frac{1}{n} \eta_{ij}^{(2)}}{\frac{1}{n}\trace(\mtx{V}_\gamma^{-1}) }  + \frac{\left(\frac{1}{n} \eta_{ij}^{(2)} \right)^2}{\left(  \frac{1}{n}\trace(\mtx{V}_\gamma^{-1})\right)^2 }.
\end{align}
The following results are implied by \cite{jiang2016high} in the supplementary material.

\begin{proposition}[\cite{jiang2016high}]\label{thm:leave_two_out}
  For any $i \neq j$ and $i,j \geq 1$, we have 
    $$
\begin{aligned}
\eta_{ij}^{(1)} &= \bar{a}_{1 ; i j}^{(1)} \eta_{i j ; i j}^{\left(1\right)},\\
\eta_{ij}^{(2)} & = \bar{a}_{1 ; i j}^{(2)} \eta_{i j ; i j}^{\left(1\right)} +  \bar{a}_{2 ; i j}^{(2)} \eta_{i j ; i j}^{\left(2\right)},
\end{aligned}
$$
with
$$
\begin{aligned}
\bar{a}_{1 ; i j}^{(1)}=& \frac{1}{\left(1+\frac{\gamma}{p}\eta_{i i ; i}^{(1)}\right)\left(1+\frac{\gamma}{p}\eta_{j j ; i j}^{(1)}\right)}, \\
\bar{a}_{1 ; i j}^{(2)}=& \frac{-\frac{\gamma}{p}\eta_{i i ; i}^{(2)}}{\left(1+\frac{\gamma}{p}\eta_{i i ; i}^{(1)}\right)^{2}\left(1+\frac{\gamma}{p}\eta_{j j ; i j}^{(1)}\right)}+\frac{-\frac{\gamma}{p}\eta_{j j ; i j}^{(2)}}{\left(1+\frac{\gamma}{p}\eta_{i i ; i}^{(1)}\right)\left(1+\frac{\gamma}{p}\eta_{j j ; i j}^{(1)}\right)^{2}}, \\
\bar{a}_{2 ; i j}^{(2)}=& \bar{a}_{1 ; i j}^{(1)}.
\end{aligned}
$$
And
\[
\max _{1 \leq i \neq j \leq p} \max _{1 \leq l \leq 2} \max _{1 \leq q_{1} \leq l}\left|\bar{a}_{q_{1} ; i j}^{(l)}-\bar{a}_{q_{1}}^{(l)}(\gamma, \tau)\right|=O_{\mathrm{P}}\left(\sqrt{\frac{\log p}{n}}\right) .
\]
Furthermore,
\begin{align}\label{eq:diff_aij_a}
   &~~~\max _{1 \leq i \neq j \leq p} \max _{1 \leq l,m \leq 2} \max _{\substack{1 \leq q_{1} \leq l\\1 \leq q_2 \leq m}}\left|\bar{a}_{q_{1} ; i j}^{(l)}\bar{a}_{q_{2} ; i j}^{(m)}-\bar{a}_{q_{1}}^{(l)}(\gamma, \tau)\bar{a}_{q_{2}}^{(m)}(\gamma, \tau)\right|=O_{\mathrm{P}}\left(\sqrt{\frac{\log p}{n}}\right) 
\end{align}
\end{proposition}

\begin{proposition}[\cite{jiang2016high}]
\label{thm:diff_tr_V_Vij}
   For any $l \geq 1$ and $ 1 \leq i \neq j$,
   \begin{align}\label{eq:diff_tr_V_Vij}
    \frac{1}{n} \left|\trace\left( \mtx{V}_{\gamma}^{-l}\right)- \trace\left( \mtx{V}_{\gamma,-ij}^{-l}\right)  \right| \leq \frac{1}{n} 2^{l+1}.
\end{align}
\end{proposition}

\begin{proposition}[\cite{jiang2016high}]\label{thm:mean_dij}
For any $1 \leq q_1, q_2 \leq 2$, define 
\[
d_{ij}^{(q_1,q_2)} \coloneqq  \frac{1}{n}\eta_{ij;ij}^{(q_1)}\eta_{ij;ij}^{(q_2)}  - \frac{1}{n} \trace \left(\mtx{V}_{\gamma,-ij}^{-(q_1 + q_2)}\right).
\]
Then the following statements are true.
\begin{itemize}
    \item [1)] For some constant $K_1 > 0$,
    \begin{align*}
        \max_{1\leq i \neq j \leq p} \E \left[ \left(d_{ij}^{(q_1,q_2)}\right)^2 \right] \leq K_1. 
    \end{align*}
    \item [2)] For any $ i \neq j \neq i'$, $j' \neq i'$ (either $j=j'$ or not) and some constant $K_2 > 0$,
    \begin{align*}
         \max_{ i \neq j \neq i',\, j' \neq i'} \left\vert \E \left[ d_{ij}^{(q_1,q_2)}d_{i'j'}^{(q_1,q_2)}\right] \right\vert \leq \frac{K_2}{\sqrt{n}}. 
    \end{align*}
\end{itemize}
    
\end{proposition}

By \eqref{eq:zi_B_zj}, let's first show that for $1 \leq l,m \leq 2$,
\begin{align}\label{eq:diff_beta_eta_kappa}
\frac{1}{p(p-1)} \sum_{i\neq j}\frac{1}{n} \eta_{ij}^{(l)}\eta_{ij}^{(m)}  \stackrel{P}{\longrightarrow} \kappa_{m,l}(\gamma,
\tau), \text{ and } \sum_{i\neq j} \beta_i^2 \beta_j^2 \frac{1}{n} \eta_{ij}^{(l)}\eta_{ij}^{(m)}  \stackrel{P}{\longrightarrow} \| \vct{\beta}\|^4 \kappa_{m,l}(\gamma,\tau).
\end{align}
Since $\eta_{ij;ij}^{(q_1)}\eta_{ij;ij}^{(q_2)} > 0$ for any $q_1,q_2 = 1,2,\ldots$, by Proposition \ref{thm:leave_two_out} we can have
\begin{align}\label{eq:sum_beta_eta_ij_ij_leave_two}
   \sum_{i\neq j} \beta_i^2 \beta_j^2 \frac{1}{n} \eta_{ij}^{(l)}\eta_{ij}^{(m)}  &=\sum_{q_1 = 1}^l\sum_{q_2 = 1}^m\left(\sum_{i\neq j}\bar{a}_{q_{1} ; i j}^{(l)}\bar{a}_{q_{2} ; i j}^{(m)} \beta_i^2 \beta_j^2 \frac{1}{n} \eta_{ij;ij}^{(q_1)}\eta_{ij;ij}^{(q_2)}  \right) \notag\\
   &=\sum_{q_1 = 1}^l\sum_{q_2 = 1}^m\bar{a}_{q_{1}}^{(l)}(\gamma, \tau)\bar{a}_{q_{2}}^{(m)}(\gamma, \tau)\left(\sum_{i\neq j} \beta_i^2 \beta_j^2 \frac{1}{n} \eta_{ij;ij}^{(q_1)}\eta_{ij;ij}^{(q_2)}  \right) + o_P(1),
\end{align}
and
\begin{align}\label{eq:sum_eta_ij_ij_leave_two}
    \frac{1}{p(p-1)} \sum_{i\neq j} \frac{1}{n} \eta_{ij}^{(l)}\eta_{ij}^{(m)}
    =\sum_{q_1 = 1}^l\sum_{q_2 = 1}^m\bar{a}_{q_{1}}^{(l)}(\gamma, \tau)\bar{a}_{q_{2}}^{(m)}(\gamma, \tau)\frac{1}{p(p-1)} \sum_{i\neq j}\frac{1}{n} \eta_{ij;ij}^{(q_1)}\eta_{ij;ij}^{(q_2)} + o_P(1).
\end{align}

We know that 
\begin{align*}
    \E\left[  \sum_{i\neq j}  \beta_i^2 \beta_j^2 d_{ij}^{(q_1,q_2)}  \right] = 0,
\end{align*}
and by Proposition \ref{thm:mean_dij}, we have that
\begin{align*}
    \E \left( \sum_{i\neq j}  \beta_i^2 \beta_j^2 d_{ij}^{(q_1,q_2)} \right)^2  & = \sum_{ i \neq j} \beta_i^4 \beta_j^4 \E \left[ \left(d_{ij}^{(q_1,q_2)}\right)^2 \right] + 2 \sum_{i\neq i' \neq j} \beta_i^2 \beta_{i'}^2 \beta_j^4  \E \left[ \left(d_{ij}^{(q_1,q_2)}d_{i'j}^{(q_1,q_2)}\right) \right]\\
    &~~~ +\sum_{i \neq j\neq i'\neq j'} \beta_i^2 \beta_j^2  \beta_{i'}^2 \beta_{j'}^2 \E \left[ \left(d_{ij}^{(q_1,q_2)}d_{i'j'}^{(q_1,q_2)}\right) \right]\\
    &\leq  K_1 \|\vct{\beta}\|_4^8 + 2 \frac{K_2}{\sqrt{n}} \| \vct{\beta}\|_2^4\|\vct{\beta}\|^4_4 + \frac{K_2}{\sqrt{n}}  \| \vct{\beta}\|_2^8\\
    & = o_P(1).
\end{align*}
This implies that
\begin{align*}
    \sum_{i\neq j} \beta_i^2 \beta_j^2 \frac{1}{n} \eta_{ij;ij}^{(q_1)}\eta_{ij;ij}^{(q_2)} =  \sum_{i\neq j}  \beta_i^2 \beta_j^2 \frac{1}{n}  \trace\left(  \mtx{V}_{\gamma,-ij}^{-(q_1+q_2)}\right) + o_P(1),
\end{align*}
where by Proposition \ref{thm:diff_tr_V_Vij} and Lemma \ref{thm:diff_tr_h},
\begin{align*}
    \left| \frac{1}{n}  \trace\left(  \mtx{V}_{\gamma,-ij}^{-(q_1+q_2)}\right) - h_{q_1 + q_2} (\gamma,\tau) \right| = o_P(1),
\end{align*}
and $ \sum_{i\neq j}  \beta_i^2 \beta_j^2 h_{q_1 + q_2} (\gamma,\tau) = \| \vct{\beta} \|^4  h_{q_1 + q_2} (\gamma,\tau) + o_P(1)$ since $\sum_{i=1}^p \beta_i^4 = o_P(1)$. Thus
\begin{align}\label{eq:sum_eta_ij_ij_leave_two_h_beta}
     \sum_{i\neq j} \beta_i^2 \beta_j^2 \frac{1}{n} \eta_{ij;ij}^{(q_1)}\eta_{ij;ij}^{(q_2)} =  \| \vct{\beta} \|^4  h_{q_1 + q_2} (\gamma,\tau) + o_P(1).
\end{align}
Thus by \eqref{eq:sum_beta_eta_ij_ij_leave_two} and \eqref{eq:sum_eta_ij_ij_leave_two_h_beta}, we can have 
\begin{align*}
     \sum_{i\neq j} \beta_i^2 \beta_j^2 \frac{1}{n} \eta_{ij}^{(l)}\eta_{ij}^{(m)} = \| \vct{\beta} \|^4 \kappa_{m,l}(\gamma,\tau)+ o_P(1).
\end{align*}
Similarly,
\begin{align*}
    \E\left[   \frac{1}{p(p-1)} \sum_{i\neq j} d_{ij}^{(q_1,q_2)}  \right] = 0,
\end{align*}
and
\begin{align*}
    &\E\left(  \frac{1}{p(p-1)} \sum_{i\neq j} d_{ij}^{(q_1,q_2)}  \right)^2 \\
    & = \frac{1}{p^2 (p-1)^2} \left(\sum_{i\neq j}  \E \left[ \left(d_{ij}^{(q_1,q_2)}\right)^2 \right] + 2 \sum_{i\neq i' \neq j} \E \left[ \left(d_{ij}^{(q_1,q_2)}d_{i'j}^{(q_1,q_2)}\right) \right] \right.\\
    &~~~ \left. +\sum_{i \neq j\neq i'\neq j'}  \E \left[ \left(d_{ij}^{(q_1,q_2)}d_{i'j'}^{(q_1,q_2)}\right) \right] \right)\\
    & \leq \frac{1}{p^2 (p-1)^2} \left( p(p-1) K_1 + 2 p(p-1)(p-2) \frac{K_2}{\sqrt{n}} +  p(p-1)(p-2)(p-3)\frac{K_2}{\sqrt{n}}  \right)\\
    & = o_P(1).
\end{align*}
Then by Proposition \ref{thm:diff_tr_V_Vij} and Lemma \ref{thm:diff_tr_h}, we can have that
\begin{align}\label{eq:sum_eta_ij_ij_leave_two_h}
 \frac{1}{p(p-1)} \sum_{i\neq j}\frac{1}{n} \eta_{ij;ij}^{(q_1)}\eta_{ij;ij}^{(q_2)}  = h_{q_1+q_2}(\gamma,\tau)+o_P(1),
\end{align}
which implies 
\begin{align*}
 \frac{1}{p(p-1)} \sum_{i\neq j}\frac{1}{n} \eta_{ij}^{(l)}\eta_{ij}^{(m)}  =\kappa_{m,l}(\gamma,\tau)+ o_P(1), 
\end{align*}
by \eqref{eq:sum_eta_ij_ij_leave_two}. 

Now we have proved \eqref{eq:diff_beta_eta_kappa}, then by \eqref{eq:zi_B_zj} and Lemma \ref{thm:diff_tr_h}, there holds that 
\begin{align*}
    \begin{dcases}
\frac{n}{p(p-1)}\sum_{i\neq j} \left(\vct{z}_i^\top \mtx{B}_{\gamma} \vct{z}_j \right)^2  &= \bar{\theta}_1(\gamma, \tau) + o_P(1),
\\
n\sum_{i\neq j} \beta_i^2 \beta_j^2 \left(\vct{z}_i^\top \mtx{B}_{\gamma} \vct{z}_j \right)^2  &= \|\vct{\beta}\|^4 \bar{\theta}_1(\gamma, \tau) + o_P(1).
\end{dcases}
\end{align*}

\subsection{Proof of Lemma \ref{lem:V_inv_diag}}
Let $\tilde{\vct{z}_i}^\top$ be the $i$th row of $\mtx{Z}$. By Sherman-Morrison-Woodbury formula, we have
\begin{align}\label{eq:V_gamma_inv1_ex}
   \mtx{V}_{\gamma}^{-1} = \left(\mtx{I}_n + \frac{\gamma}{p} \mtx{Z} \mtx{Z}^\top
 \right)^{-1} = \mtx{I}_n - \frac{\gamma}{p} \mtx{Z} \left(\mtx{I}_{p} + \frac{\gamma}{p}  \mtx{Z}^\top  \mtx{Z}\right)^{-1}\mtx{Z}^\top,
\end{align}
and
\begin{align}\label{eq:V_gamma_inv2_ex}
   \mtx{V}_{\gamma}^{-2} &=\left(\mtx{V}_{\gamma}^{-1}\right)^2\notag \\
   &= \mtx{I}_n - 2\frac{\gamma}{p} \mtx{Z} \left(\mtx{I}_p + \frac{\gamma}{p}  \mtx{Z}^\top  \mtx{Z}\right)^{-1}\mtx{Z}^\top + \left( \frac{\gamma}{p}\right)^2 \left(\mtx{Z}\left(\mtx{I}_{p} + \frac{\gamma}{p}  \mtx{Z}^\top  \mtx{Z}\right)^{-1}\mtx{Z}^\top\right)^2.
\end{align}
Combining \eqref{eq:V_gamma_inv1_ex} and \eqref{eq:V_gamma_inv2_ex} gives
\begin{align*}
    \left(\mtx{V}_{\gamma}^{-1} \right)_{ii} = 1 -  \frac{\gamma}{p} \tilde{\vct{z}}_i^\top
 \left(\mtx{I}_{p} + \frac{\gamma}{p}  \mtx{Z}^\top  \mtx{Z}\right)^{-1}\tilde{\vct{z}}_i,
\end{align*}
and
\begin{align*}
     \left(\mtx{V}_{\gamma}^{-2} \right)_{ii} = &1 - 2\frac{\gamma}{p}\tilde{\vct{z}}_i^\top
 \left(\mtx{I}_{p} + \frac{\gamma}{p}  \mtx{Z}^\top  \mtx{Z}\right)^{-1}\tilde{\vct{z}}_i
\\
&+ \left( \frac{\gamma}{p}\right)^2 \tilde{\vct{z}}_i^\top \left(\mtx{I}_{p} + \frac{\gamma}{p}  \mtx{Z}^\top  \mtx{Z}\right)^{-1}\mtx{Z}^\top \mtx{Z} \left(\mtx{I}_{p} + \frac{\gamma}{p}  \mtx{Z}^\top  \mtx{Z}\right)^{-1}\tilde{\vct{z}}_i.
\end{align*}
where $\tilde{\vct{z}}_i$ is the $i$-th column of $\mtx{Z}^\top$. Define 
\begin{align}
    \widetilde{\mtx{V}}_{\gamma} = \mtx{I}_{p} + \frac{\gamma}{p}  \mtx{Z}^\top  \mtx{Z},
\end{align}
then we can rewrite $ \left(\mtx{V}_{\gamma}^{-1} \right)_{ii}$ and $ \left(\mtx{V}_{\gamma}^{-2} \right)_{ii}$ as
\begin{align*}
     \left(\mtx{V}_{\gamma}^{-1} \right)_{ii} = 1 - \frac{\gamma}{p} \tilde{\vct{z}}_i^\top \widetilde{\mtx{V}}_{\gamma}^{-1}\tilde{\vct{z}}_i, \quad \left(\mtx{V}_{\gamma}^{-2} \right)_{ii} = 1 - \frac{\gamma}{p}\tilde{\vct{z}}_i^\top \widetilde{\mtx{V}}_{\gamma}^{-1}\tilde{\vct{z}}_i-\frac{\gamma}{p}
    \tilde{\vct{z}}_i^\top\widetilde{\mtx{V}}_{\gamma}^{-2}\tilde{\vct{z}}_i.
\end{align*}
Furthermore, by a similar argument, it can be shown that for $l = 1,2,\ldots$
\begin{equation}
\label{eq:V_tilde_l}
    \mtx{V}_{\gamma}^{-l} = \mtx{I}_n - \frac{\gamma}{p}\sum_{q=1}^l\mtx{Z} \widetilde{\mtx{V}}_{\gamma}^{-q}\mtx{Z}^\top,
\end{equation}
with
\begin{equation}
\label{eq:V_invk_ii}
     \left(\mtx{V}_{\gamma}^{-l} \right)_{ii} = 1 - \frac{\gamma}{p}\sum_{q=1}^l\tilde{\vct{z}}_i^\top \widetilde{\mtx{V}}_{\gamma}^{-q}\tilde{\vct{z}}_i.
\end{equation}
Similar to \eqref{eq:diag_homo_3} in Lemma \ref{lem:leave_one_out}, combining the leave-one-out technique and Hanson-Wright inequality, taking the uniform bound gives
\begin{align}\label{eq:tilde_zi_V_zi}
    \max_{i \in [n]} \left| \frac{1}{p} \tilde{\vct{z}}_i^\top \widetilde{\mtx{V}}_{\gamma}^{-q}\tilde{\vct{z}}_i -  \frac{1}{np}\trace\left(  \widetilde{\mtx{V}}_{\gamma}^{-q} \mtx{Z}^\top \mtx{Z}\right)  \right| = O_P\left(\sqrt{\frac{\log n}{n}} \right).
\end{align}
Together with \eqref{eq:V_tilde_l} and \eqref{eq:V_invk_ii} yields \eqref{eq:diff_Vii_tr}.\\

Note that we have 
\begin{align*}
(\mtx{B}_{\gamma})_{ii} &= \frac{1}{n}( \mtx{V}_{\gamma}^{-1} )_{ii}  - \frac{ \frac{1}{n}( \mtx{V}_{\gamma}^{-2} )_{ii}}{\frac{1}{n}\trace(\mtx{V}_{\gamma}^{-1})},
\\
(\mtx{B}_{\gamma})_{ii}^2 = & \left(\frac{1}{n} \right)^2\left( \mtx{V}_{\gamma}^{-1}\right)_{ii}^2  - \frac{2\left(\frac{1}{n} \right)^2 \left( \mtx{V}_{\gamma}^{-1} \right)_{ii}\left( \mtx{V}_{\gamma}^{-2} \right)_{ii}}{\left(\frac{1}{n}\trace(\mtx{V}_{\gamma}^{-1})\right)^2} + \frac{\left(\frac{1}{n} \right)^2 \left( \mtx{V}_{\gamma}^{-2} \right)_{ii}^2}{\left(\frac{1}{n}\trace(\mtx{V}_{\gamma}^{-1})\right)^2}.
\end{align*}
By \eqref{eq:diff_Vii_tr}, we have
\begin{align*}
\max_{i \in [n]} \left| (\mtx{B}_{\gamma})_{ii} - \frac{1}{n}\trace\left(\mtx{B}_{\gamma} \right) \right| 
&\leq  \max_{i \in [n]} \left|\frac{1}{n} \left( \mtx{V}_{\gamma}^{-1}\right)_{ii} -\frac{1}{n^2}\trace\left( \mtx{V}_{\gamma}^{-1}\right)  \right| 
\\
&~~~ + \frac{ \max_{i \in [n]} \left|\frac{1}{n} \left( \mtx{V}_{\gamma}^{-2}\right)_{ii} -\frac{1}{n^2}\trace\left( \mtx{V}_{\gamma}^{-2}\right)  \right|}{\frac{1}{n}\trace(\mtx{V}_{\gamma}^{-1})}\\
& =  O_P\left( \frac{1}{n}\sqrt{\frac{\log n}{n}}\right),
\end{align*}
and consequently
\begin{align*}
\max_{i \in [n]} \left| (\mtx{B}_{\gamma})_{ii}^2 - \left(\frac{1}{n}\trace\left(\mtx{B}_{\gamma} \right)\right)^2 \right| 
&\leq  \max_{i \in [n]} \left| (\mtx{B}_{\gamma})_{ii} - \frac{1}{n}\trace\left(\mtx{B}_{\gamma} \right) \right|^2
\\
&~~~+ 2 \frac{1}{n}\vert\trace\left(\mtx{B}_{\gamma} \right)\vert  \max_{i \in [n]} \left| (\mtx{B}_{\gamma})_{ii} - \frac{1}{n}\trace\left(\mtx{B}_{\gamma} \right) \right|
\\
&= O_P\left( \frac{1}{n^2}\sqrt{\frac{\log n}{n}}\right),
\end{align*}
which yield \eqref{eq:diff_Bii_tr} and \eqref{eq:diff_Bii_tr_sq}.\\

\subsection{Proof of Lemma \ref{lem:B_off_diag}}

From \eqref{eq:V_tilde_l}, we have
\[
\left(\mtx{V}_{\gamma}^{-l}\right)_{ij} = -\frac{\gamma}{p}\sum_{k=1}^l\tilde{\vct{z}}_i^\top \widetilde{\mtx{V}}_{\gamma}^{-k}\tilde{\vct{z}}_j.
\]
As with \eqref{eq:off_diag_rough} in Lemma \ref{lem:leave_two_out}, we can have
\[
\max_{i \neq j} \left|\tilde{\vct{z}}_i^\top \widetilde{\mtx{V}}_{\gamma}^{-k}\tilde{\vct{z}}_j\right| = O_P(\sqrt{n \log n}).
\]
Then \eqref{eq:Bij_max} follows from
\begin{align*}
    (\mtx{B}_{\gamma})_{ij} = &\frac{1}{n}\left( \mtx{V}_{\gamma}^{-1}\right)_{ij} - \frac{\frac{1}{n} \left( \mtx{V}_{\gamma}^{-2} \right)_{ij}}{\frac{1}{n}\trace(\mtx{V}_{\gamma}^{-1})}.
\end{align*}

To prove \eqref{eq:Bij_weight_sum}, we can know that
\begin{align}\label{eq:B_gamma_ij_sq}
    (\mtx{B}_{\gamma})_{ij}^2 &= \frac{1}{n^2}\left( \mtx{V}_{\gamma}^{-1}\right)_{ij}^2 -  \frac{\frac{2}{n^2} \left( \mtx{V}_{\gamma}^{-1} \right)_{ij}\left( \mtx{V}_{\gamma}^{-2} \right)_{ij}}{\frac{1}{n}\trace(\mtx{V}_{\gamma}^{-1})}  +  \frac{\frac{1}{n^2} \left( \mtx{V}_{\gamma}^{-2} \right)_{ij}^2}{\left(\frac{1}{n}\trace(\mtx{V}_{\gamma}^{-1})\right)^2} \notag\\
    & =  \frac{\gamma^2}{n^2} \frac{1}{p^2}\left( \tilde{\eta}_{ij}^{(1)} \right)^2 -2\frac{\gamma^2}{n^2} \frac{1}{p^2}\frac{\left( \tilde{\eta}_{ij}^{(1)} \right)^2 +  \tilde{\eta}_{ij}^{(1)}  \tilde{\eta}_{ij}^{(2)}  }{\frac{1}{n}\trace(\mtx{V}_{\gamma}^{-1})} +  \frac{\gamma^2}{n^2} \frac{1}{p^2} \frac{\left( \tilde{\eta}_{ij}^{(1)} \right)^2 + 2 \tilde{\eta}_{ij}^{(1)}  \tilde{\eta}_{ij}^{(2)}   + \left( \tilde{\eta}_{ij}^{(2)} \right)^2 }{\left(\frac{1}{n}\trace(\mtx{V}_{\gamma}^{-1})\right)^2},
\end{align}
where
\begin{align*}
     \tilde{\eta}_{ij}^{(k)} \coloneqq \tilde{\vct{z}}_i^\top \widetilde{\mtx{V}}_{\gamma}^{-k}\tilde{\vct{z}}_j.
\end{align*}
Define
\begin{align*}
     \bar{\theta}_2 &= \gamma^2 \tau \left(  \tilde{\kappa}_{1,1} (\gamma,\tau) - 2 \frac{ \tilde{\kappa}_{1,1} (\gamma,\tau)  +  \tilde{\kappa}_{1,2} (\gamma,\tau) }{h_1(\gamma,\tau)}  + \frac{ \tilde{\kappa}_{1,1} (\gamma,\tau) + 2  \tilde{\kappa}_{1,2} (\gamma,\tau) +  \tilde{\kappa}_{2,2} (\gamma,\tau)}{\left(h_1(\gamma,\tau)\right)^2}\right).
\end{align*}
where
\[
\tilde{\kappa}_{m,l}(\gamma,\tau) = \sum_{q_1 = 1}^l\sum_{q_2 = 1}^m \tilde{\bar{a}}_{q_1}^{(l)}(\gamma,\tau) \tilde{\bar{a}}_{q_2}^{(m)}(\gamma,\tau) \tilde{h}_{q_1+q_2}(\gamma,\tau),
\]
and
\[
\tilde{\bar{a}}_{1}^{(1)}(\gamma, \tau)=\frac{1}{\left(1+\gamma \tilde{h}_{1}(\gamma, \tau) \right)^{2}}, \quad \tilde{\bar{a}}_{1}^{(2)}(\gamma, \tau)=\frac{-2 \gamma   \tilde{h}_{2}(\gamma, \tau)}{\left(1+\gamma  \tilde{h}_{1}(\gamma, \tau) \right)^{3}}, \quad \tilde{\bar{a}}_{2}^{(2)}(\gamma, \tau)=\frac{1}{\left(1+\gamma  \tilde{h}_{1}(\gamma, \tau) \right)^{2}}.
\]
Similar to the definition of $h_l(\gamma, \tau)$, $\tilde{h}_{l}(\gamma, \tau)$ is the limit of $\frac{1}{p}\trace\left(  \widetilde{\mtx{V}}_{\gamma}^{-l}  \right)$.\\
Then similar to Proposition \ref{thm:leave_two_out}, using the leave-two-out technique, there holds that
$$
\begin{aligned}
\tilde{\eta}_{ij}^{(1)} &= \tilde{\bar{a}}_{1 ; i j}^{(1)} \tilde{\eta}_{i j ; i j}^{\left(1\right)},\\
\tilde{\eta}_{ij}^{(2)} & = \tilde{\bar{a}}_{1 ; i j}^{(2)} \tilde{\eta}_{i j ; i j}^{\left(1\right)} +  \tilde{\bar{a}}_{2 ; i j}^{(2)} \tilde{\eta}_{i j ; i j}^{\left(2\right)},
\end{aligned}
$$
with
$$
\begin{aligned}
\tilde{\bar{a}}_{1 ; i j}^{(1)}=& \frac{1}{\left(1+\frac{\gamma}{p}\tilde{\eta}_{i i ; i}^{(1)}\right)\left(1+\frac{\gamma}{p}\tilde{\eta}_{j j ; i j}^{(1)}\right)}, 
\\
\tilde{\bar{a}}_{1 ; i j}^{(2)}=& \frac{-\frac{\gamma}{p}\tilde{\eta}_{i i ; i}^{(2)}}{\left(1+\frac{\gamma}{p}\tilde{\eta}_{i i ; i}^{(1)}\right)^{2}\left(1+\frac{\gamma}{p}\tilde{\eta}_{j j ; i j}^{(1)}\right)}+\frac{-\frac{\gamma}{p}\tilde{\eta}_{j j ; i j}^{(2)}}{\left(1+\frac{\gamma}{p}\tilde{\eta}_{i i ; i}^{(1)}\right)\left(1+\frac{\gamma}{p}\tilde{\eta}_{j j ; i j}^{(1)}\right)^{2}}, \\
\tilde{\bar{a}}_{2 ; i j}^{(2)}=& \tilde{\bar{a}}_{1 ; i j}^{(1)}.
\end{aligned}
$$
And
\[
\max _{1 \leq i \neq j \leq p} \max _{1 \leq l \leq 2} \max _{1 \leq q_{1} \leq l}\left|\tilde{\bar{a}}_{q_{1} ; i j}^{(l)}-\tilde{\bar{a}}_{q_{1}}^{(l)}(\gamma, \tau)\right|=O_{\mathrm{P}}\left(\sqrt{\frac{\log p}{n}}\right) .
\]
Furthermore,
\begin{align}\label{eq:diff_aij_a_tilde}
   &~~~\max _{1 \leq i \neq j \leq p} \max _{1 \leq l,m \leq 2} \max _{\substack{1 \leq q_{1} \leq l\\1 \leq q_2 \leq m}}\left|\tilde{\bar{a}}_{q_{1} ; i j}^{(l)}\tilde{\bar{a}}_{q_{2} ; i j}^{(m)}-\tilde{\bar{a}}_{q_{1}}^{(l)}(\gamma, \tau)\tilde{\bar{a}}_{q_{2}}^{(m)}(\gamma, \tau)\right|=O_{\mathrm{P}}\left(\sqrt{\frac{\log p}{n}}\right) 
\end{align}
Similar to \eqref{eq:sum_eta_ij_ij_leave_two_h} and \eqref{eq:sum_eta_ij_ij_leave_two_h_beta} , we can have that
\begin{align*}
    \frac{1}{n(n-1)}  \sum_{i \neq j} \frac{1}{p} \tilde{\eta}_{ij;ij}^{(q_1)}  \tilde{\eta}_{ij;ij}^{(q_2)} = \frac{1}{p} \trace\left(   \widetilde{\mtx{V}}_{\gamma}^{-(q_1 + q_2)} \right) + o_P(1),
\end{align*}
and
\begin{align*}
    \frac{1}{n^2}  \sum_{i \neq j}\varepsilon_i^2 \varepsilon_j^2 \frac{1}{p} \tilde{\eta}_{ij;ij}^{(q_1)}  \tilde{\eta}_{ij;ij}^{(q_2)} = \sigma_0^4\frac{1}{p} \trace\left(   \widetilde{\mtx{V}}_{\gamma}^{-(q_1 + q_2)} \right) + o_P(1),
\end{align*}
by \eqref{eq:noise_var_concent} in Lemma \ref{lem:noise_var_bound}. Then similar to the proof of Lemma \ref{lem:leave_two_out}, \eqref{eq:Bij_weight_sum} can be obtained from \eqref{eq:B_gamma_ij_sq}.

\subsection{Proof of Lemma \ref{lem:V_inv_diag_loo}}
For any $k = 1,\cdots, p$, denote $\mtx{Z}_{-k} = [\vct{z}_1,\cdots,\vct{z}_{k-1},\vct{z}_{k+1},\cdots]$, then 
\begin{align*}
\mtx{V}_{\gamma,-k} = \mtx{V}_\gamma - \frac{\gamma}{p} \vct{z}_k \vct{z}_k^\top = \mtx{I}_n +   \frac{\gamma}{p} \mtx{Z}_{-k} \mtx{Z}_{-k}^\top.
\end{align*}
Similar to the proof of \eqref{eq:diff_Vii_tr} in Lemma \ref{lem:V_inv_diag}, we can define
\begin{align}\label{eq:V_tilde_l_k}
    \widetilde{\mtx{V}}_{\gamma,-k} = \mtx{I}_{p-1} + \frac{\gamma}{p} \mtx{Z}_{-k}^\top \mtx{Z}_{-k},
\end{align}
and $\tilde{\vct{z}}_{i,-k}^\top$ is the i-th row of $\mtx{Z}_{-k}$. Then we can have that
\begin{align}\label{eq:V_invk_ii_k}
    \left(  \mtx{V}_{\gamma,-k}^{-l} \right)_{ii} = 1 - \frac{\gamma}{p} \sum_{q=1}^l \tilde{\vct{z}}_{i,-k}^\top \widetilde{\mtx{V}}_{\gamma,-k}^{-q}\tilde{\vct{z}}_{i,-k}.
\end{align}
Again, similar to \eqref{eq:diag_homo_3} in Lemma \ref{lem:leave_one_out}, combining the leave-one-out technique and Hanson-Wright inequality (taking $t = \sqrt{ \frac{ 3 \log p}{c}} \| \widetilde{\mtx{V}}_{\gamma,-ik}^{-q} \|_F$), taking the uniform bound gives
\begin{align}\label{eq:tilde_zi_V_zi_k}
    \max_{k \in [p]} \max_{i \in [n]} \left| \frac{1}{p} \tilde{\vct{z}}_{i,-k}^\top \widetilde{\mtx{V}}_{\gamma,-k}^{-q}\tilde{\vct{z}}_{i,-k} -  \frac{1}{np}\trace\left(  \widetilde{\mtx{V}}_{\gamma,-k}^{-q} \mtx{Z}_{-k}^\top \mtx{Z}_{-k}\right)  \right| = O_P\left(\sqrt{\frac{\log n}{n}} \right).
\end{align}
Then \eqref{eq:V_tilde_l} and \eqref{eq:V_invk_ii} implies
\begin{equation}
\label{eq:diff_Vii_tr_k}
 \max_{k \in [p]} \max_{i \in [n]} \left|   \left(\mtx{V}_{\gamma,-k}^{-l} \right)_{ii}  - \frac{1}{n} \trace\left(  \mtx{V}_{\gamma,-k}^{-l} \right)  \right| = O_P\left(\sqrt{\frac{\log n}{n}} \right), \quad l = 1,2,3,4.
\end{equation}
By \eqref{eq:diff_Vinv_Vinvk} in Lemma \ref{lem:leave_one_out} and \eqref{eq:diff_Vii_tr_k} we can get \eqref{eq:diff_Vii_tr_loo}.

\subsection{Proof of Lemma \ref{lem:conditional_variance_formula}}
In this section, we focus on the conditional variance $\Var[\Delta(\gamma) \vert \mtx{Z}, \vct{\varepsilon} ]$. 

With $\vct{y}$ defined in \eqref{eq:fixed_effects_rad_Zeps}, write
\begin{align*}
    \Delta(\gamma)
    &=\underbrace{ \vct{\xi}^\top \mtx{A}_\gamma \vct{\xi}}_{M_1}
    +\underbrace{2\vct{\xi}^\top \mtx{C}_\gamma \vct{\zeta}}_{M_2}
    +\underbrace{\vct{\zeta}^\top \mtx{D}_\gamma \vct{\zeta}}_{M_3},
\end{align*}
where
\[
\mtx{A}_\gamma
=\mtx{\Lambda}_\beta \mtx{Z}^\top \mtx{B}_\gamma \mtx{Z}\mtx{\Lambda}_\beta,
\qquad
\mtx{C}_\gamma
=\mtx{\Lambda}_\beta \mtx{Z}^\top \mtx{B}_\gamma \mtx{\Lambda}_\varepsilon,
\qquad
\mtx{D}_\gamma
=\mtx{\Lambda}_\varepsilon \mtx{B}_\gamma \mtx{\Lambda}_\varepsilon.
\]
Conditionally on $(\mtx{Z},\vct{\varepsilon})$, the vectors $\vct{\xi}$ and $\vct{\zeta}$ are independent Rademacher vectors. Hence
\[
\E[M_1\mid \mtx{Z},\vct{\varepsilon}]=\trace(\mtx{A}_\gamma),
\qquad
\E[M_2\mid \mtx{Z},\vct{\varepsilon}]=0,
\qquad
\E[M_3\mid \mtx{Z},\vct{\varepsilon}]=\trace(\mtx{D}_\gamma),
\]
and therefore
\[
\E[\Delta(\gamma)\mid \mtx{Z},\vct{\varepsilon}]
=\trace(\mtx{A}_\gamma)+\trace(\mtx{D}_\gamma)
=\sum_{k=1}^p \beta_k^2\vct{z}_k^\top \mtx{B}_\gamma\vct{z}_k
+\trace(\mtx{\Lambda}_\varepsilon^2\mtx{B}_\gamma),
\]
which is exactly $\widetilde{\Delta}_*(\gamma)$.

We next use the following elementary identity: if $\vct{r}$ is a vector of independent Rademacher variables and $\mtx{Q}$ is a deterministic symmetric matrix, then
\begin{equation}
\label{eq:rad_quad_variance_identity_revision}
\Var(\vct{r}^\top \mtx{Q}\vct{r})
=2\sum_{i\neq j} Q_{ij}^2.
\end{equation}
Moreover, if $\vct{r}$ and $\vct{s}$ are independent Rademacher vectors, then for any deterministic matrix $\mtx{C}$,
\begin{equation}
\label{eq:rad_bilinear_variance_identity_revision}
\Var(2\vct{r}^\top\mtx{C}\vct{s})=4\|\mtx{C}\|_F^2.
\end{equation}
The cross-covariances between $M_2$ and $M_1$, between $M_2$ and $M_3$, and between $M_1$ and $M_3$ are zero: the first two vanish by odd Rademacher moments, while the last one vanishes because $M_1$ and $M_3$ depend on independent Rademacher vectors. Thus
\begin{align*}
\Var[\Delta(\gamma)\mid \mtx{Z},\vct{\varepsilon}]
&=\Var(M_1\mid \mtx{Z},\vct{\varepsilon})
+\Var(M_2\mid \mtx{Z},\vct{\varepsilon})
+\Var(M_3\mid \mtx{Z},\vct{\varepsilon})  \\
&=2\sum_{k\neq j}(\mtx{A}_\gamma)_{kj}^2
+4\|\mtx{C}_\gamma\|_F^2
+2\sum_{i\neq j}(\mtx{D}_\gamma)_{ij}^2.
\end{align*}
Expanding the three terms gives
\[
2\sum_{k\neq j}(\mtx{A}_\gamma)_{kj}^2
=2\sum_{k\neq j}\beta_k^2\beta_j^2
\left(\vct{z}_k^\top\mtx{B}_\gamma\vct{z}_j\right)^2,
\]
\[
4\|\mtx{C}_\gamma\|_F^2
=4\sum_{k=1}^p\beta_k^2\vct{z}_k^\top
\mtx{B}_\gamma\mtx{\Lambda}_{\varepsilon}^2\mtx{B}_\gamma\vct{z}_k,
\]
and
\[
2\sum_{i\neq j}(\mtx{D}_\gamma)_{ij}^2
=2\sum_{i\neq j}\varepsilon_i^2\varepsilon_j^2(\mtx{B}_\gamma)_{ij}^2.
\]
Multiplying by $n$ yields
\begin{align*}
    &~~~\Var\left[\sqrt{n}(\Delta(\gamma) ) \vert \mtx{Z},\vct{\varepsilon} \right] \\
    &= 2n \sum_{1 \leq k\neq j \leq p} \beta_k^2 \beta_j^2 \left(\vct{z}_k^\top \mtx{B}_{\gamma} \vct{z}_j \right)^2
    + 4 n \sum_{k=1}^p \beta_k^2 \vct{z}_k^\top \mtx{B}_{\gamma}\mtx{\Lambda}_{\varepsilon}^2\mtx{B}_{\gamma} \vct{z}_k
    + 2n \sum_{1 \leq k \neq j \leq n} \varepsilon_k^2 \varepsilon_j^2 (\mtx{B}_{\gamma})_{kj}^2,
\end{align*}
which proves the stated conditional variance formula.

\subsection{Proof of Lemma \ref{lem:Delta_starstar}}
Let
\[
H_{1,n}(\gamma)=\trace(\mtx{V}_\gamma^{-1}),
\qquad
H_{2,n}(\gamma)=\trace(\mtx{V}_\gamma^{-2}).
\]
Since
\[
\mtx{V}_{\gamma_0}
=
\mtx{I}_n+\gamma_0 p^{-1}\mtx{Z}\mtx{Z}^\top
=
\frac{\gamma_0}{\gamma}\mtx{V}_\gamma
+\left(1-\frac{\gamma_0}{\gamma}\right)\mtx{I}_n,
\]
we can evaluate $\Delta_{**}(\gamma)$ without treating the
$\mtx{Z}\mtx{Z}^\top$ traces separately:
\begin{align*}
\frac{\Delta_{**}(\gamma)}{\sigma_0^2}
&=
\trace\left\{\left(\frac{\mtx{V}_\gamma^{-1}}{n}
-
\frac{\mtx{V}_\gamma^{-2}}{H_{1,n}(\gamma)}\right)
\left(\frac{\gamma_0}{\gamma}\mtx{V}_\gamma
+\left(1-\frac{\gamma_0}{\gamma}\right)\mtx{I}_n\right)\right\} \\
&=
\left(1-\frac{\gamma_0}{\gamma}\right)
\left\{\frac{H_{1,n}(\gamma)}{n}
-
\frac{H_{2,n}(\gamma)}{H_{1,n}(\gamma)}\right\} \\
&=
\left(\frac{\gamma_0}{\gamma}-1\right)
\left\{\frac{H_{2,n}(\gamma)}{H_{1,n}(\gamma)}
-
\frac{H_{1,n}(\gamma)}{n}\right\}.
\end{align*}
By Corollary \ref{thm:MP_trace} and the definition of $h_k$ in \eqref{eq:h_MP},
\[
\frac{H_{1,n}(\gamma)}{n}\stackrel{a.s.}{\longrightarrow} h_1(\gamma,\tau),
\qquad
\frac{H_{2,n}(\gamma)}{n}\stackrel{a.s.}{\longrightarrow} h_2(\gamma,\tau).
\]
Therefore
\[
\Delta_{**}(\gamma)
\stackrel{a.s.}{\longrightarrow}
c_\gamma
\coloneqq
\sigma_0^2
\left(\frac{\gamma_0}{\gamma}-1\right)
\frac{h_2(\gamma,\tau)-h_1^2(\gamma,\tau)}{h_1(\gamma,\tau)}.
\]
It remains only to check the sign of the last factor. Let $\mu_\tau$ be the Mar$\check{c}$enko-Pastur probability measure appearing in \eqref{eq:h_MP}, including the possible atom at zero when $\tau>1$, and set $g_\gamma(x)=(1+\gamma x)^{-1}$. Then
\[
h_1(\gamma,\tau)=\int g_\gamma(x)\,\mathrm{d}\mu_\tau(x),
\qquad
h_2(\gamma,\tau)=\int g_\gamma^2(x)\,\mathrm{d}\mu_\tau(x).
\]
Hence
\[
h_2(\gamma,\tau)-h_1^2(\gamma,\tau)
=
\Var_{\mu_\tau}\{g_\gamma(X)\}>0,
\]
because $\mu_\tau$ has nonzero continuous mass on $[b_-(\tau),b_+(\tau)]$ and $g_\gamma$ is not constant there for $\gamma>0$. Also $h_1(\gamma,\tau)>0$. Thus the sign of $c_\gamma$ is exactly the sign of $\gamma_0/\gamma-1$, proving the claim.

\subsection{Proof of Lemma \ref{lem:Taylor_approx}}
Recall that $\Delta(\gamma) = \vct{y}^\top \mtx{B}_\gamma \vct{y}$, where
\[
\mtx{B}_\gamma
=
\frac{\mtx{V}_\gamma^{-1}}{n}
-
\frac{\mtx{V}_\gamma^{-2}}{\trace(\mtx{V}_\gamma^{-1})}.
\]
For $l=1,2,\ldots$,
\[
\frac{\mathrm{d}}{\mathrm{d}\gamma} \mtx{V}_\gamma^{-l} = -l\mtx{V}_\gamma^{-(l+1)} \left(\frac{1}{p}\mtx{Z}\mtx{Z}^\top \right) = - \frac{l}{\gamma} \left( \mtx{V}_\gamma^{-l} - \mtx{V}_\gamma^{-(l+1)}\right),
\]
and hence
\[
\frac{\mathrm{d}}{\mathrm{d}\gamma} \trace(\mtx{V}_\gamma^{-1}) = - \frac{1}{\gamma}\trace\left( \mtx{V}_\gamma^{-1}-\mtx{V}_\gamma^{-2}  \right),
\]
so
\begin{align}\label{eq:B_1st_der}
    \frac{\mathrm{d}}{\mathrm{d}\gamma} \mtx{B}_\gamma &=-\frac{1}{\gamma} \left(\frac{\mtx{V}_\gamma^{-1} -  \mtx{V}_\gamma^{-2}  }{n} + \frac{-2\mtx{V}_\gamma^{-2} + 2\mtx{V}_\gamma^{-3} }{\trace(\mtx{V}_\gamma^{-1})} + \frac{\mtx{V}_\gamma^{-2}\trace\left( \mtx{V}_\gamma^{-1} - \mtx{V}_\gamma^{-2}  \right) }{\left( \trace(\mtx{V}_\gamma^{-1})\right)^2 } \right).
\end{align}
The deterministic equivalent used in the proof of Theorem \ref{thm:consistency} gives, for every fixed $\gamma>0$ and $l=1,2,\ldots$,
\begin{align}\label{eq:y_V_y}
\frac{1}{n} \vct{y}^\top \mtx{V}_{\gamma}^{-l}\vct{y}
=
\frac{\sigma_0^2}{n}\trace\left( \mtx{V}_{\gamma}^{-l} \mtx{V}_{\gamma_0} \right)+o_P(1),
\end{align}
while $\mtx{V}_{\gamma_0}=(\gamma_0/\gamma)\mtx{V}_\gamma-(\gamma_0/\gamma-1)\mtx{I}_n$ implies
\begin{align}\label{eq:y_V_y_tr}
\frac{\sigma_0^2}{n}\trace\left( \mtx{V}_{\gamma}^{-l} \mtx{V}_{\gamma_0} \right)
\stackrel{a.s.}{\longrightarrow}
\sigma_0^2 \frac{\gamma_0}{\gamma} h_{l-1}(\gamma,\tau)
-
\sigma_0^2\left( \frac{\gamma_0}{\gamma} - 1 \right)h_l(\gamma,\tau),
\end{align}
where $h_0(\gamma,\tau)=1$. In particular, at $\gamma=\gamma_0$,
\begin{align}\label{eq:y_V0_y}
\frac{1}{n} \vct{y}^\top \mtx{V}_{\gamma_0}^{-l}\vct{y}  \stackrel{P}{\longrightarrow}\sigma_0^2h_{l-1} (\gamma_0,\tau).
\end{align}
Combining this with \eqref{eq:B_1st_der} and $n^{-1}\trace(\mtx{V}_{\gamma_0}^{-l})\to h_l(\gamma_0,\tau)$, and writing $h_l=h_l(\gamma_0,\tau)$, we obtain
\[
\Delta'(\gamma_0)
\stackrel{P}{\longrightarrow}
-\frac{\sigma_0^2}{\gamma_0}
\left\{
(1-h_1)-\frac{h_1-h_2}{h_1}
\right\}
=
\frac{\sigma_0^2}{\gamma_0}\frac{h_1^2-h_2}{h_1}.
\]
Thus
\begin{align}\label{eq:Delta_inf_prime}
    \Delta'_\infty(\gamma_0) =  \frac{\sigma_0^2}{\gamma_0} \frac{h_1^2(\gamma_0,\tau) - h_2(\gamma_0,\tau)}{h_1(\gamma_0,\tau)}.
\end{align}
It remains to control the Taylor remainder. Let
$K=[\gamma_0/2,3\gamma_0/2]$. Since $\hat\gamma\stackrel{P}{\to}\gamma_0$, the intermediate point $\gamma_\delta$ between $\gamma_0$ and $\hat\gamma$ lies in $K$ with probability tending to one. On this event, repeated differentiation of $\mtx{B}_\gamma$ shows that every term of $\mtx{B}_\gamma''$ is a finite product of powers of $\mtx{V}_\gamma^{-1}$ multiplied by coefficients of order $O_P(n^{-1})$. Indeed, uniformly over $\gamma\in K$, $\|\mtx{V}_\gamma^{-r}\|\leq 1$, $\trace(\mtx{V}_\gamma^{-1})$ is bounded above and below by constant multiples of $n$ with probability tending to one, and derivatives of $\trace(\mtx{V}_\gamma^{-r})$ are $O_P(n)$ for each fixed $r$. Hence
\begin{align}\label{eq:delta_2md_der}
    \sup_{\gamma\in K}\left\|\frac{\mathrm{d}^2}{\mathrm{d}\gamma^2}\mtx{B}_\gamma\right\| = O_P(n^{-1}).
\end{align}
Since $n^{-1}\|\vct y\|^2=O_P(1)$, \eqref{eq:delta_2md_der} implies
\[
\Delta''(\gamma_\delta)
=
\vct y^\top \mtx{B}_{\gamma_\delta}''\vct y
=O_P(1).
\]
Taylor's expansion at $\gamma_0$ gives
\[
0=\Delta(\hat\gamma)
=
\Delta(\gamma_0)
+(\hat\gamma-\gamma_0)\Delta'(\gamma_0)
+\frac{1}{2}(\hat\gamma-\gamma_0)^2\Delta''(\gamma_\delta),
\]
and therefore
\begin{align}\label{eq:diff_gamma0_taylor_rm}
   \sqrt{n}\left( \hat{\gamma} - \gamma_0 \right)  = - \frac{ \sqrt{n}\Delta(\gamma_0)}{\Delta'(\gamma_0) + \frac{1}{2}\left( \hat{\gamma} - \gamma_0 \right)\Delta''(\gamma_\delta)}.
\end{align}
Using $\hat\gamma-\gamma_0=o_P(1)$ and $\Delta''(\gamma_\delta)=O_P(1)$, the denominator satisfies
\begin{align}\label{eq:taylor_rm}
\Delta'(\gamma_0) + \frac{1}{2}\left( \hat{\gamma} - \gamma_0 \right)\Delta''(\gamma_\delta) = \Delta_{\infty}'(\gamma_0) + o_P(1).
\end{align}
Substituting \eqref{eq:taylor_rm} into \eqref{eq:diff_gamma0_taylor_rm} yields
\begin{align}\label{eq:diff_gamma0}
     \sqrt{n}\left( \hat{\gamma} - \gamma_0 \right)  = - \frac{\sqrt{n} \Delta(\gamma_0)  }{\Delta_{\infty}'(\gamma_0) } + o_P(1),
\end{align}
where $\Delta_{\infty}'(\gamma_0) $ is defined in \eqref{eq:Delta_inf_prime}.

\subsection{Proof of Lemma \ref{lem:noise_diag_approx_fine}}

We prove the result for fixed $\gamma>0$. Since
\[
(\mtx B_\gamma)_{ii}
=
\frac{1}{n}(\mtx V_\gamma^{-1})_{ii}
-
\frac{
\frac{1}{n}(\mtx V_\gamma^{-2})_{ii}
}{
\frac{1}{n}\trace(\mtx V_\gamma^{-1})
},
\]
it is enough to prove, for $l=1,2$,
\[
\sum_{i=1}^n
(\varepsilon_i^2-\sigma_0^2)
\left\{
\frac{1}{n}(\mtx V_\gamma^{-l})_{ii}
-
\frac{1}{n^2}\trace(\mtx V_\gamma^{-l})
\right\}
=
o_P(n^{-1/2}).
\]
We prove this display.

Let $\tilde{\vct z}_i^\top$ be the $i$th row of $\mtx Z$ and define
\[
\widetilde{\mtx V}_\gamma
=
\mtx I_p+\frac{\gamma}{p}\mtx Z^\top\mtx Z.
\]
The identity
\[
\mtx V_\gamma^{-l}
=
\mtx I_n
-
\frac{\gamma}{p}\sum_{q=1}^l
\mtx Z\widetilde{\mtx V}_\gamma^{-q}\mtx Z^\top
\]
implies
\[
(\mtx V_\gamma^{-l})_{ii}
-
\frac{1}{n}\trace(\mtx V_\gamma^{-l})
=
-
\frac{\gamma}{p}\sum_{q=1}^l
\left\{
\tilde{\vct z}_i^\top\widetilde{\mtx V}_\gamma^{-q}\tilde{\vct z}_i
-
\frac{1}{n}\trace(\widetilde{\mtx V}_\gamma^{-q}\mtx Z^\top\mtx Z)
\right\}.
\]
Thus it suffices to show that, for each fixed $q$,
\[
T_q
\coloneqq
\sum_{i=1}^n
\frac{\varepsilon_i^2-\sigma_0^2}{n}
\frac{1}{p}
\left\{
\tilde{\vct z}_i^\top\widetilde{\mtx V}_\gamma^{-q}\tilde{\vct z}_i
-
\frac{1}{n}\trace(\widetilde{\mtx V}_\gamma^{-q}\mtx Z^\top\mtx Z)
\right\}
=
o_P(n^{-1/2}).
\]

Condition on $\vct\varepsilon$. The weights
\[
a_i=\frac{\varepsilon_i^2-\sigma_0^2}{n}
\]
are then deterministic. In the independent-noise case,
\[
\sum_{i=1}^n a_i^2
=
\frac{1}{n^2}\sum_{i=1}^n(\varepsilon_i^2-\sigma_0^2)^2
=O_P(n^{-1}).
\]
The same bound holds under the correlated Gaussian assumptions, since
$\max_i\varepsilon_i^2=O_P(\log n)$ and
$n^{-1}\sum_i\varepsilon_i^2=O_P(1)$. In both cases,
$\sum_i |a_i|=O_P(1)$.

Now apply the row version of Lemma \ref{lem:diag_homo_moment}. Namely, for
the row leave-one-out matrices obtained by deleting row $i$, the centered
quadratic terms
\[
D_i^{(q)}
=
\frac{1}{p}
\left\{
\tilde{\vct z}_i^\top\widetilde{\mtx V}_\gamma^{-q}\tilde{\vct z}_i
-
\frac{1}{n}\trace(\widetilde{\mtx V}_\gamma^{-q}\mtx Z^\top\mtx Z)
\right\}
\]
satisfy
\[
\max_i \E[(D_i^{(q)})^2]\leq \frac{C}{n},
\qquad
\max_{i\neq j}|\E[D_i^{(q)}D_j^{(q)}]|\leq \frac{C}{n^{3/2}}.
\]
Therefore, conditionally on $\vct\varepsilon$,
\[
\E_{\mtx Z}T_q^2
\leq
\frac{C}{n}\sum_{i=1}^n a_i^2
+
\frac{C}{n^{3/2}}\left(\sum_{i=1}^n |a_i|\right)^2.
\]
Because $\sum_i a_i^2=O_P(n^{-1})$ and
$\sum_i |a_i|=O_P(1)$, the right-hand side is
$O_P(n^{-2})+O_P(n^{-3/2})=o_P(n^{-1})$.
Thus
\[
T_q=o_P(n^{-1/2}).
\]
Summing over the finite set $q=1,\ldots,l$ proves the desired claim and hence
\eqref{eq:noise_diag_approx_fine}.

\section{MoM Interval Formulas}
\label{app:mom_derivation}

This appendix records the explicit MoM confidence intervals used in Section \ref{sec:experiments} and the algebra behind their covariance calibrations. The derivation is conditional on the realized design matrix and is intended only to document the finite-sample calibration used for the simulation benchmark.

Let $A_1=\mtx{I}_n$ and $A_2=K=p^{-1}\mtx{Z}\mtx{Z}^{\top}$, so that
\[
m_1=\frac{1}{n}\vct{y}^{\top}A_1\vct{y},
\qquad
m_2=\frac{1}{n}\vct{y}^{\top}A_2\vct{y}.
\]
Let $c_1=n^{-1}\trace(K)$ and $c_2=n^{-1}\trace(K^2)$.
The MoM point estimate is the smooth transformation
\[
\hat{\eta}_{\mathrm{MoM}}=\frac{m_2-c_1m_1}{c_2-c_1},
\qquad
\hat{\sigma}_{\mathrm{MoM}}^2=m_1-\hat{\eta}_{\mathrm{MoM}},
\qquad
\hat{\gamma}_{\mathrm{MoM}}=\frac{\hat{\eta}_{\mathrm{MoM}}}{\hat{\sigma}_{\mathrm{MoM}}^2}.
\]
Writing $\vct{m}=(m_1,m_2)^\top$, the Jacobian of $(\hat{\eta}_{\mathrm{MoM}},\hat{\sigma}^2_{\mathrm{MoM}})$ with respect to $\vct{m}$ is
\[
\mtx{J}=
\begin{pmatrix}
-c_1/(c_2-c_1) & 1/(c_2-c_1)\\
1+c_1/(c_2-c_1) & -1/(c_2-c_1)
\end{pmatrix}.
\]
Since the gradient of $\eta/\sigma^2$ with respect to $(\eta,\sigma^2)$ is
\[
\vct{g}=
\begin{pmatrix}
1/\hat{\sigma}^2_{\mathrm{MoM}}\\
-\hat{\eta}_{\mathrm{MoM}}/\hat{\sigma}^4_{\mathrm{MoM}}
\end{pmatrix},
\]
the gradient with respect to $(m_1,m_2)$ is $\hat{\vct{\ell}}=\mtx{J}^{\top}\vct{g}$.

Equivalently, writing
\[
\hat{\eta}=\hat{\eta}_{\mathrm{MoM}},\qquad
\hat{\eta}_{+}=\max\{\hat{\eta},0\},\qquad
\hat{v}=\hat{\sigma}^2_{\mathrm{MoM}},\qquad
T_j=\trace(K^j),\qquad
k_i^{(j)}=(K^j)_{ii},
\]
and $d=c_2-c_1$, the plug-in gradient is
\[
\hat{\vct{\ell}}
=
\begin{pmatrix}
\displaystyle
-\frac{c_1}{d\hat{v}}
-\left(1+\frac{c_1}{d}\right)
\frac{\hat{\eta}}{\hat{v}^2}
\\[1.2em]
\displaystyle
\frac{1}{d}
\left(
\frac{1}{\hat{v}}
+\frac{\hat{\eta}}{\hat{v}^2}
\right)
\end{pmatrix}.
\]
For any symmetric $2\times 2$ covariance matrix with entries $(\omega_{11},\omega_{12},\omega_{22})$, define
\[
\mathcal{I}(\omega_{11},\omega_{12},\omega_{22})
=
\hat{\gamma}_{\mathrm{MoM}}
\pm
z_{0.975}
\left(\hat{\ell}_1^2\omega_{11}+2\hat{\ell}_1\hat{\ell}_2\omega_{12}
+\hat{\ell}_2^2\omega_{22}\right)^{1/2}.
\]
The homo-Gauss MoM interval is
\[
\mathrm{CI}^{\mathrm{homo\text{-}Gaussian}}_{\mathrm{MoM}}
=
\mathcal{I}(\omega_{11}^{\mathrm{hG}},\omega_{12}^{\mathrm{hG}},\omega_{22}^{\mathrm{hG}}),
\]
where
\[
\omega_{11}^{\mathrm{hG}}
=
\frac{2}{n^2}\left(\hat{\eta}_{+}^2T_2+2\hat{\eta}_{+}\hat{v}T_1+n\hat{v}^2\right),
\qquad
\omega_{12}^{\mathrm{hG}}
=
\frac{2}{n^2}\left(\hat{\eta}_{+}^2T_3+2\hat{\eta}_{+}\hat{v}T_2+\hat{v}^2T_1\right),
\]
\[
\omega_{22}^{\mathrm{hG}}
=
\frac{2}{n^2}\left(\hat{\eta}_{+}^2T_4+2\hat{\eta}_{+}\hat{v}T_3+\hat{v}^2T_2\right).
\]
With $(a_1,\ldots,a_n)$ denoting the normalized variance-profile shape used in the simulation, $n^{-1}\sum_{i=1}^n a_i=1$, the hetero-Gauss MoM interval is
\[
\mathrm{CI}^{\mathrm{hetero\text{-}Gaussian}}_{\mathrm{MoM}}
=
\mathcal{I}(\omega_{11}^{\mathrm{hetG}},\omega_{12}^{\mathrm{hetG}},\omega_{22}^{\mathrm{hetG}}),
\]
where
\[
\omega_{11}^{\mathrm{hetG}}
=
\frac{2}{n^2}\left\{
\hat{\eta}_{+}^2T_2
+2\hat{\eta}_{+}\hat{v}\sum_{i=1}^n a_i k_i^{(1)}
+\hat{v}^2\sum_{i=1}^n a_i^2
\right\},
\]
\[
\omega_{12}^{\mathrm{hetG}}
=
\frac{2}{n^2}\left\{
\hat{\eta}_{+}^2T_3
+2\hat{\eta}_{+}\hat{v}\sum_{i=1}^n a_i k_i^{(2)}
+\hat{v}^2\sum_{i=1}^n a_i^2 k_i^{(1)}
\right\},
\]
\[
\omega_{22}^{\mathrm{hetG}}
=
\frac{2}{n^2}\left\{
\hat{\eta}_{+}^2T_4
+2\hat{\eta}_{+}\hat{v}\sum_{i=1}^n a_i k_i^{(3)}
+\hat{v}^2\sum_{i,j=1}^n a_i a_j K_{ij}^2
\right\}.
\]
For the homo-non-Gauss row, set
\[
\hat{\kappa}^{\mathrm{MoM}}_{\varepsilon,\mathrm{H}}
=
\frac{1}{2}\left\{
\frac{1}{n\hat{v}^2}\sum_{i=1}^n y_i^4
-3\hat{\gamma}_{\mathrm{MoM}}^2-6\hat{\gamma}_{\mathrm{MoM}}-1
\right\},
\qquad
\hat{\delta}_{\varepsilon}
=
2\left(\max\{\hat{\kappa}^{\mathrm{MoM}}_{\varepsilon,\mathrm{H}},0\}-1\right)
\hat{v}^2.
\]
The homo-non-Gauss MoM interval is
\[
\mathrm{CI}^{\mathrm{homo\text{-}Non\text{-}Gaussian}}_{\mathrm{MoM}}
=
\mathcal{I}(\omega_{11}^{\mathrm{hNG}},\omega_{12}^{\mathrm{hNG}},\omega_{22}^{\mathrm{hNG}}),
\]
where
\[
\omega_{11}^{\mathrm{hNG}}
=
\omega_{11}^{\mathrm{hG}}+\frac{\hat{\delta}_{\varepsilon}}{n},
\qquad
\omega_{12}^{\mathrm{hNG}}
=
\omega_{12}^{\mathrm{hG}}+\frac{\hat{\delta}_{\varepsilon}T_1}{n^2},
\qquad
\omega_{22}^{\mathrm{hNG}}
=
\omega_{22}^{\mathrm{hG}}+\frac{\hat{\delta}_{\varepsilon}}{n^2}\sum_{i=1}^n (k_i^{(1)})^2.
\]

We next compute the covariance of the moment vector. Under a Gaussian working model
\[
\vct{y}\mid \mtx{Z}\sim \mathcal{N}\{\vct{0},\mtx{\Sigma}_y(\mtx{D})\},
\qquad
\mtx{\Sigma}_y(\mtx{D})=\eta K+\mtx{D},
\]
the standard covariance identity for quadratic forms gives, for symmetric matrices $A$ and $B$,
\[
\Cov\left(\vct{y}^{\top}A\vct{y},\vct{y}^{\top}B\vct{y}\mid \mtx{Z}\right)
=
2\trace(A\mtx{\Sigma}_y(\mtx{D})B\mtx{\Sigma}_y(\mtx{D})).
\]
Dividing by $n^2$ and taking $(A,B)$ over $(A_1,A_1)$, $(A_1,A_2)$, and $(A_2,A_2)$ yields
\[
\mtx{\Omega}_{\mathrm{G}}(\mtx{D})
=
\frac{2}{n^2}
\begin{pmatrix}
\trace\{\mtx{\Sigma}_y(\mtx{D})^2\}
&
\trace\{\mtx{\Sigma}_y(\mtx{D})K\mtx{\Sigma}_y(\mtx{D})\}\\
\trace\{\mtx{\Sigma}_y(\mtx{D})K\mtx{\Sigma}_y(\mtx{D})\}
&
\trace\{K\mtx{\Sigma}_y(\mtx{D})K\mtx{\Sigma}_y(\mtx{D})\}
\end{pmatrix}.
\]
The homo-Gauss and hetero-Gauss covariance estimators in Section \ref{sec:experiments} are obtained from this formula after plugging in $\eta=\hat{\eta}_{+}$ and taking
\[
\mtx{D}=\hat{\sigma}^2_{\mathrm{MoM}}\mtx{I}_n
\qquad\text{or}\qquad
\mtx{D}=\hat{\sigma}^2_{\mathrm{MoM}}\Diag(a_1,\ldots,a_n),
\quad n^{-1}\sum_{i=1}^n a_i=1.
\]

For the homo-non-Gauss calibration, suppose the noise coordinates are independent and homogeneous with variance $\sigma^2$. Let
\[
\kappa_{\varepsilon}=\frac{\Var(\varepsilon_i^2)}{2\sigma^4}.
\]
The fourth cumulant of a noise coordinate is
\[
\E[\varepsilon_i^4]-3\sigma^4
=
\Var(\varepsilon_i^2)-2\sigma^4
=
2(\kappa_{\varepsilon}-1)\sigma^4
\eqqcolon \delta_{\varepsilon}.
\]
Adding this independent-coordinate fourth-cumulant contribution to the Gaussian quadratic-form covariance gives
\[
\Cov\left(\vct{y}^{\top}A\vct{y},\vct{y}^{\top}B\vct{y}\mid \mtx{Z}\right)
=
2\trace(A\mtx{\Sigma}_y B\mtx{\Sigma}_y)
+
\delta_{\varepsilon}\sum_{i=1}^n A_{ii}B_{ii}.
\]
With $A,B\in\{\mtx{I}_n,K\}$ and after division by $n^2$, the additional covariance matrix is
\[
\frac{\delta_{\varepsilon}}{n^2}
\begin{pmatrix}
n & \trace(K)\\
\trace(K) & \sum_{i=1}^n K_{ii}^2
\end{pmatrix}.
\]
Replacing $\delta_{\varepsilon}$ by
\[
\hat{\delta}_{\varepsilon}
=
2\left(\max\{\hat{\kappa}^{\mathrm{MoM}}_{\varepsilon,\mathrm{H}},0\}-1\right)
\hat{\sigma}^4_{\mathrm{MoM}}
\]
gives the homo-non-Gauss MoM covariance estimator used in the simulations. Finally, applying the delta method gives
\[
\widehat{\Var}(\hat{\gamma}_{\mathrm{MoM}}\mid \mtx{Z})
=
\hat{\vct{\ell}}^{\top}\widehat{\mtx{\Omega}}\,\hat{\vct{\ell}},
\]
with $\widehat{\mtx{\Omega}}$ chosen as the homo-Gauss, hetero-Gauss, or homo-non-Gauss moment covariance matrix.

\section{Additional Simulations}
\label{app:additional_simulation_figures}

This appendix collects supplementary design-distribution checks and diagnostic plots that support the simulation tables in Section \ref{sec:experiments}. The main text focuses on tabular comparisons of confidence intervals under heterogeneous and non-Gaussian noise.
The results below are not intended as separate benchmark studies; rather, they document that the same qualitative patterns persist across several non-Gaussian designs and that the normal approximations used in the main inference tables are empirically reasonable.

\subsection{Design-Distribution Checks}

Table \ref{tab:ci_gamma_design} reports the same confidence-interval metrics under three non-Gaussian design distributions: Rademacher, standardized $t_7$, and standardized genotype.
The setting is otherwise matched to the main heterogeneous-noise comparison, so the table should be read as a robustness check with respect to the covariate distribution rather than a new noise experiment.
Across the three designs, the hetero-Gauss calibration gives coverage much closer to the nominal level than the homo-Gauss calibration, while the Monte Carlo means remain close to the target SNR.

\begin{table}[H]
\centering
\small
\resizebox{\textwidth}{!}{%
\begin{tabular}{llcccc}
\hline
Design & Method & Mean & Bias & Coverage & Length \\
\hline
Rademacher & MLE (homo-Gauss) & 2.077 & 0.077 & 0.835 & 1.105 \\
Rademacher & MLE (hetero-Gauss) & 2.077 & 0.077 & 0.950 & 1.665 \\
Rademacher & MoM (homo-Gauss) & 2.071 & 0.071 & 0.890 & 1.771 \\
Rademacher & MoM (hetero-Gauss) & 2.071 & 0.071 & 0.950 & 2.186 \\
\hline
Standardized $t_7$ & MLE (homo-Gauss) & 2.066 & 0.066 & 0.795 & 1.099 \\
Standardized $t_7$ & MLE (hetero-Gauss) & 2.066 & 0.066 & 0.940 & 1.652 \\
Standardized $t_7$ & MoM (homo-Gauss) & 2.094 & 0.094 & 0.890 & 1.806 \\
Standardized $t_7$ & MoM (hetero-Gauss) & 2.094 & 0.094 & 0.950 & 2.227 \\
\hline
Standardized genotype & MLE (homo-Gauss) & 2.055 & 0.055 & 0.795 & 1.094 \\
Standardized genotype & MLE (hetero-Gauss) & 2.055 & 0.055 & 0.935 & 1.653 \\
Standardized genotype & MoM (homo-Gauss) & 2.116 & 0.116 & 0.920 & 1.823 \\
Standardized genotype & MoM (hetero-Gauss) & 2.116 & 0.116 & 0.970 & 2.246 \\
\hline
\end{tabular}
}
\caption{Supplementary comparison of $95\%$ confidence intervals for $\gamma_0$ under non-Gaussian designs. All rows use $n=2400$, $p=4000$, $\gamma_0=2$, $\sigma_0^2=0.5$, $g=0.1$, and strong independent geometric heterogeneity with $\kappa_0=30.77$. Each row uses $200$ Monte Carlo replications.}
\label{tab:ci_gamma_design}
\end{table}

\subsection{Consistency Diagnostics}

\begin{figure}[H]
\centering
\subfigure[$\hat{\gamma}$ of Simulation (i)] {
\includegraphics[width=0.3\columnwidth]{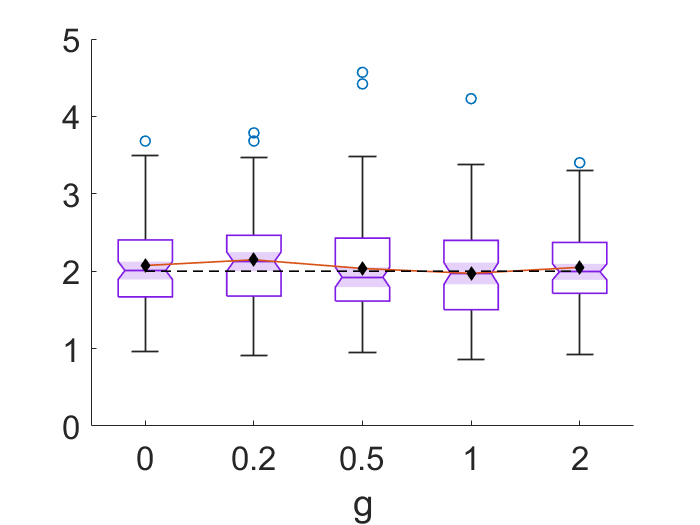}
}
\subfigure[$\hat{\gamma}$ of Simulation (ii)] {
\includegraphics[width=0.3\columnwidth]{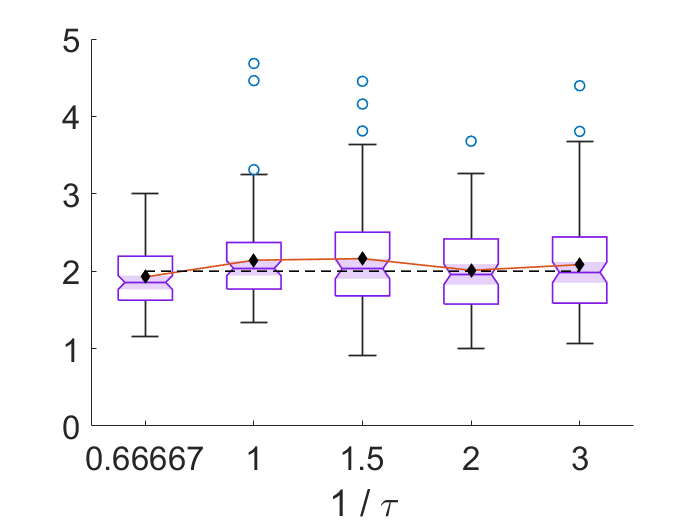}
}
\subfigure[$\hat{\gamma}$ of Simulation (iii)] {
\includegraphics[width=0.3\columnwidth]{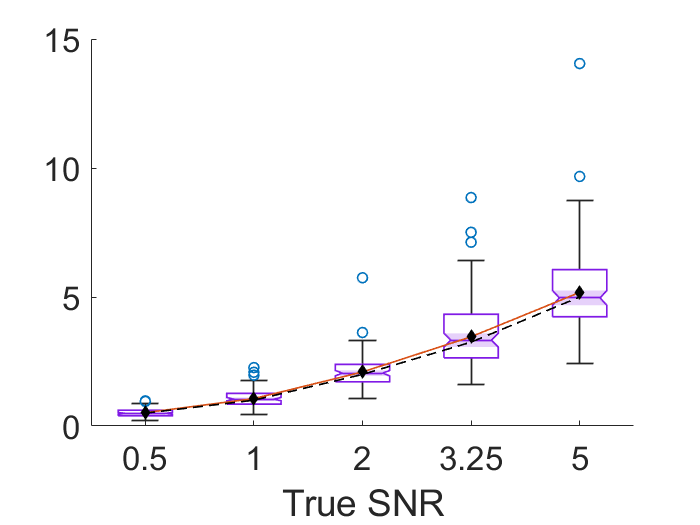}
}
\subfigure[$\hat{\sigma}^2$ of Simulation (i)] {
\includegraphics[width=0.3\columnwidth]{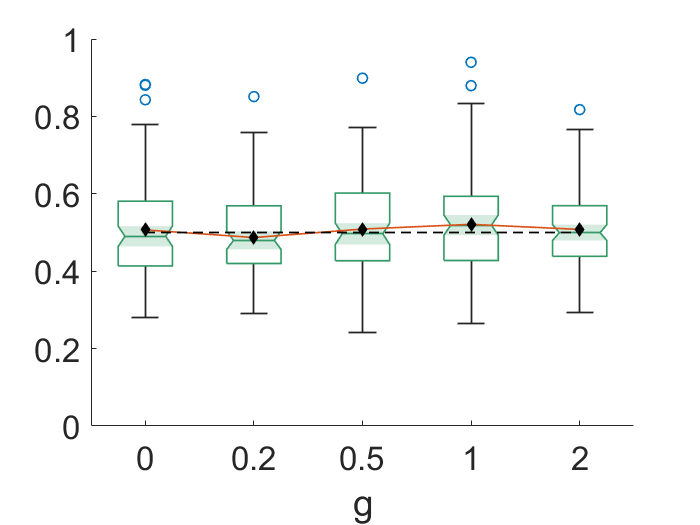}
}
\subfigure[$\hat{\sigma}^2$ of Simulation (ii)] {
\includegraphics[width=0.3\columnwidth]{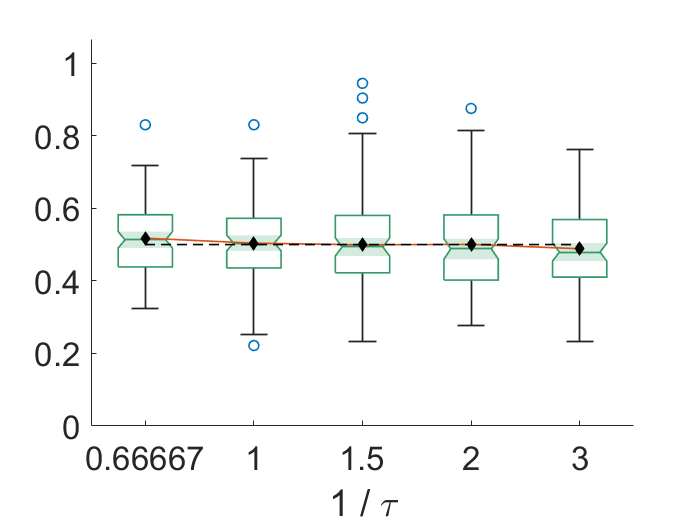}
}
\subfigure[$\hat{\sigma}^2$ of Simulation (iii)] {
\includegraphics[width=0.3\columnwidth]{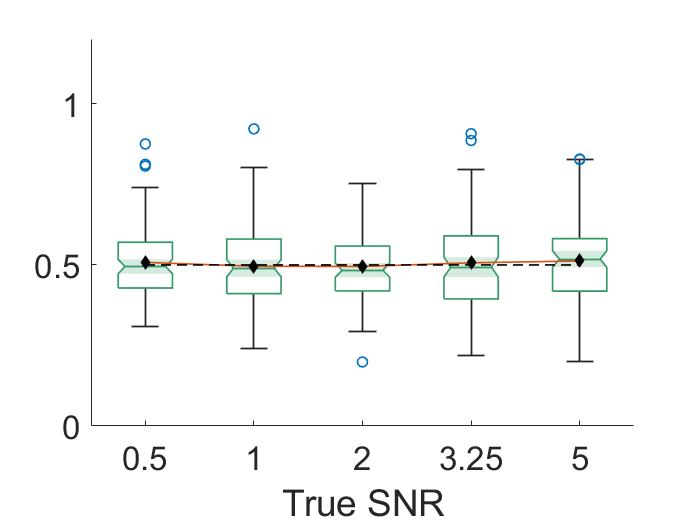}
}
\subfigure[$\hat{\kappa}$ of Simulation (i)] {
\includegraphics[width=0.3\columnwidth]{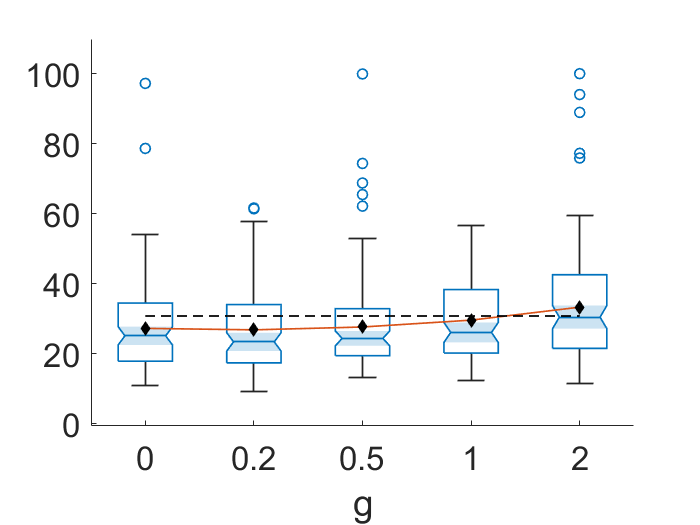}
}
\subfigure[$\hat{\kappa}$ of Simulation (ii)] {
\includegraphics[width=0.3\columnwidth]{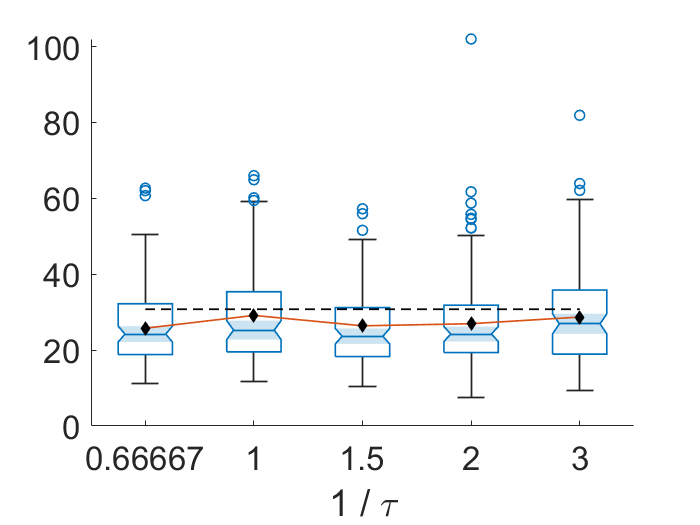}
}
\subfigure[$\hat{\kappa}$ of Simulation (iii)] {
\includegraphics[width=0.3\columnwidth]{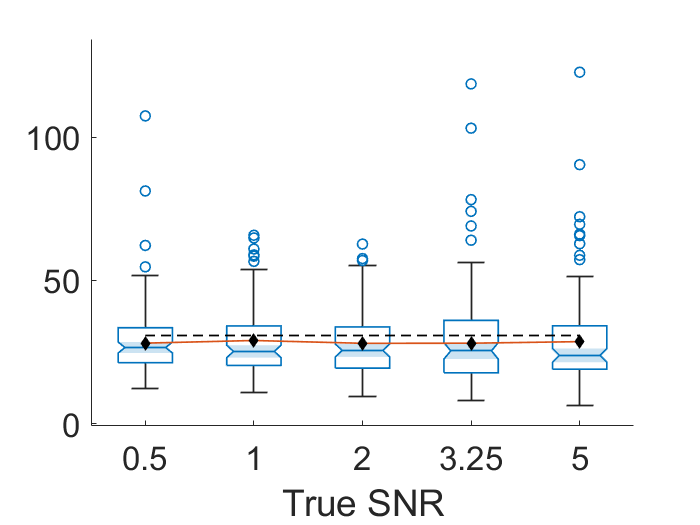}
}
\caption{Estimates of SNR, noise variance, and heterogeneity index for the sensitivity simulations under the standardized $t_5$ design. Each simulation uses $100$ independent Monte Carlo samples. The dashed lines mark the true values, and the black diamonds mark Monte Carlo averages.}
\label{fig:app_single_design_t5}
\end{figure}

Figure \ref{fig:app_single_design_t5} shows the behavior of the point estimators when the signal decay, aspect ratio, and SNR level are varied one at a time under a standardized $t_5$ design.
The panels are included mainly as diagnostics: the SNR and noise-variance estimates track their population targets over the parameter ranges considered, and the heterogeneity estimate is stable enough for the confidence-interval comparisons in the main text.

\subsection{Distribution Diagnostics}

\begin{figure}[H]
\centering
\subfigure[Rademacher design]
{
\includegraphics[width=1\columnwidth, height=5cm]{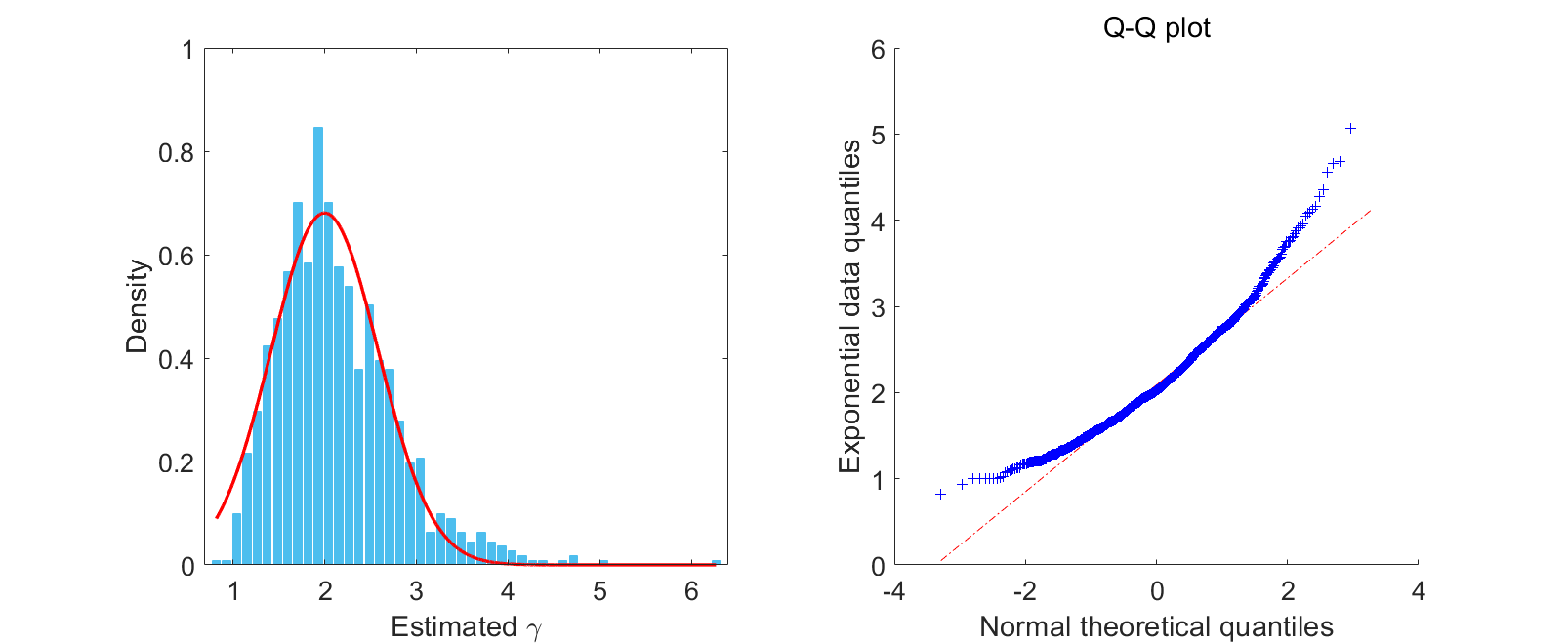}
}
\subfigure[Standardized $t_5$ design]
{
\includegraphics[width=1\columnwidth, height=5cm]{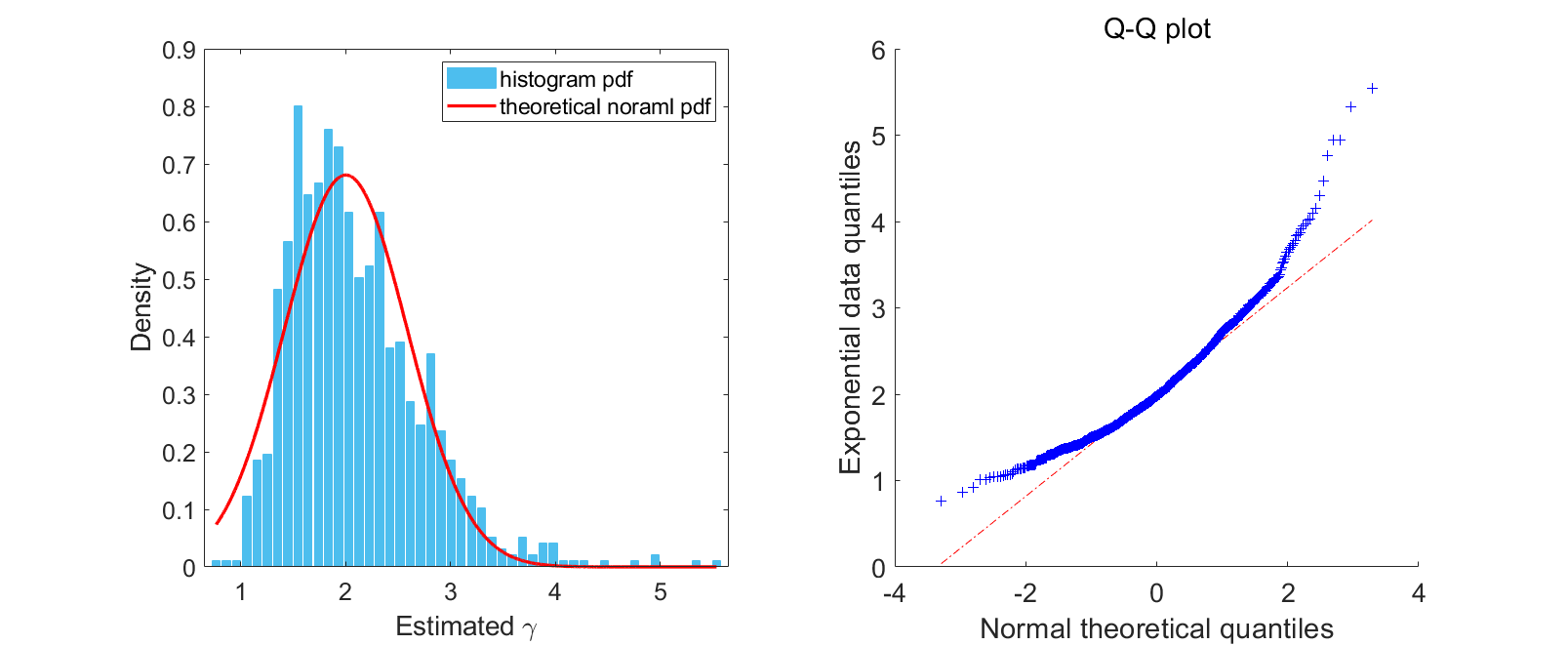}
}
\subfigure[Standardized genotype design]
{
\includegraphics[width=1\columnwidth, height=5cm]{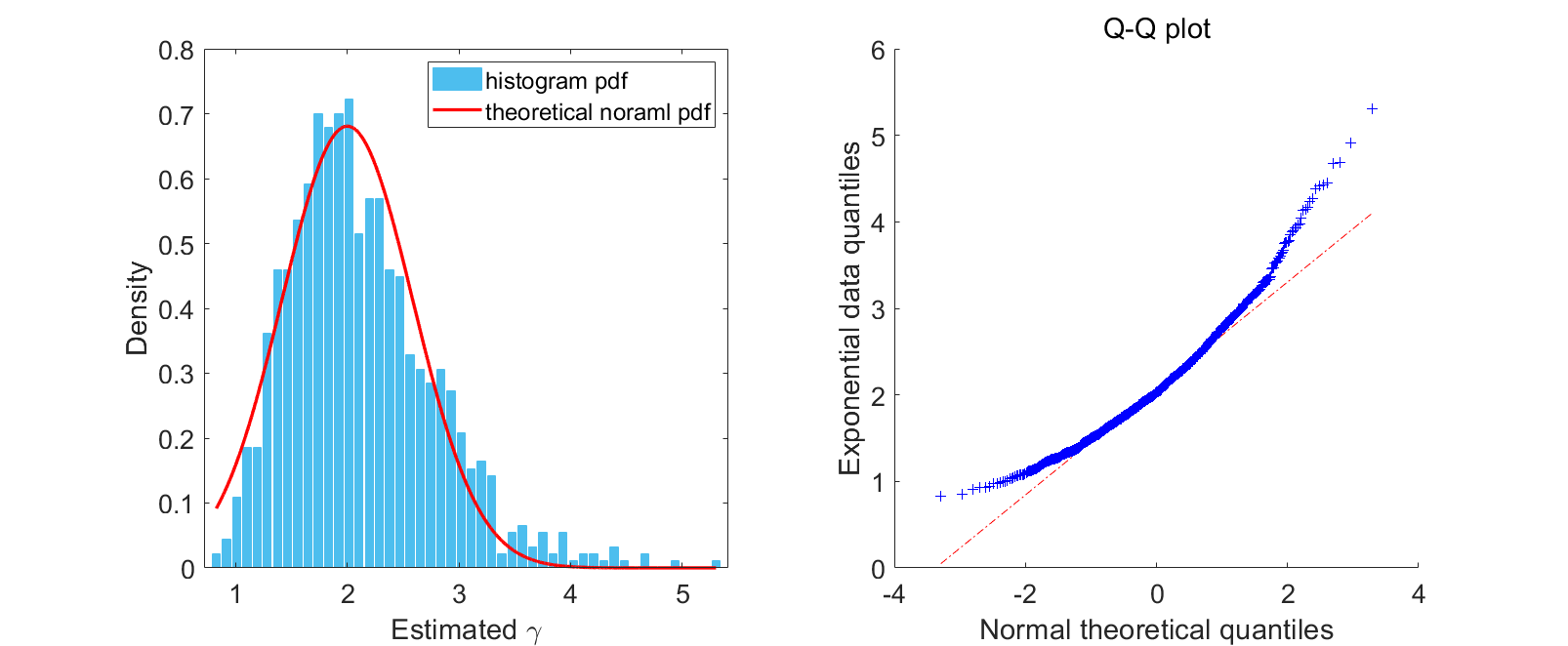}
}
\caption{Probability density of the estimated SNR $\hat{\gamma}$ and the normal Q--Q plot of the corresponding $\hat{\gamma}$ samples under heterogeneous and correlated noise. In the density plots, the purple curve uses the Monte Carlo mean and variance, while the red curve uses the theoretical asymptotic normal distribution.}
\label{fig:app_test_nor_Rad_hete_geo}
\end{figure}

Figure \ref{fig:app_test_nor_Rad_hete_geo} compares the empirical distribution of $\hat{\gamma}$ with the corresponding normal approximation under several non-Gaussian designs.
The density and Q--Q plots show the expected finite-sample deviations in the tails, but the central part of the distribution is well captured by the normal approximation, which is the feature most relevant for the reported confidence intervals.

\end{appendices}

\end{document}